\documentclass[12pt]{amsart}
\usepackage[all,dvips,arc,curve,color,frame]{xy}
\usepackage{graphics}
\usepackage{amsmath,amssymb}
\usepackage{color}
\usepackage{url}
\usepackage{epsfig}
\title{$(\{2,3\}, 6)$-spheres and their generalizations}

\author{Michel Deza}
\address{M. Deza, \'Ecole Normale Sup\'erieure, Paris}
\email{Michel.Deza@ens.fr}

\author{Mathieu Dutour Sikiri\'c}
\address{M. Dutour Sikiri\'c, Rudjer Boskovi\'c Institute, Bijenicka 54, 10000 Zagreb, Croatia}
\email{mdsikir@irb.hr}

\def\QuotS#1#2{\leavevmode\kern-.0em\raise.2ex\hbox{$#1$}\kern-.1em/\kern-.1em\lower.25ex\hbox{$#2$}}

\begin{document}
\newcommand{\KK}{\ensuremath{\mathbb{K}}}
\newcommand{\RR}{\ensuremath{\mathbb{R}}}
\newcommand{\NN}{\ensuremath{\mathbb{N}}}
\newcommand{\QQ}{\ensuremath{\mathbb{Q}}}
\newcommand{\ZZ}{\ensuremath{\mathbb{Z}}}

\newtheorem{theorem}{Theorem}
\newtheorem{proposition}[theorem]{Theorem}
\newtheorem{corollary}[theorem]{Corollary}
\newtheorem{lemma}[theorem]{Lemma}
\newtheorem{problem}[theorem]{Problem}
\newtheorem{conjecture}[theorem]{Conjecture}
\newtheorem{claim}[theorem]{Claim}
\newtheorem{remark}[theorem]{Remark}
\newtheorem{definition}[theorem]{Definition}

\begin{abstract}
We consider here $6$-regular plane graphs whose faces have size $1$,
$2$ or $3$.
In Section \ref{GenerationMethod} a practical enumeration method
is given that allowed us to enumerate them up to $53$ vertices.
Subsequently, in Section \ref{SymmetryGroupsSection} we enumerate
all possible symmetry groups of the spheres that showed up.
In Section \ref{GoldbergCoxeterConstruction} we introduce a new
Goldberg-Coxeter construction that takes a $6$-regular plane graph $G_0$,
two integers $k$ and $l$ and returns two $6$-regular plane graphs.

Then in the final section, we consider the notions of zigzags and central
circuits for the considered graphs.
We introduced the notions of tightness and weak tightness 
for them and we prove an upper bound on the number of zigzags and central
circuits of such tight graphs.
We also classify the tight and weakly tight graphs with simple zigzags or 
central circuits.
\end{abstract}

\maketitle

\section{Introduction}
By a {\em $(S, k)$-sphere} we call a plane $k$-regular graph such that
any face has size in $S$.

If $G$ is a $6$-regular plane graph, then by Euler formula it satisfies 
the equality:
\begin{equation*}
\sum_{k\geq 2} p_k(3-k)=6
\end{equation*}
with $p_k$ the number of $k$-gons, i.e. faces of size $k$.
So, if, moreover, $G$ has only $2$- and $3$-gonal faces, 
then it has exactly six $2$-gons.

Note that  a $(\{2,3\}, 6)$-sphere with $p_3$ $3$-gons
has $n=2+\frac{p_3}{2}$ vertices.
In \cite{book3} (Theorem 2.0.1)
we proved that for any $n\geq 2$ there exist a $(\{2,3\}, 6)$-sphere with 
$n$ vertices.
If $1$-gons are permitted, then $2p_1+p_2$ being $6$, all possible pairs
$(p_1,p_2)$, besides $(0,6)$, are $(1,4)$, $(2,2)$ and $(3,0)$. 

The only possible $(\{s-1,s\}, k)$-spheres have $(s,k)= (6,3)$ (well-known
geometrical {\em fullerenes}), $(4,4)$ (considered in \cite{octa,oct2,covcent,selfdual})
and $(3,6)$ (the object of this paper). $(\{2,3\}, 6)$-spheres are spherical 
analog of the $6$-regular partition $\{3^6\}$ of the Euclidean plane by 
regular triangles, with six $2$-gons playing role of ``defects'', 
disclinations needed to increase the curvature zero to the one of sphere.
The problem of existence of plane graphs with a fixed $p$-vector is
an active subject of research, see for example \cite{chinese4valent}.

In Section \ref{GenerationMethod} we expose a practical method for
generating $(\{1,2,3\}, 6)$-spheres. The main idea is to use a reduction
to $3$-regular graphs for which very efficient programs
exist \cite{Heidemeier}.
Then in Section \ref{SymmetryGroupsSection} we determine the
possible symmetry groups of $(\{1,2,3\},6)$-spheres with $i$ $1$-gons.
The methods are reasonably easy except for the $(\{1,3\},6)$-spheres
for which the symmetry groups are $C_3$, $C_{3h}$ or $C_{3v}$.

In Section \ref{GoldbergCoxeterConstruction} we introduce a new Goldberg-Coxeter
construction. It takes a $6$-regular sphere $G_0$, two integers $k,l$
and returns two $6$-regular spheres $G_1$, $G_2$
with $GC_{k,l}(G_0)=\{G_1, G_2\}$.
The construction satisfies a multiplicativity property based on the ring of
Eisenstein integers.
In the case $k=l=1$ we call the construction {\em oriented tripling} and
we have a more explicit description of it.
The Goldberg-Coxeter construction defined here generalizes the one
introduced in \cite{Gold,Cox71,octa} for $3$- or $4$-regular plane
graphs and allows to describe explicitly
all $(\{1,3\},6)$-spheres. It also allows to describe all
$(\{2,3\},6)$-spheres of symmetry $D_6$, $D_{6h}$, $T$, $T_h$, or $T_d$.

In a plane graph $G$, a {\em zigzag} is a circuit of edges such that any
two but no three consecutive edges are contained in the same face.
In an {\em Eulerian} (i.e. degree of any vertex is even) plane graph, a 
{\em central circuit} is a circuit
of edges such that any edge entering a vertex is followed by the edge 
opposite to the entering one.

A zigzag is called {\em simple} if no two edges
occur two times and a central circuit
is called {\em simple} if no two vertices occur two times.
Let $Z$ and $Z'$ be (possibly, $Z=Z'$) zigzags of a plane
graph $G$ and let an orientation be selected on them. An edge $e$ of
intersection $Z\cap Z'$ is called of {\em type I} or {\em type II},
if $Z$ and $Z'$ traverse $e$ in opposite or same direction, respectively
Let $C$ and $C'$ be (possibly, $C=C'$) central circuits of a
$6$-regular plane graph and let an orientation be selected on them.
A vertex $v$ of intersection $C\cap C'$ is called of {\em type I} or
{\em type II} if $C$ and $C'$ pass by $v$ with orientation shifted by
$60^{\circ}$, respectively, $120^{\circ}$.
We prove in Theorem \ref{TypeIIalways} that the intersection type
is always of type II.

We then introduce the notions of tighness and weak tightness for zigzags
and central circuits and we prove upper bound on the maximal number
of zigzags and central circuit. The results are summarized in Table
\ref{TableMaximalNrKnownResult} and Figure \ref{WeaklyIrreducibleC0123_zz}.
Then we determine completely the weakly tight spheres
with simple zigzags or central circuits.

\section{Generation method}\label{GenerationMethod}
In any $(\{2,3\}, 6)$-sphere, one can collapse its $2$-gons into 
simple edges.
By doing so one obtains a graph with vertices of degree at most $6$
and with faces of size $3$ only.
So, the dual will be a $3$-regular graph with faces of size at most $6$.

\begin{theorem}
With the exception of the following $(\{2,3\},6)$-spheres
\begin{center}
\begin{minipage}[b]{2.3cm}
\centering
\epsfig{height=16mm, file=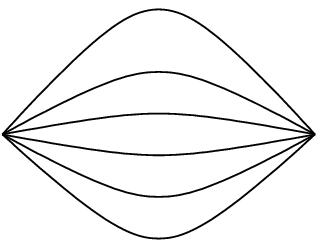}\par
$6\times K_2$: $D_{6h}$, $n=2$
\end{minipage}
\begin{minipage}[b]{2.3cm}
\centering
\epsfig{height=16mm, file=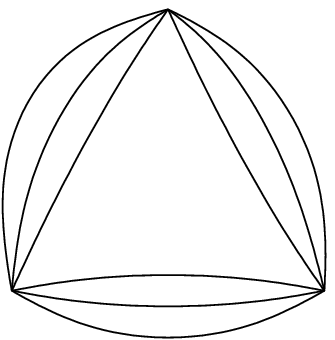}\par
$3\times K_3$: $D_{3h}$, $n=3$
\end{minipage}
\end{center}
any $(\{2,3\},6)$-sphere is obtained from a $(\{3,4,5,6\}, 3)$-sphere
by adding vertices of degree $2$ and taking the dual.

\end{theorem}
\proof Let $G$ be a $(\{2,3\},6)$-sphere and let $G^{*}$ be its dual.
Then, by removing from $G^{*}$ its vertices of degree $2$, one gets 
a $3$-regular graph $G_1$.
It can happen that $G_1$ has no vertices and is reduced to a simple
circular edge $e$. In this case, if one adds six vertices on $e$
and take the dual, 
one will get the first exceptional graph with $2$ vertices.
If $G_1$ has one face $F$ which is a $1$-gon, then we have to add $5$
vertices of degree $2$ on the edge $e$ of $F$.
Necessarily, any face adjacent to $F$ has to be a $1$-gon, but this
is, clearly, impossible.
If $F$ is a $2$-gon and $F$ is adjacent to at least one $2$-gon, 
then $G_1$ is reduced to a graph with two vertices and three edges.
The corresponding $(\{2,3\},6)$-sphere is the second exceptional graph.
Assume that $F$ is adjacent to $F_1$, $F_2$ with $F_i$ being a $a_i$-gon
and $a_i\geq 3$.
If one of $a_i$ is $3$, then the other is $6$ and this gives a $1$-gon.
Thus, the only possibility is $a_1=a_2=4$.
This implies that we have a graph with $4$ vertices, two $4$-gonal
and two $2$-gonal faces.
But consideration of all possibilities rules out this option.
So, $G_1$ is a $3$-regular plane graph with faces of size
within $\{3,4,5,6\}$. \qed

The method can be generalized (in Theorem \ref{Generation123graphs_C123})
to deal with graphs with $1$-gons.
Note that for most $(\{3,4,5,6\}, 3)$-spheres one cannot add
those vertices of degree $2$, in order to get the required spheres, because
whenever we add such a $2$-gon, we have two faces of size lower
than $6$ that are adjacent. Graphs admitting such adjacency are
relatively rare in the set of $(\{3,4,5,6\},3)$-spheres.
Some such graphs are the $(\{5,6\},3)$-spheres with the $5$-gons
organized in pairs, they are part of the class of {\em face-regular}
maps \cite{book3}.

The above theorem gives a method to enumerate $(\{2,3\}, 6)$-spheres.
First enumerate the $(\{3,4,5,6\}, 3)$-spheres using the program
{\tt CPF}, which is available from \cite{CaGe} and whose algorithm
has been described in \cite{Heidemeier}.
After such enumeration is done, the trick is to add the six vertices
of degree $2$ in all possibilities. This is relatively easy to do
and thus we have an efficient enumeration method. The numbers of graphs 
are shown in Table \ref{Enumeration123_6spheres_C0123} for $2\leq n\leq 41$.
We should point out that this algorithm while reasonable for our purpose
is very far from being optimal. A better method would be to adapt the 
algorithm from \cite{Heidemeier} although this is not easy to do.

\begin{table}
\caption{Number $N_i$ of $(\{1,2,3\}, 6)$-spheres with $n$ vertices and $(p_1,p_2)=(i,6-2i)$}
\label{Enumeration123_6spheres_C0123}
{\scriptsize
\begin{equation*}
\begin{array}{||c|c|c|c|c||c|c|c|c|c||c|c|c|c|c||}
\hline
\hline
n & N_0 & N_1 & N_2 & N_3 & n & N_0 & N_1 & N_2 & N_3 & n & N_0 & N_1 & N_2 & N_3\\
\hline
\hline
1 & 0 & 0 & 1 & 1 &19 & 69 & 36 & 13 & 1 &37 & 436 & 133 & 24 & 1\\
2 & 1 & 0 & 1 & 0 &20 & 100 & 34 & 28 & 0 &38 & 581 & 118 & 37 & 0\\
3 & 1 & 1 & 3 & 1 &21 & 86 & 46 & 19 & 1 &39 & 495 & 159 & 32 & 1\\
4 & 3 & 1 & 5 & 1 &22 & 133 & 33 & 23 & 0 &40 & 677 & 112 & 59 & 0\\
5 & 2 & 3 & 5 & 0 &23 & 112 & 62 & 16 & 0 &41 & 582 & 187 & 26 & 0\\
6 & 7 & 2 & 8 & 0 &24 & 165 & 44 & 37 & 0 &42 & 758 & 133 & 53 & 0\\
7 & 5 & 6 & 6 & 1 &25 & 144 & 57 & 20 & 1 &43 & 679 & 180 & 27 & 1\\
8 & 12 & 5 & 12 & 0 &26 & 205 & 54 & 27 & 0 &44 & 869 & 172 & 53 & 0\\
9 & 10 & 8 & 8 & 1 &27 & 176 & 75 & 22 & 1 &45 & 749 & 199 & 43 & 0\\
10 & 19 & 6 & 12 & 0 &28 & 251 & 61 & 36 & 1 &46 & 1000 & 149 & 44 & 0\\
11 & 16 & 14 & 9 & 0 &29 & 214 & 95 & 19 & 0 &47 & 868 & 250 & 30 & 0\\
12 & 29 & 11 & 17 & 1 &30 & 299 & 61 & 40 & 0 &48 & 1101 & 182 & 72 & 1\\
13 & 24 & 17 & 10 & 1 &31 & 265 & 96 & 20 & 1 &49 & 989 & 235 & 35 & 2\\
14 & 42 & 16 & 16 & 0 &32 & 360 & 89 & 43 & 0 &50 & 1259 & 194 & 57 & 0\\
15 & 35 & 23 & 15 & 0 &33 & 305 & 111 & 28 & 0 &51 & 1076 & 270 & 40 & 0\\
16 & 59 & 18 & 22 & 1 &34 & 429 & 80 & 33 & 0 &52 & 1410 & 210 & 61 & 1\\
17 & 48 & 33 & 12 & 0 &35 & 375 & 134 & 31 & 0 &53 & 1228 & 313 & 33 & 0\\
18 & 79 & 22 & 22 & 0 &36 & 488 & 105 & 50 & 1 & & & & & \\
\hline
\hline
\end{array}
\end{equation*}
}
\end{table}

\begin{figure}
\begin{center}
\begin{minipage}[b]{3.2cm}
\centering
\epsfig{height=20mm, file=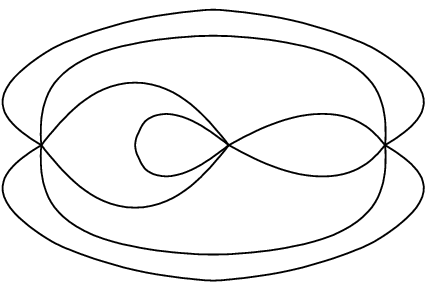}\par
$C_{s}$, $n=3$
\end{minipage}
\begin{minipage}[b]{3.2cm}
\centering
\epsfig{height=20mm, file=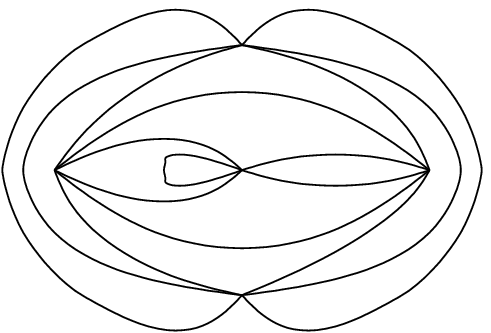}\par
$C_{s}$, $n=5$
\end{minipage}
\begin{minipage}[b]{3.2cm}
\centering
\epsfig{height=20mm, file=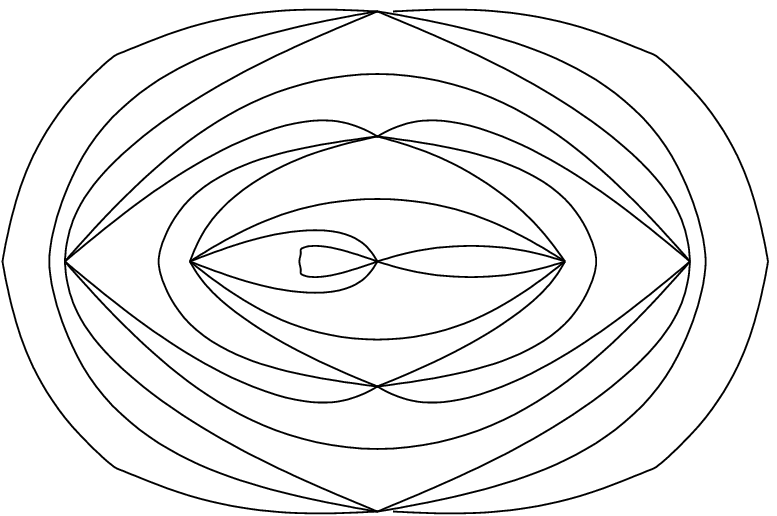}\par
$C_{s}$, $n=7$
\end{minipage}

\end{center}
\caption{First terms of an infinite sequence of $(\{1,2,3\},6)$-spheres $R_{2i+1}$ with $(p_1,p_2)=(1,4)$}
\label{FirstTermInfiniteSequence_C1}
\end{figure}

\begin{figure}
\begin{center}
\begin{minipage}[b]{3.2cm}
\centering
\epsfig{height=20mm, file=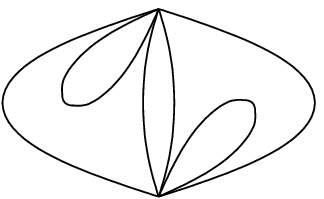}\par
$C_{2h}$, $n=2$
\end{minipage}
\begin{minipage}[b]{3.2cm}
\centering
\epsfig{height=20mm, file=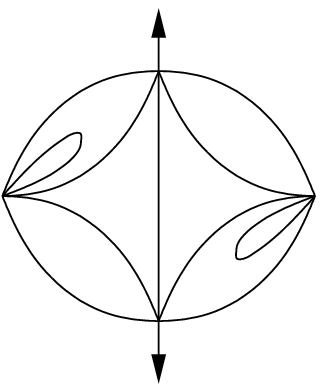}\par
$C_{2h}$, $n=4$
\end{minipage}
\begin{minipage}[b]{3.2cm}
\centering
\epsfig{height=20mm, file=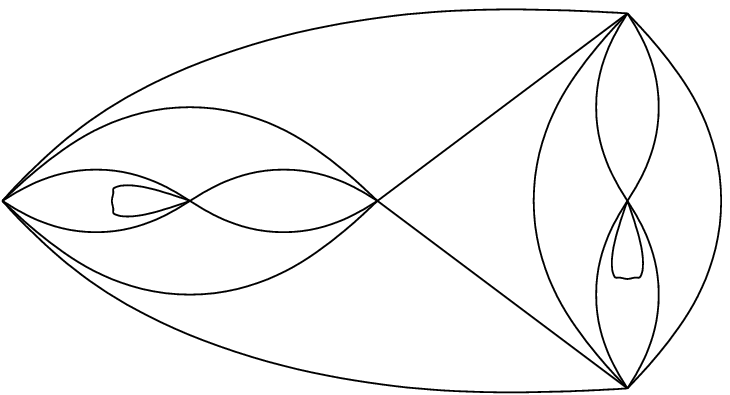}\par
$C_{2h}$, $n=6$
\end{minipage}

\end{center}
\caption{First terms of an infinite sequence of $(\{1,2,3\},6)$-spheres $S_{2i}$ with $(p_1,p_2)=(2,2)$}
\label{FirstTermInfiniteSequence_C2}
\end{figure}

\begin{figure}
\begin{center}
\begin{minipage}[b]{3.2cm}
\centering
\epsfig{height=15mm, file=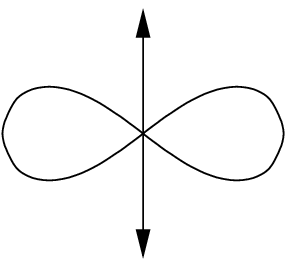}\par
$C_{2v}$, $n=1$
\end{minipage}
\begin{minipage}[b]{3.2cm}
\centering
\epsfig{height=15mm, file=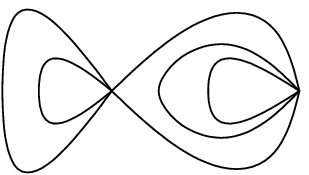}\par
$C_{2h}$, $n=2$
\end{minipage}
\begin{minipage}[b]{3.2cm}
\centering
\epsfig{height=15mm, file=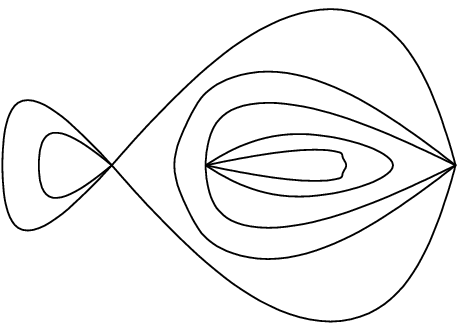}\par
$C_{2v}$, $n=3$
\end{minipage}

\end{center}
\caption{First terms of an infinite sequence of $(\{1,2,3\},6)$-spheres $A_{i}$ with $(p_1,p_2)=(2,2)$}
\label{FirstTermInfiniteSequence_C2_adj}
\end{figure}

\begin{figure}
\begin{center}
\begin{minipage}[b]{3.2cm}
\centering
\epsfig{height=20mm, file=InfSer_C2_0.eps}\par
$C_{2h}$, $n=2$
\end{minipage}
\begin{minipage}[b]{3.2cm}
\centering
\epsfig{height=20mm, file=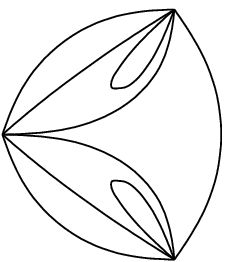}\par
$C_{s}$, $n=3$
\end{minipage}
\begin{minipage}[b]{3.2cm}
\centering
\epsfig{height=20mm, file=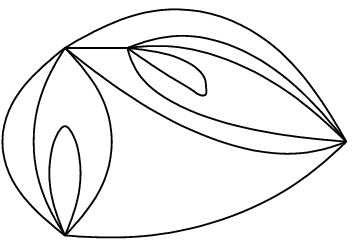}\par
$C_{s}$, $n=4$
\end{minipage}
\end{center}

\begin{center}
\begin{minipage}[b]{3.2cm}
\centering
\epsfig{height=15mm, file=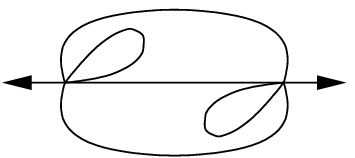}\par
$C_{2}$, $n=2$
\end{minipage}
\begin{minipage}[b]{3.2cm}
\centering
\epsfig{height=20mm, file=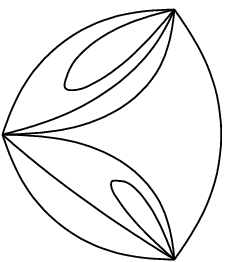}\par
$C_{2}$, $n=3$
\end{minipage}
\begin{minipage}[b]{3.2cm}
\centering
\epsfig{height=20mm, file=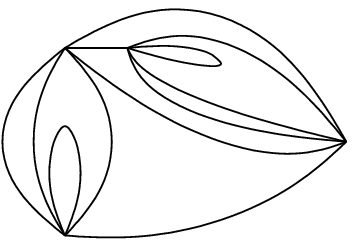}\par
$C_{2}$, $n=4$
\end{minipage}

\end{center}
\caption{First terms of two infinite sequences of $(\{1,2,3\},6)$-spheres $B_{i+1}$, $C_{i+1}$ with $(p_1,p_2)=(2,2)$}
\label{FirstTermInfiniteSequence_C2_132}
\end{figure}

\begin{theorem}\label{Generation123graphs_C123}
With the exception of the following graphs $T_1$, $T_2$
\begin{center}
\begin{minipage}[b]{5.2cm}
\centering
\epsfig{height=20mm, file=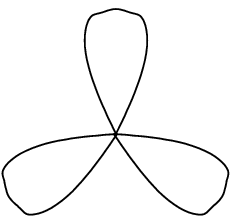}\par
Trifolium $T_1$: $C_{3v}$, $n=1$
\end{minipage}
\begin{minipage}[b]{3.2cm}
\centering
\epsfig{height=20mm, file=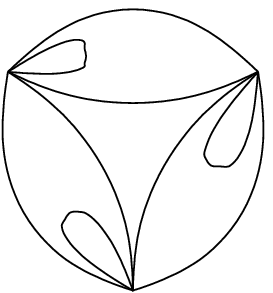}\par
$T_2$: $C_{3h}$, $n=3$
\end{minipage}
\end{center}
and the spheres of the infinite series depicted in Figures \ref{FirstTermInfiniteSequence_C1}, \ref{FirstTermInfiniteSequence_C2}, \ref{FirstTermInfiniteSequence_C2_adj} and \ref{FirstTermInfiniteSequence_C2_132},
any $(\{1,2,3\},6)$-sphere with at least one $1$-gon is obtained from a 
$(\{3,4,5,6\}, 3)$-sphere by taking the dual and then splitting some edges according to following two schemes:
\begin{center}
\begin{minipage}[b]{3.8cm}
\centering
\epsfig{width=30mm, file=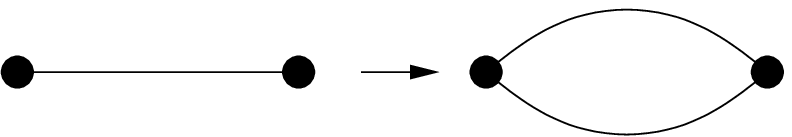}\par
\end{minipage}
\begin{minipage}[b]{3.8cm}
\centering
\epsfig{width=30mm, file=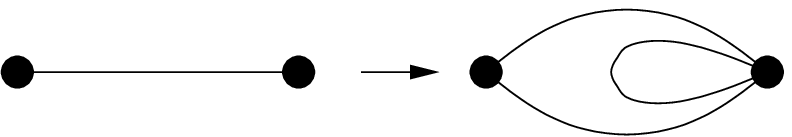}\par
\end{minipage}

\end{center}

\end{theorem}
\proof Let us take a $(\{1,2,3\}, 6)$-sphere $G$ with at least one $1$-gon $F$ in its face-set.
Clearly, $F$ cannot be adjacent to another $1$-gon.
If $F$ is adjacent to a $2$-gon, then simple considerations yield that $G$ belongs to the infinite series of Figure \ref{FirstTermInfiniteSequence_C2_adj}.
So, we can assume in the following that all $1$-gons, say $F_1$, \dots ,$F_s$ are adjacent to $3$-gons $G_1$,\dots, $G_s$.
If one of the $G_i$ is adjacent to two $2$-gons, then we get the
sphere $B_2$ ($C_{2h}$, $n=2$) depicted in
Figure \ref{FirstTermInfiniteSequence_C2_132}.
If one of the $G_i$ is adjacent to exactly one $1$-gon, then we get the following partial diagram:
\begin{center}
\begin{minipage}[b]{3.8cm}
\centering
\epsfig{width=25mm, file=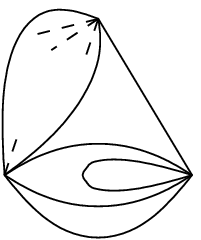}\par
\end{minipage}
\end{center}
Clearly, such diagram extends to one of the graphs of the 
infinite series depicted in
Figure \ref{FirstTermInfiniteSequence_C2_132}.

So, we can now assume that the $G_i$ are adjacent to $3$-gons only.
If one of the $3$-gons adjacent to a $G_i$ turns out to be another $G_j$, then
we get the map $C_2$ from Figure \ref{FirstTermInfiniteSequence_C2_132}.
So, we assume further that those $3$-gons are not of the type $G_i$.

The faces $G_i$ contains two vertices $v_i$, $v'_i$ with $v_i$ being contained
in $F_i$.
If $v_i=v'_i$, then we get the exceptional sphere Trifolium.
So, we assume further that $v_i\not=v'_i$.
If $v_i=v'_j$ for $i\not=j$, then some easy considerations gives the sphere
$T_2$ as the only possibility.
So, let us assume now that the vertex $v_i$ is contained in a $2$-gon.
Then we have the following local configuration:
\begin{center}
\begin{minipage}[b]{3.8cm}
\centering
\epsfig{width=30mm, file=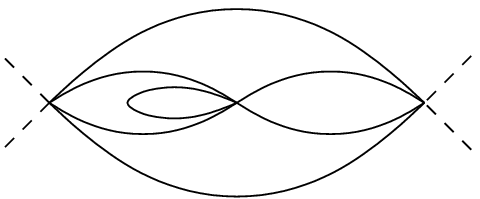}\par
\end{minipage}
\end{center}
From that point, after enumeration of all possibilities we get the infinite series of Figures \ref{FirstTermInfiniteSequence_C1} and \ref{FirstTermInfiniteSequence_C2}.
So, now we have that all vertices $v_i$ are contained in four $3$-gons.
This implies that $G$ is obtained from a $(\{3,4,5,6\}, 3)$-sphere by
taking the dual and then splitting some edges according to mentioned
above schemes. \qed

Obviously, the above theorem gives us a method to enumerate the
$(\{1,2,3\},6)$-spheres.
The enumeration results are shown in Table \ref{Enumeration123_6spheres_C0123}.
Like for $(\{2,3\},6)$-spheres, it would be interesting to have a faster
enumeration method.

\section{Symmetry groups}\label{SymmetryGroupsSection}
We now give the possible groups of the considered spheres.
Note that we are using the terminology of points groups
in chemistry as explained, for example, in \cite{pointgroup}.

\begin{figure}
\begin{center}
\begin{minipage}[b]{3.0cm}
\centering
\epsfig{height=20mm, file=Bundle6.eps}\par
$D_{6h}$, $n=2$
\end{minipage}
\begin{minipage}[b]{3.0cm}
\centering
\epsfig{height=20mm, file=Example23_6val_1.eps}\par
$D_{3h}$, $n=3$
\end{minipage}
\begin{minipage}[b]{3.0cm}
\centering
\epsfig{height=20mm, file=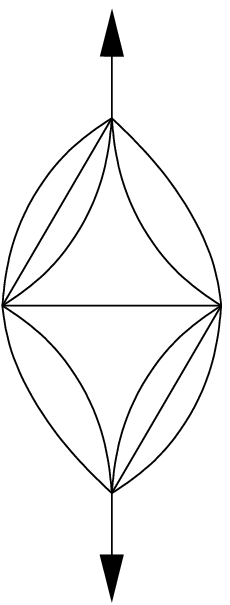}\par
$D_{2}$, $n=4$
\end{minipage}
\begin{minipage}[b]{3.0cm}
\centering
\epsfig{height=20mm, file=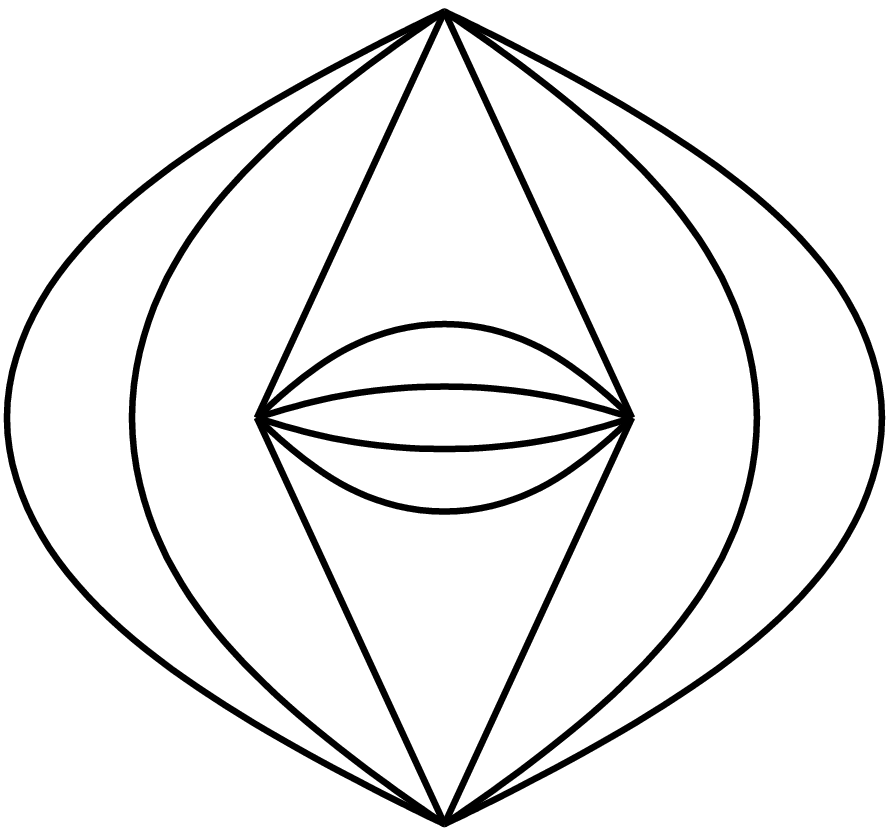}\par
$D_{2d}$, $n=4$
\end{minipage}
\begin{minipage}[b]{3.0cm}
\centering
\epsfig{height=20mm, file=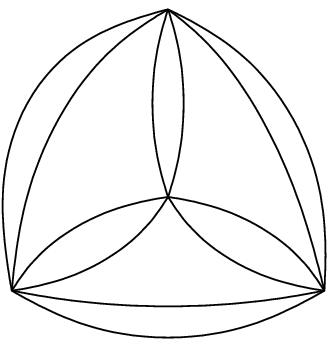}\par
$T_d$, $n=4$
\end{minipage}
\begin{minipage}[b]{3.0cm}
\centering
\epsfig{height=20mm, file=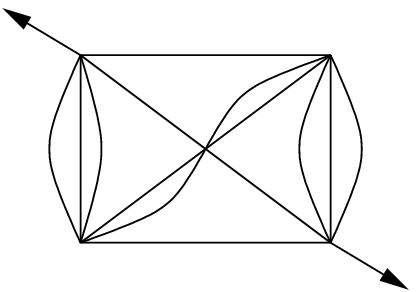}\par
$C_2$, $n=5$
\end{minipage}
\begin{minipage}[b]{3.0cm}
\centering
\epsfig{height=20mm, file=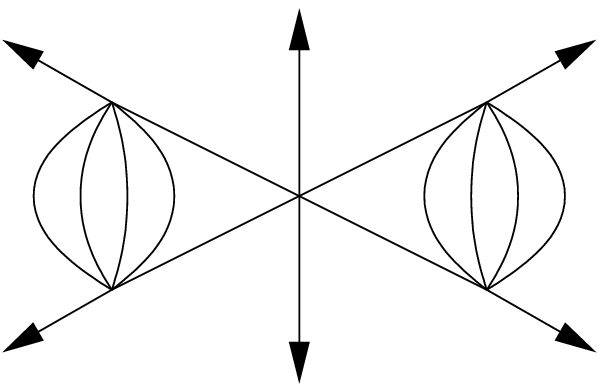}\par
$D_{2h}$, $n=6$
\end{minipage}
\begin{minipage}[b]{3.0cm}
\centering
\epsfig{height=20mm, file=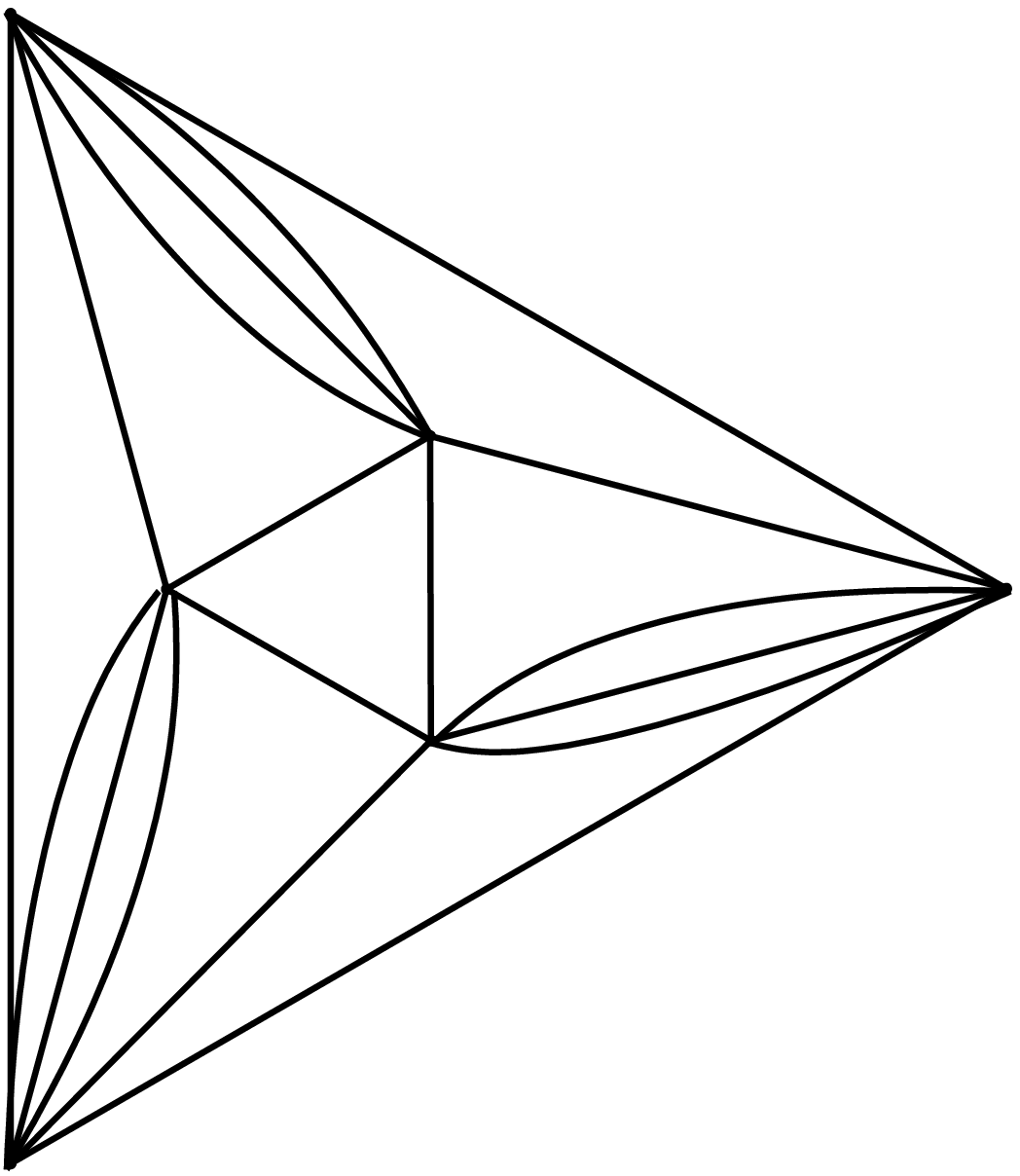}\par
$D_3$, $n=6$
\end{minipage}
\begin{minipage}[b]{3.0cm}
\centering
\epsfig{height=20mm, file=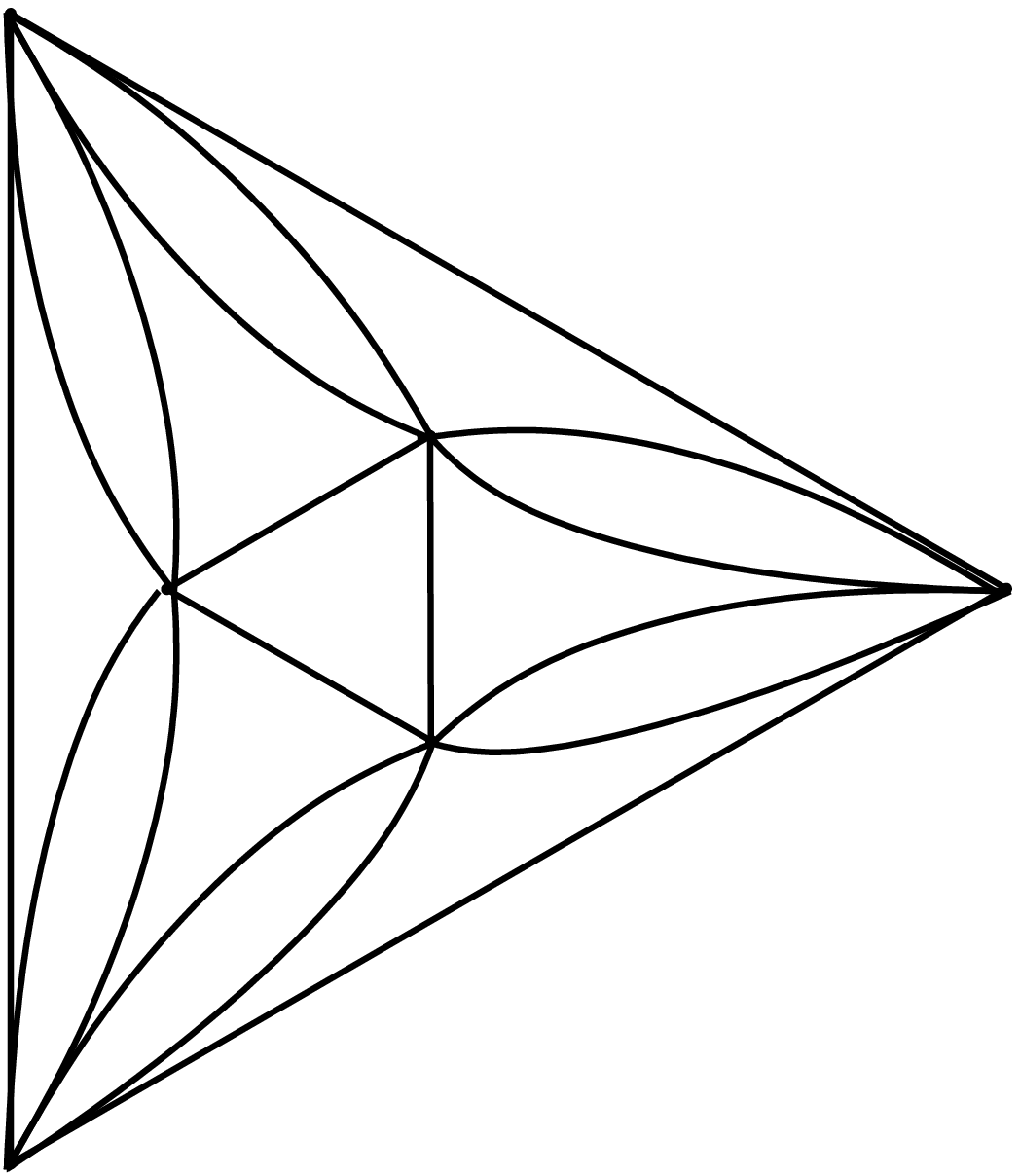}\par
$D_{3d}$, $n=6$
\end{minipage}
\begin{minipage}[b]{3.0cm}
\centering
\epsfig{height=20mm, file=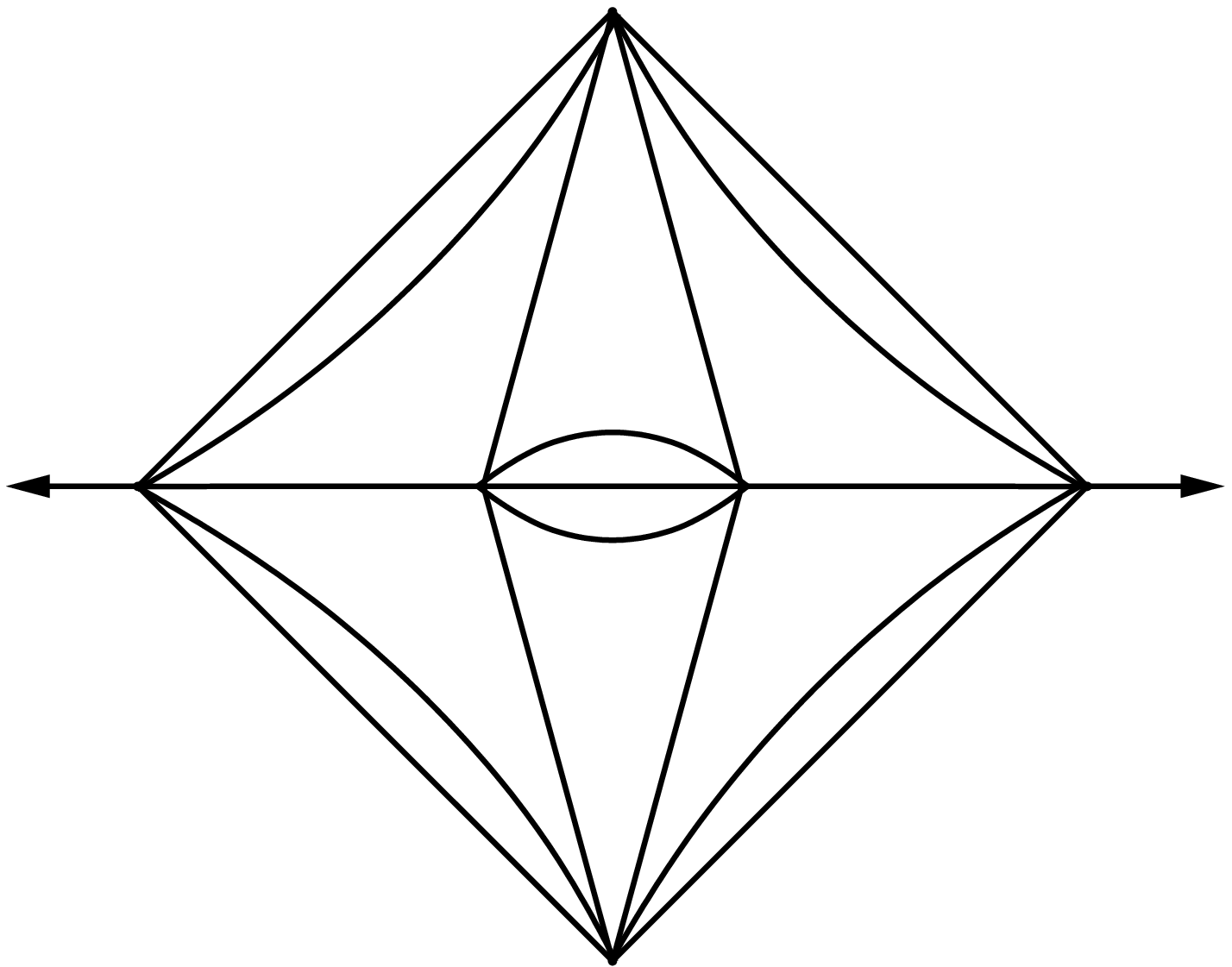}\par
$C_{2v}$, $n=6$
\end{minipage}
\begin{minipage}[b]{3.0cm}
\centering
\epsfig{height=20mm, file=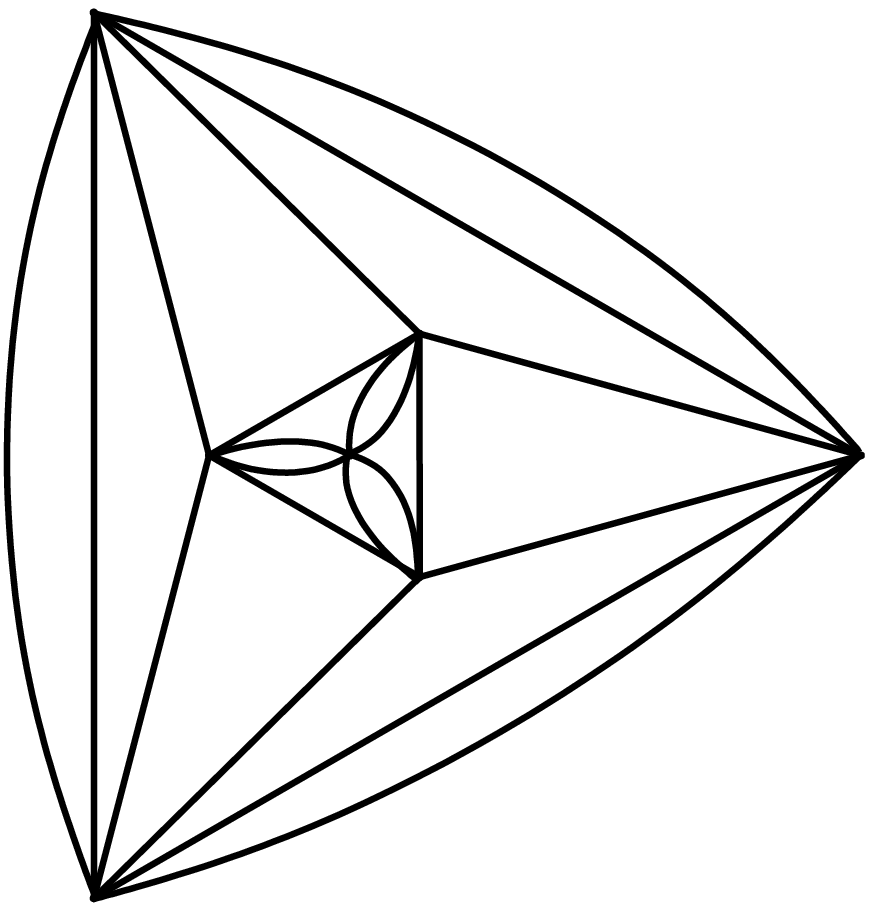}\par
$C_{3v}$, $n=7$
\end{minipage}
\begin{minipage}[b]{3.0cm}
\centering
\epsfig{height=20mm, file=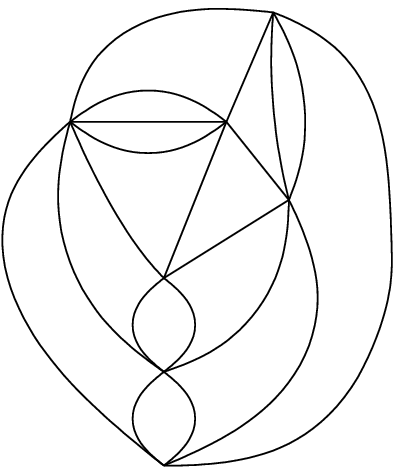}\par
$C_{1}$, $n=8$
\end{minipage}
\begin{minipage}[b]{3.0cm}
\centering
\epsfig{height=20mm, file=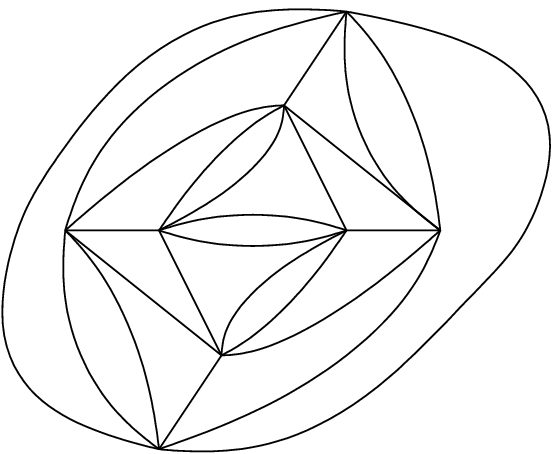}\par
$S_{4}$, $n=8$
\end{minipage}
\begin{minipage}[b]{3.0cm}
\centering
\epsfig{height=20mm, file=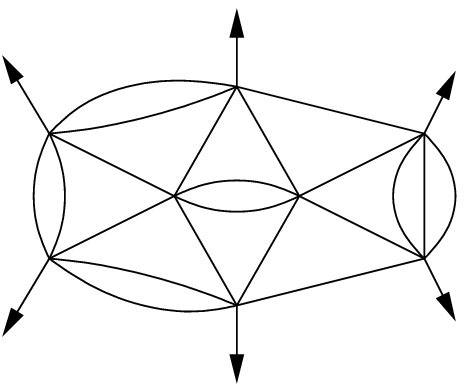}\par
$C_{s}$, $n=9$
\end{minipage}
\begin{minipage}[b]{3.0cm}
\centering
\epsfig{height=20mm, file=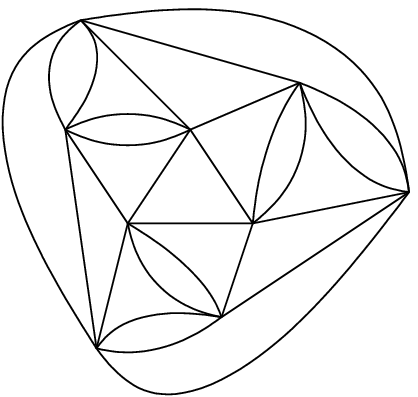}\par
$C_{3h}$, $n=9$
\end{minipage}
\begin{minipage}[b]{3.0cm}
\centering
\epsfig{height=20mm, file=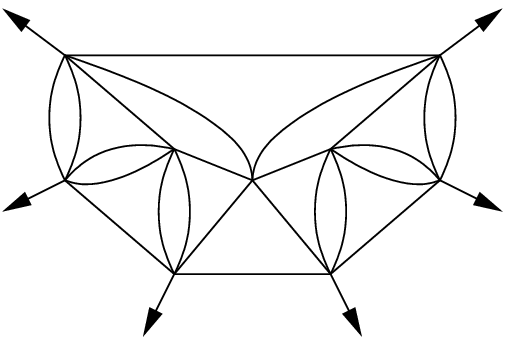}\par
$C_{2h}$, $n=10$
\end{minipage}
\begin{minipage}[b]{3.0cm}
\centering
\epsfig{height=20mm, file=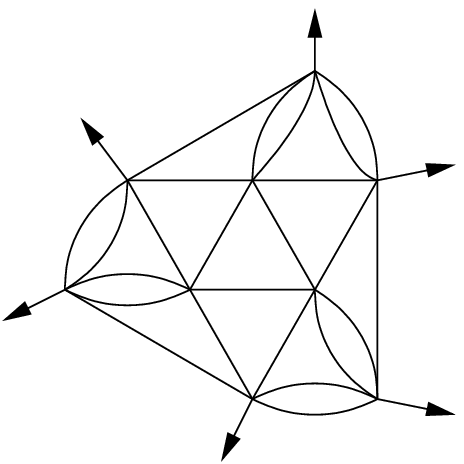}\par
$C_3$, $n=10$
\end{minipage}
\begin{minipage}[b]{3.0cm}
\centering
\epsfig{height=20mm, file=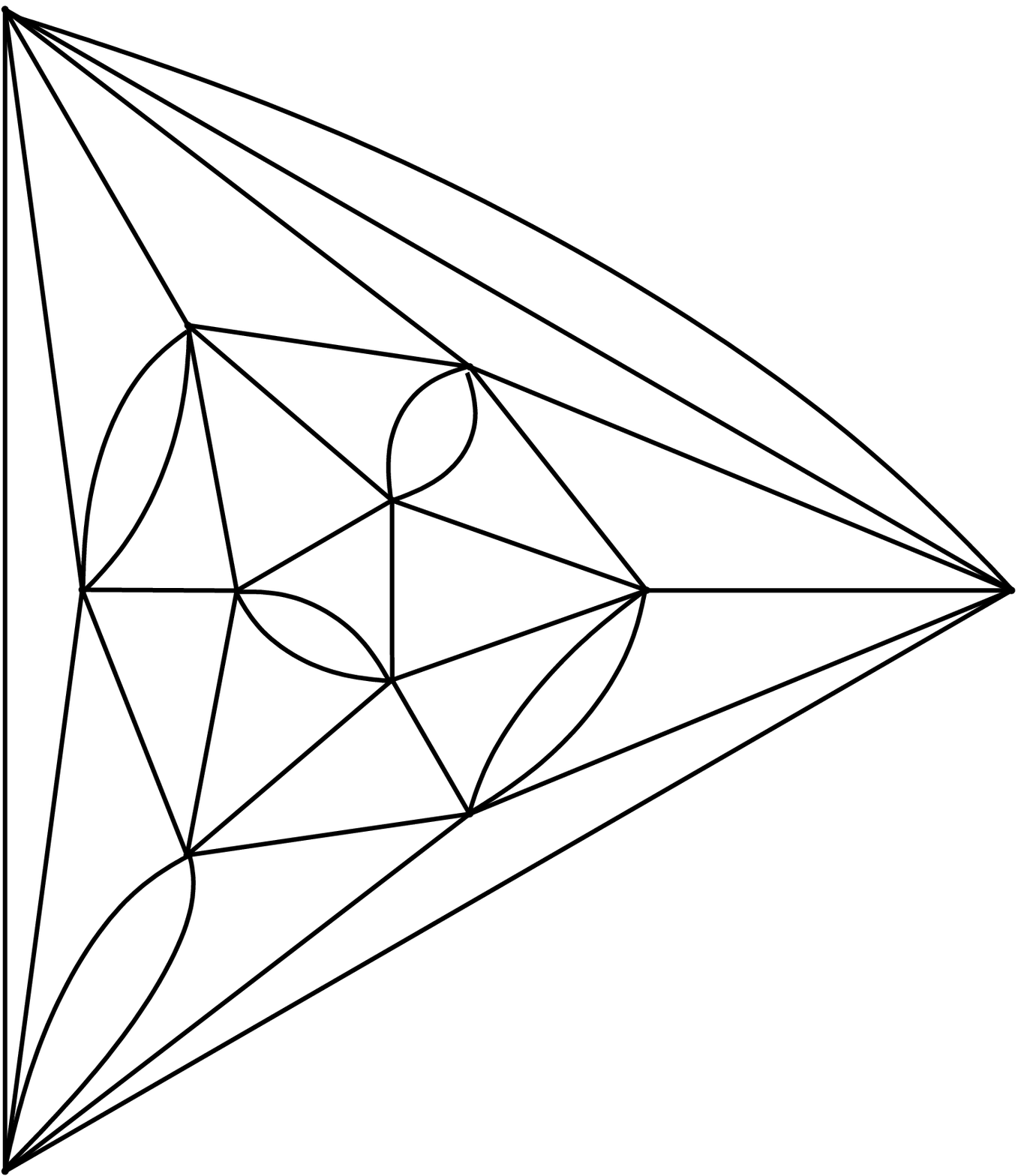}\par
$T_h$, $n=12$
\end{minipage}
\begin{minipage}[b]{3.0cm}
\centering
\epsfig{height=20mm, file=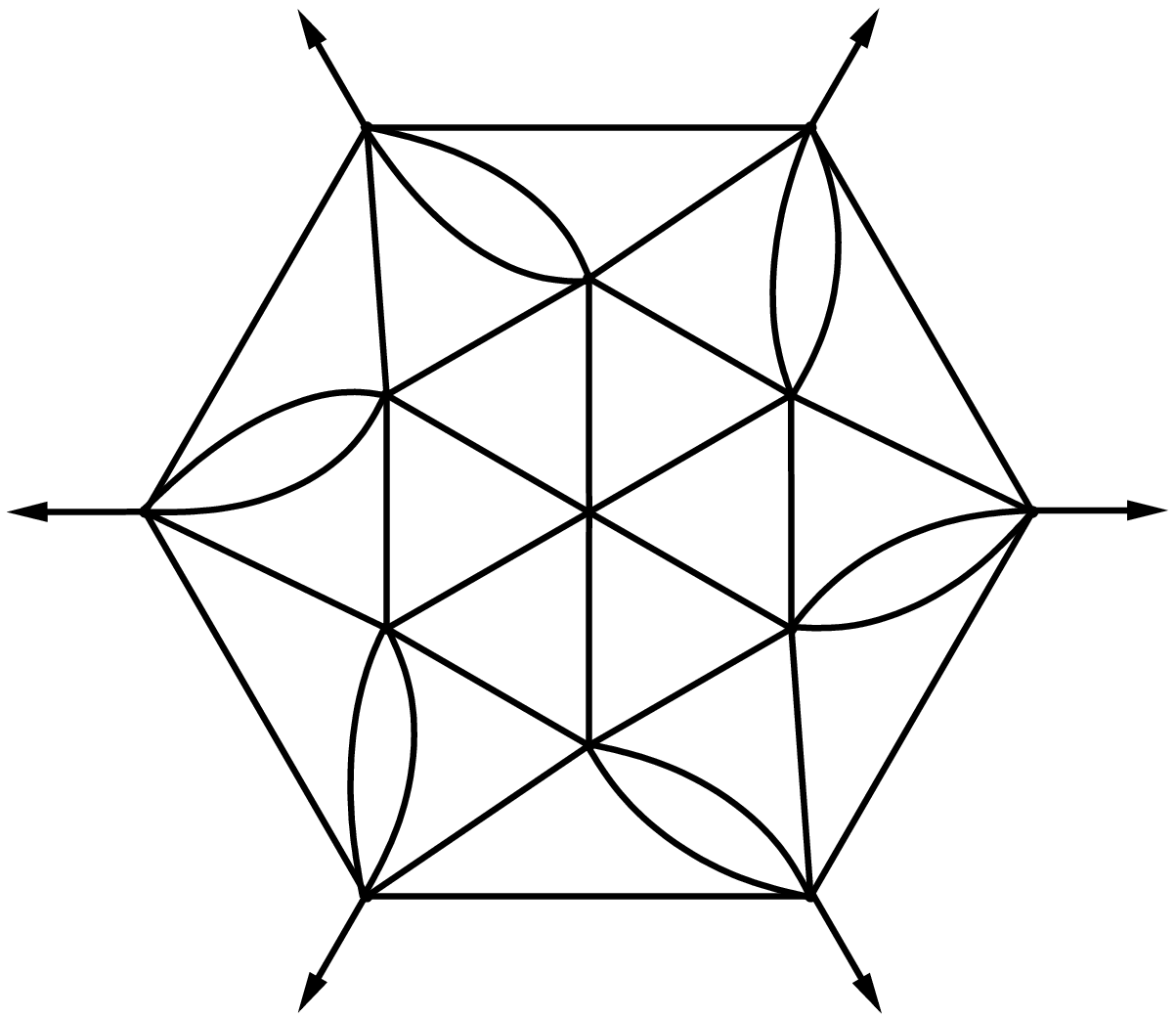}\par
$D_6$, $n=14$
\end{minipage}
\begin{minipage}[b]{3.0cm}
\centering
\epsfig{height=20mm, file=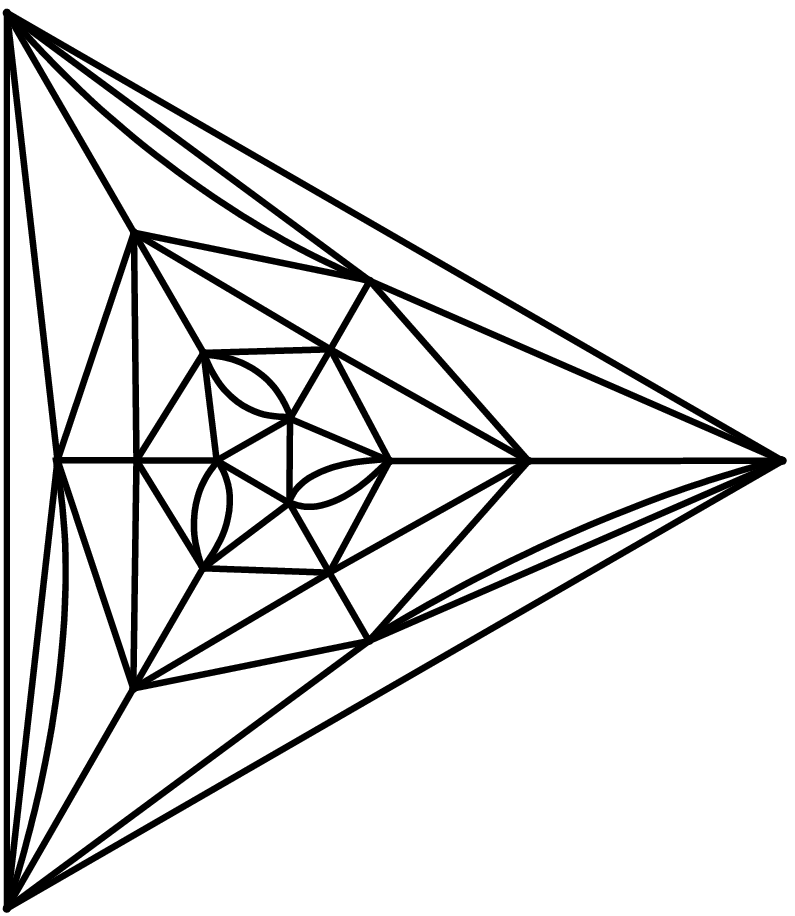}\par
$S_6$, $n=18$
\end{minipage}
\begin{minipage}[b]{3.0cm}
\centering
\epsfig{height=20mm, file=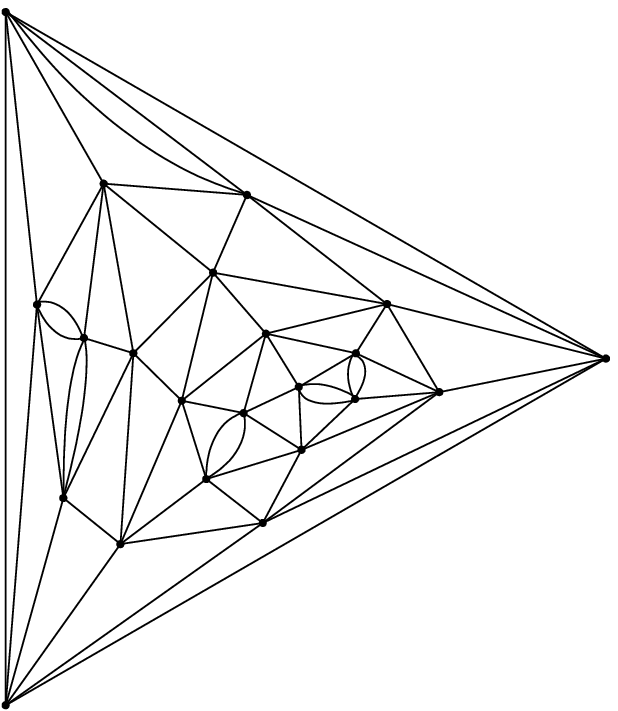}\par
$C_{i}$, $n=22$
\end{minipage}
\begin{minipage}[b]{3.0cm}
\centering
\epsfig{height=20mm, file=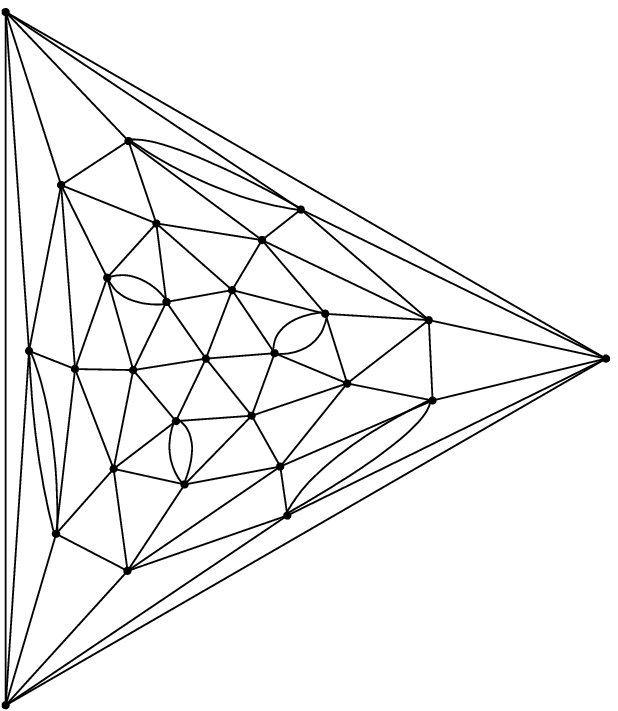}\par
$T$, $n=28$
\end{minipage}

\end{center}
\caption{Minimal representatives for each possible symmetry group of a $(\{2,3\}, 6)$-sphere}
\label{MinimalRepresentative}
\end{figure}

\begin{theorem}
The possible symmetry group of a $(\{2,3\}, 6)$-sphere are
$C_1$, $C_2$, $C_{2h}$, $C_{2v}$, $C_3$, $C_{3h}$, $C_{3v}$, $C_i$,
$C_s$, $D_2$, $D_{2d}$, $D_{2h}$, $D_3$, $D_{3d}$, $D_{3h}$, $D_6$,
$D_{6h}$, $S_4$, $S_6$, $T$, $T_h$ and $T_d$.
The minimal possible representatives are given in
Figure \ref{MinimalRepresentative}.
\end{theorem}
\proof The method is to consider the possible axes of symmetry;
they are passing through faces, edges or vertices. As a consequence,
the possibilities for a $k$-fold axis of symmetry are $2$, $3$ or $6$.
The only groups that could occur, besides $22$ given in the theorem,
are $D_{6d}$, $C_6$, $C_{6h}$ or $C_{6v}$.

If a $6$-fold axis occurs, then it necessarily passes through two vertices,
say, $v_1$ and $v_2$.
Around this vertex one can add successive rings of triangles as in the
classical structure of the triangular lattice.
At some point one gets a $2$-gon and thus, by the $6$-fold symmetry,
six $2$-gons.
Then, one can continue the structure uniquely and the structure is
defined uniquely.
This completion is the same as the one around $v$ and it implies
the existence of a mapping that inverts $v$ with the 
transformation inverting $v_1$ and $v_2$ and the group are $D_6$, $D_{6h}$
or $D_{6d}$.
The group $D_{6d}$ is ruled out
because $2$-fold axis passes through the $2$-gons. \qed

\begin{figure}
\begin{center}
\begin{minipage}[b]{3.2cm}
\centering
\epsfig{height=20mm, file=Class_C1_CsB.eps}\par
$C_{s}$, $n=3$
\end{minipage}
\begin{minipage}[b]{3.2cm}
\centering
\epsfig{height=20mm, file=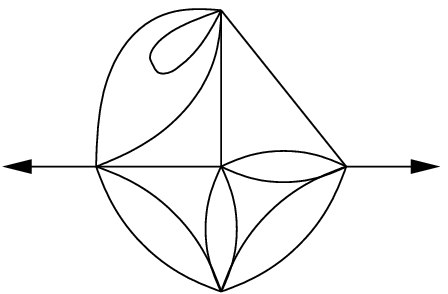}\par
$C_{1}$, $n=5$
\end{minipage}

\end{center}
\caption{Minimal representatives for each possible symmetry group of a $(\{1,2,3\}, 6)$-sphere with $(p_1,p_2)=(1,4)$}
\label{MinimalRepresentative_C1}
\end{figure}

\begin{figure}
\begin{center}
\begin{minipage}[b]{3.2cm}
\centering
\epsfig{height=15mm, file=Class_C2_C2vB.eps}\par
$C_{2v}$, $n=1$
\end{minipage}
\begin{minipage}[b]{3.2cm}
\centering
\epsfig{height=15mm, file=InfSer_C2_0.eps}\par
$C_{2h}$, $n=2$
\end{minipage}
\begin{minipage}[b]{3.2cm}
\centering
\epsfig{height=15mm, file=Class_C2_C2B.eps}\par
$C_{2}$, $n=2$
\end{minipage}
\begin{minipage}[b]{3.2cm}
\centering
\epsfig{height=20mm, file=Class_C2_CsB.eps}\par
$C_{s}$, $n=3$
\end{minipage}
\begin{minipage}[b]{3.2cm}
\centering
\epsfig{height=20mm, file=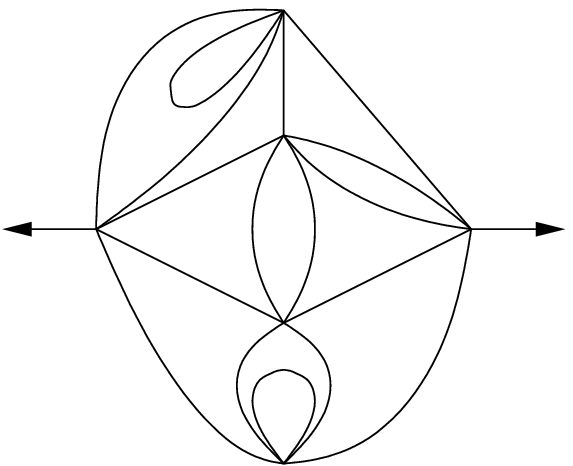}\par
$C_{1}$, $n=6$
\end{minipage}
\begin{minipage}[b]{3.2cm}
\centering
\epsfig{height=20mm, file=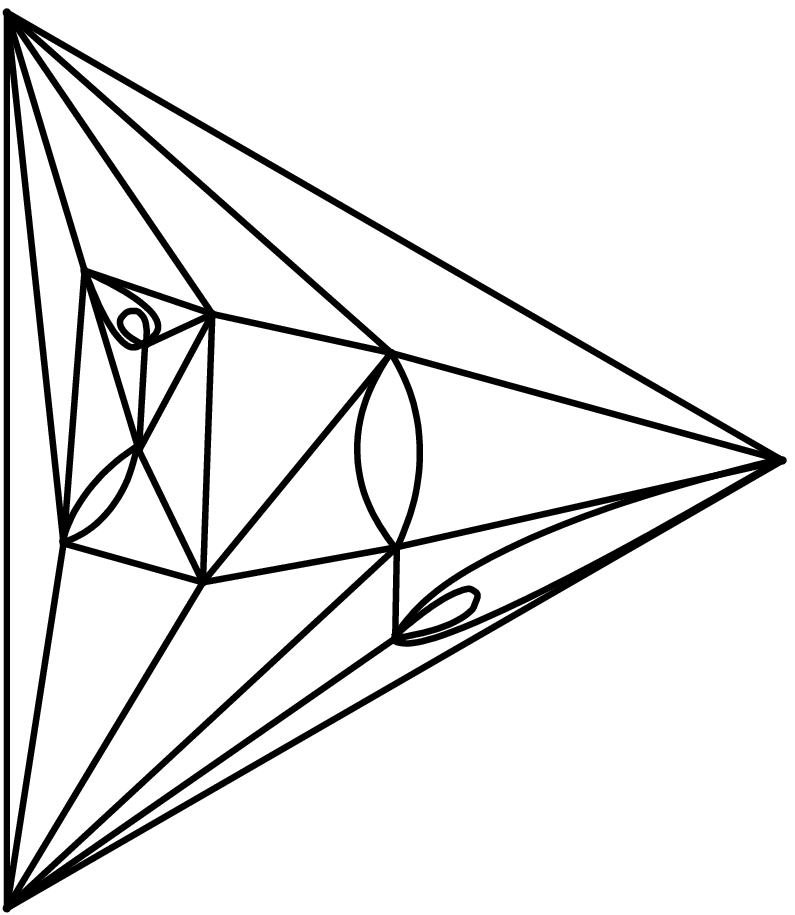}\par
$C_{i}$, $n=12$
\end{minipage}

\end{center}
\caption{Minimal representatives for each possible symmetry group of a $(\{1,2,3\}, 6)$-sphere with $(p_1,p_2)=(2,2)$}
\label{MinimalRepresentative_C2}
\end{figure}

\begin{figure}
\begin{center}
\begin{minipage}[b]{3.2cm}
\centering
\epsfig{height=20mm, file=Class_C3_C3vB.eps}\par
$C_{3v}$, $n=1$
\end{minipage}
\begin{minipage}[b]{3.2cm}
\centering
\epsfig{height=20mm, file=Singular13_3vert.eps}\par
$C_{3h}$, $n=3$
\end{minipage}
\begin{minipage}[b]{3.2cm}
\centering
\epsfig{height=20mm, file=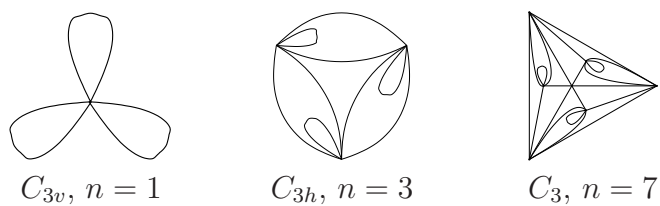}\par
$C_3$, $n=7$
\end{minipage}

\end{center}
\caption{Minimal representatives for each possible symmetry group of a $(\{1,3\}, 6)$-sphere}
\label{MinimalRepresentative_C3}
\end{figure}

\begin{theorem}\label{GroupSpheres_C1_C2_C3}
The possible symmetry group of a $(\{1,2,3\}, 6)$-sphere with $p_1>0$ are

(i) $C_1$ or $C_s$ if $p_1=1$.

(ii) $C_{1}$, $C_{2}$, $C_{i}$, $C_{s}$, $C_{2v}$ or $C_{2h}$
if $p_1=2$.

(iii) $C_3$, $C_{3v}$ or $C_{3h}$ if $p_1=3$.

The minimal possible representative are given in
Figures \ref{MinimalRepresentative_C1},
\ref{MinimalRepresentative_C2}
and \ref{MinimalRepresentative_C3}.
\end{theorem}
\proof For (i), the $1$-gon has to be preserved by any symmetry which leaves
$C_1$ and $C_s$ as the only possibilities.
They are both realized.
For (ii), we proceed in the same way.
(iii) is proved in Theorem \ref{Structure13spheres}. \qed

An interesting question is to consider whether a $(\{1,2,3\},6)$-sphere
can be mapped onto the projective plane $\mathbb{P}^2$.
This is clearly equivalent to the map having a central inversion.
$(\{1,3\},6)$-maps on the projective plane do not exist since no
centrally-symmetric $(\{1,3\},6)$-sphere exist.
All $(\{2,3\},6)$-maps on the $\mathbb{P}^2$ are antipodal quotients
(i.e. with halved $p$-vector and number of vertices)
of $(\{2,3\},6)$-spheres whose groups contain the inversion, i.e. 
those of symmetry $C_i$, $C_{2h}$, $D_{2h}$, 
$D_{3d}$, $D_{6h}$, $S_6$ and $T_h$.
In the next section we will describe explicitly the $(\{2,3\},6)$-sphere
of symmetry $D_{6h}$ and $T_h$.

\section{The Goldberg-Coxeter construction}\label{GoldbergCoxeterConstruction}
In \cite{goldberg} a construction is given,
generalizing Goldberg-Coxeter construction given in \cite{Gold,Cox71}
for $3$-regular graphs with $6$-gonal and $5$-, $4$-, $3$-gonal faces only. 
For the particular case when $G$ is a geometrical fullerene,
there is a large body of literature, see bibliography of \cite{goldberg}.
It takes a $3$- or $4$-regular plane graph $G$ and return a $3$-
or $4$-regular plane graph.
There the first step was to take the dual and get a triangulation
or a quadrangulation of the sphere.
The respective triangles and squares were subdivided, then put together
and the dual was taken.
An instrument in this operation was that the Eisenstein and Gaussian
integers are best represented on the tiling of the plane by equilateral
triangles, respectively, squares.
We are able to generalize this construction to the $6$-regular case but
there are differences.

\begin{figure}
\begin{center}

\begin{minipage}[b]{6.2cm}
\centering
\epsfig{height=46mm, file=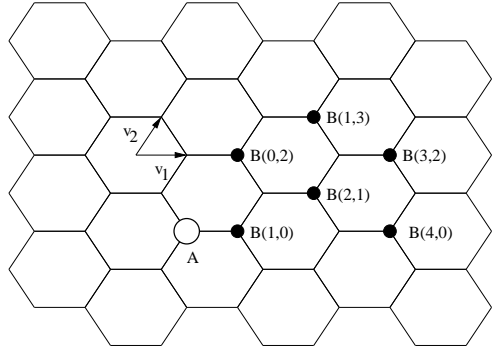}\par
\end{minipage}

\end{center}
\caption{The tiling by hexagons, the point $A$ and some points in the other bipartite component}
\label{HexagonalSystem}
\end{figure}

First, if $G$ is a $(\{2,3\}, 6)$-sphere, then the dual $G^{*}$
is a plane graph with faces of size $6$ and thus, bipartite.
The tessellation of Euclidean plane by regular hexagons is represented on
Figure \ref{HexagonalSystem}.
We use there two vectors $v_1$, $v_2$ to represent the coordinate of
the points. In complex coordinates $v_1=1$ and $v_2=j$ with $j=e^{i\pi/3}$.
The lattice $L=\ZZ v_1 +\ZZ v_2$ is called the {\em Eisenstein ring}.
The point $A$ is the origin and the point $B(k,l)$ is the point $k+l j$.
The points in the bipartite component of $A$ are $L_A=(1+j)L$, while
the points in the component of $B(1,0)$ are $L_B=1 + (1+j)L$.
Both sets $L_A$ and $L_B$ are stable under multiplication.
We will first define the Goldberg-Coxeter construction for $k+lj\in L_B$.
Then we will extend it to any $(k,l)\not= 0$.

\begin{figure}
\begin{center}
\begin{minipage}[b]{3.2cm}
\centering
\epsfig{height=20mm, file=Bundle6.eps}\par
$GC_{1,0}(6\times K_2)$, $D_{6h}$, $n=2$
\end{minipage}
\begin{minipage}[b]{3.2cm}
\centering
\epsfig{height=27mm, file=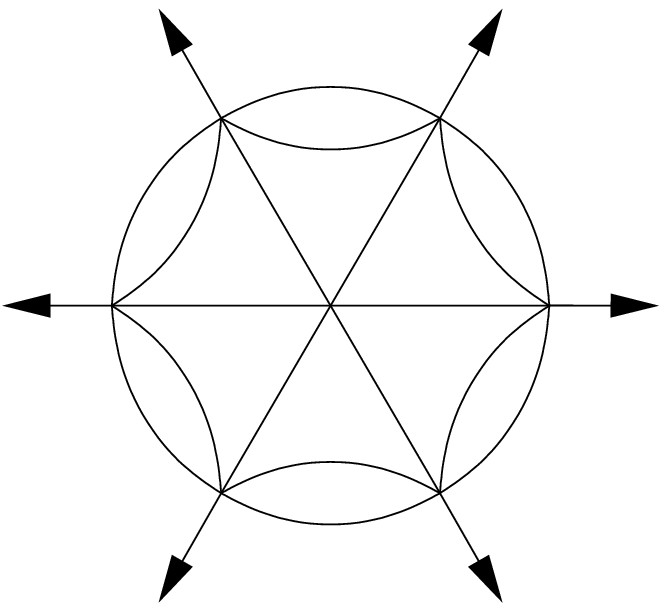}\par
$GC_{2,0}(6\times K_2)$, $D_{6h}$, $n=8$
\end{minipage}
\begin{minipage}[b]{3.2cm}
\centering
\epsfig{height=27mm, file=SmallestD6thi.eps}\par
$GC_{2,1}(6\times K_2)$, $D_6$, $n=14$
\end{minipage}
\begin{minipage}[b]{3.2cm}
\centering
\epsfig{height=27mm, file=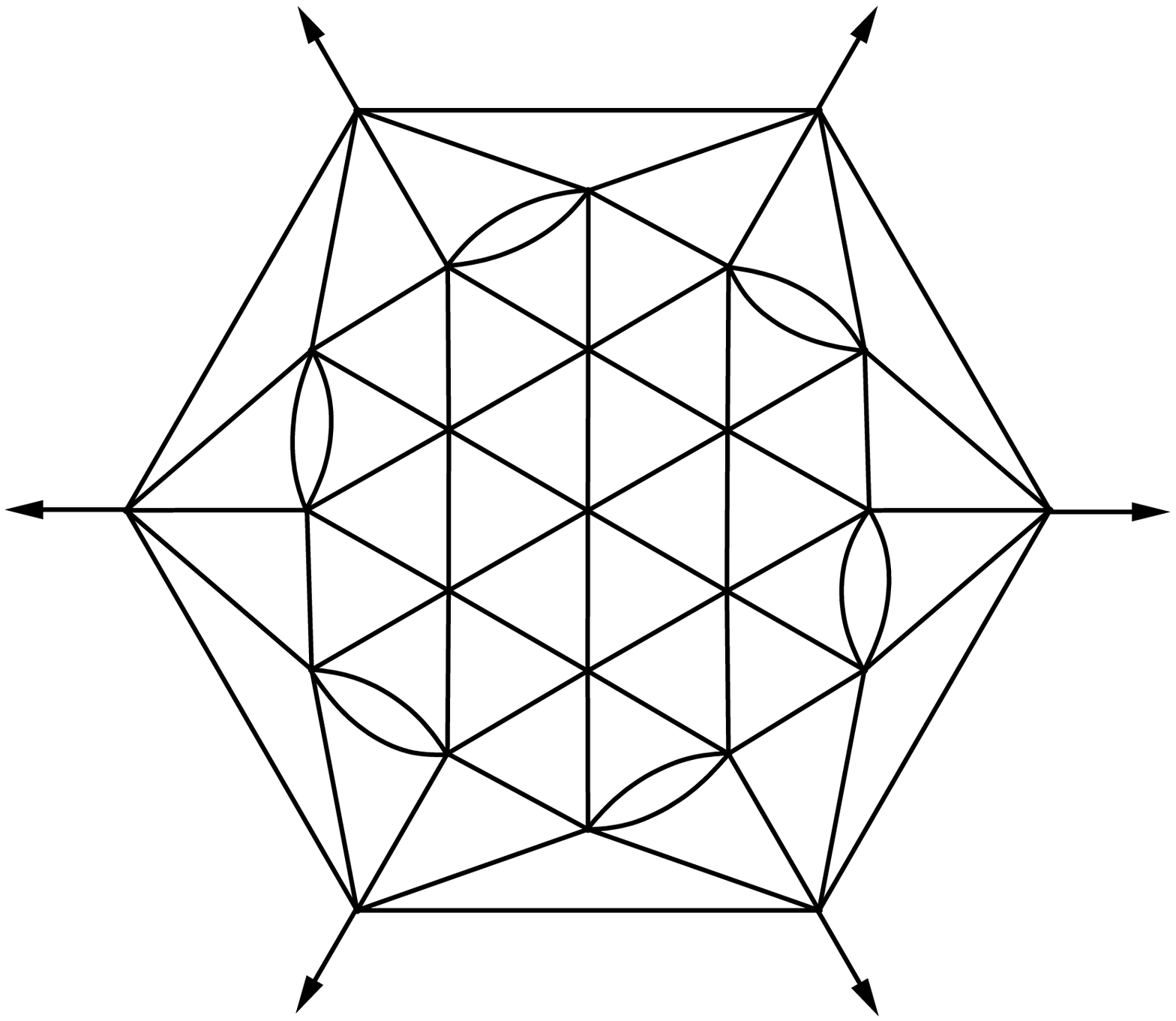}\par
$GC_{3,1}(6\times K_2)$, $D_6$, $n=26$
\end{minipage}
\begin{minipage}[b]{3.2cm}
\centering
\epsfig{height=27mm, file=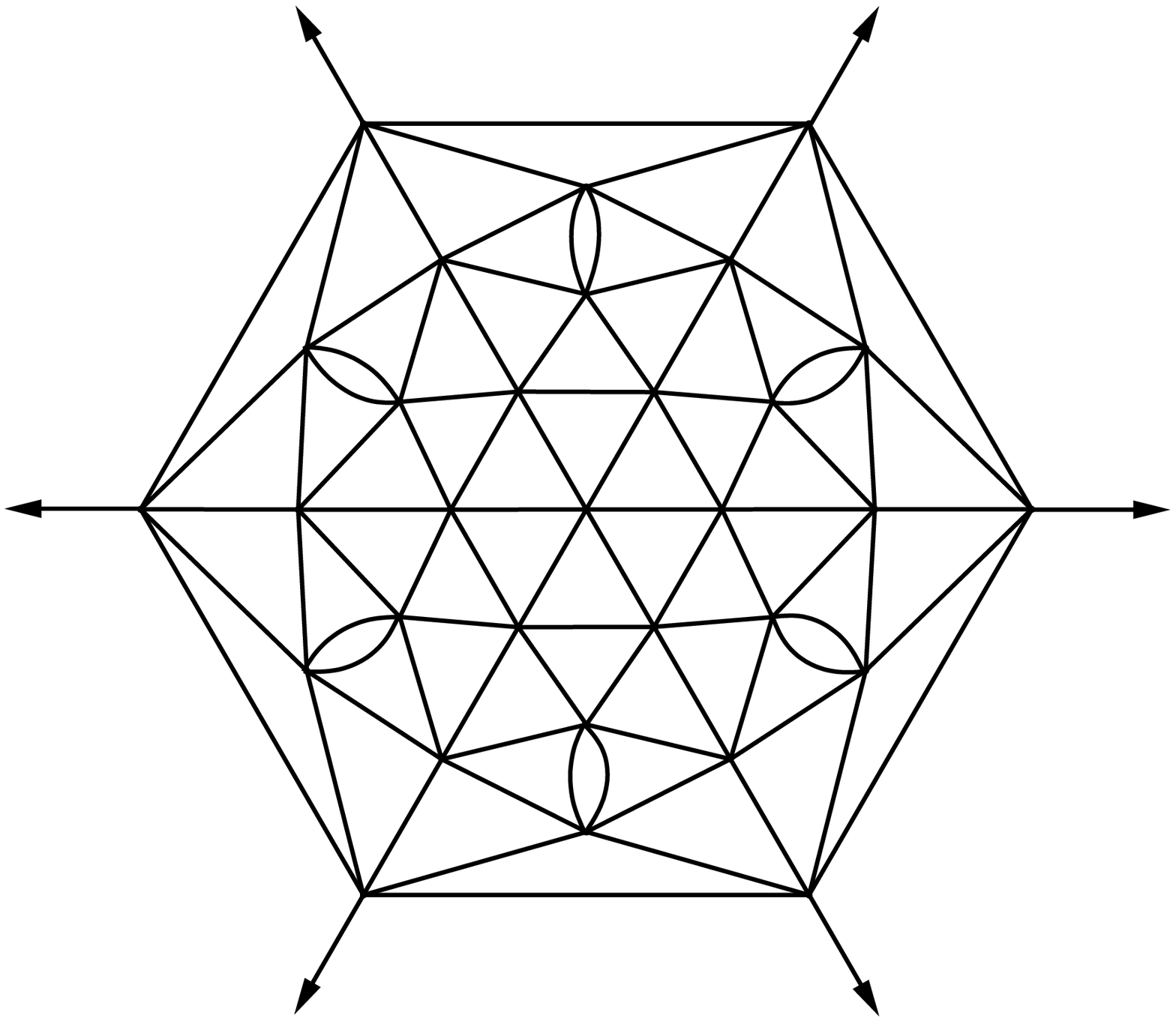}\par
$GC_{4,0}(6\times K_2)$, $D_{6h}$, $n=32$
\end{minipage}
\begin{minipage}[b]{3.2cm}
\centering
\epsfig{height=27mm, file=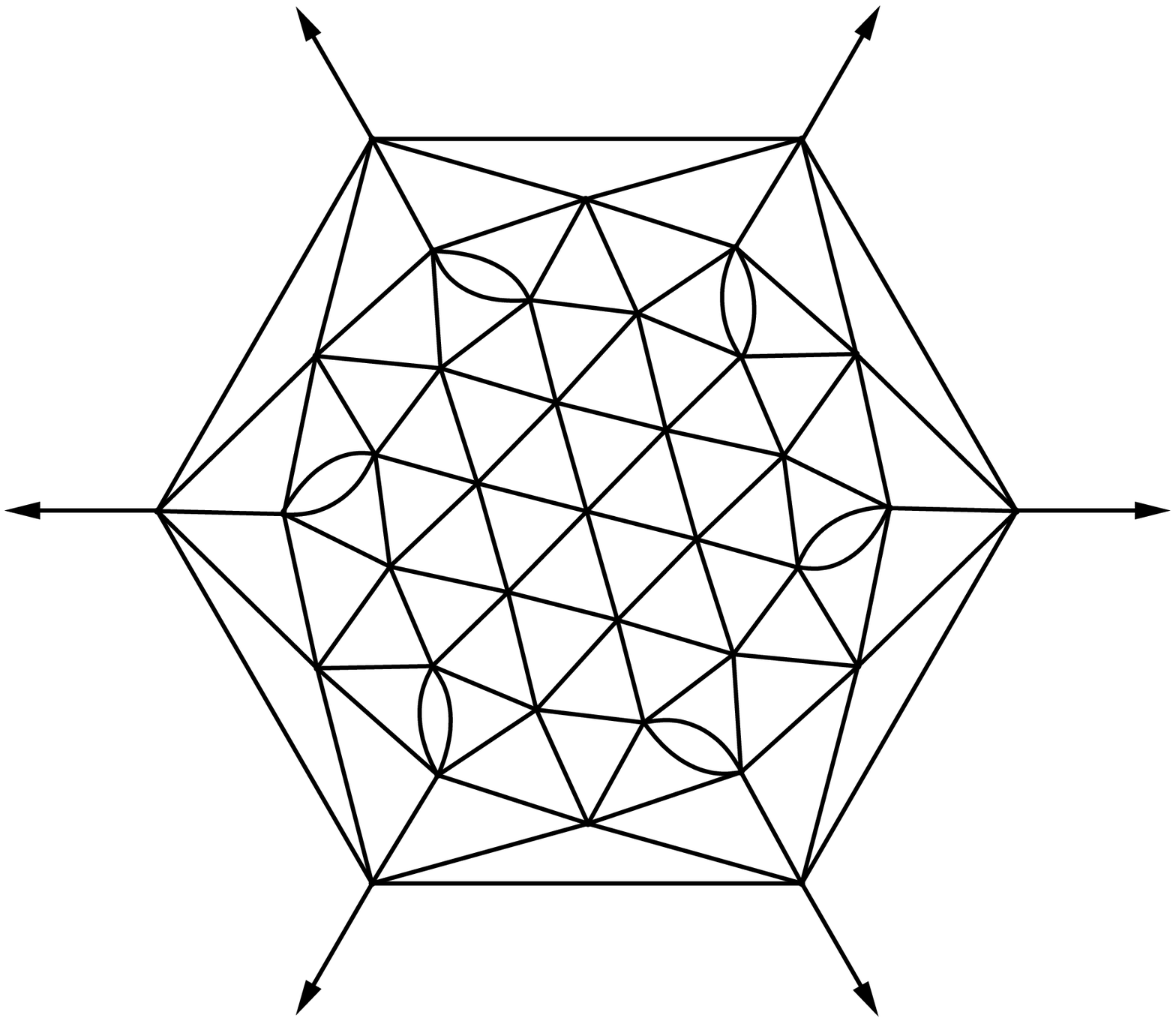}\par
$GC_{3,2}(6\times K_2)$, $D_6$, $n=38$
\end{minipage}

\end{center}
\caption{Smallest $(\{2,3\}, 6)$-spheres of symmetry $D_6$, $D_{6h}$
in terms of the Goldberg-Coxeter construction}
\label{GCconstructionParam}
\end{figure}

\begin{theorem}\label{FirstDefinitionGCkl}
If $z=k + l j\in L_B$ and $G_0$ is a $6$-regular plane graph with $|G_0|$
vertices, then it is possible to define a plane
graph $G'=GC_{z}(G_0)=GC_{k,l}(G_0)$ such that the following holds:

(i) $G'$ is a $6$-regular plane graph with $|G_0| (k^2 + kl + l^2)$ vertices.

(ii) Every face of $G_0$ corresponds to a face of $G'$ with all new faces
of $G'$ being $3$-gons.

(iii) $G'$ has all rotational symmetries of $G_0$ and all symmetries
as well if $l=0$ or $k=0$.

(iv) $GC_{1,0}(G_0)=G_0$ and $GC_{z}(G_0)=GC_{z j^2} (G_0)$.

(v) $GC_{z}(GC_{z'}(G_0)) = GC_{zz'}(G_0)$.

(vi) $GC_{z}(G_0)=GC_{\overline{z}}(\overline{G_0})$ where $\overline{G_0}$
is the graph that differs from $G_0$ only by a plane symmetry.

\end{theorem}
\proof Let $G_0$ be a $6$-regular graph.
The dual $G_0^{*}$ is a plane graph with all faces being $6$-gons.
If $z=k + l j \in L_B$, then the point $B(k,l)$ belongs to
the same connected component as $B$.

The point $c=-j^2 z$ is the center of an hexagon and we build
around it six points $P_q$:
\begin{equation*}
P_q = c - j^q c \mbox{~for~} 0\leq q\leq 5.
\end{equation*}
Those six points form a {\em master hexagon} that correspond to the
original hexagon.
Every hexagon of $G_0^{*}$ can be thus modified and we can arrange them
together at the boundary between adjacent hexagons. We can thus
obtain another plane graph with $6$-regular faces. By taking the dual
one more time, we get $GC_{k,l}(G_0)$.
Checking the remaining properties is relatively easy. \qed

\begin{figure}
\begin{center}
\epsfig{height=20mm, file=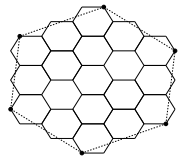}\par
\end{center}
\caption{Local structure of the Goldberg-Coxeter construction $GC_{3,2}$}
\label{LocalStructureGC32}
\end{figure}

The above theorem is similar to Proposition 3.1 in \cite{goldberg}.
But there are some differences. In the $3$-regular case, we have $GC_{k,l}(G_0)$
with all symmetries if $k=l$, while here the case $k=l$ is impossible.
See in Figure \ref{LocalStructureGC32} the local structure of the
Goldberg-Coxeter construction $GC_{3,2}$ and
in Figure \ref{LocalStructureGC40}, $GC_{4,0}$.

\begin{theorem}
The $(\{2,3\}, 6)$-spheres of symmetry $D_6$, $D_{6h}$ are obtained as $GC_{k,l}(6\times K_2)$ with $k+lj \in L_B$.
\end{theorem}
\proof Let us take a $(\{2,3\}, 6)$-sphere $G$ of symmetry $D_6$ or $D_{6h}$
and let us take the dual $G^{*}$.
The $6$-fold axis passes through a $6$-gon $F$ and the $2$-gons of $G$ correspond
to $6$-regular vertices.
But the position of those $2$-gons define a master hexagon around $F$ and
thus we get exactly the structure of a graph $GC_{k,l}(6\times K_2)$. \qed

In Theorem \ref{FirstDefinitionGCkl} we have defined the Goldberg-Coxeter
construction $GC_{k,l}$ for $k+l j\in L_B$. Now we want to define it
for any $k,l\not= 0$.
For that we first introduce the notion of 
{\em oriented tripling}.

\begin{figure}
\begin{center}
\begin{minipage}[b]{9.2cm}
\centering
\epsfig{height=40mm, file=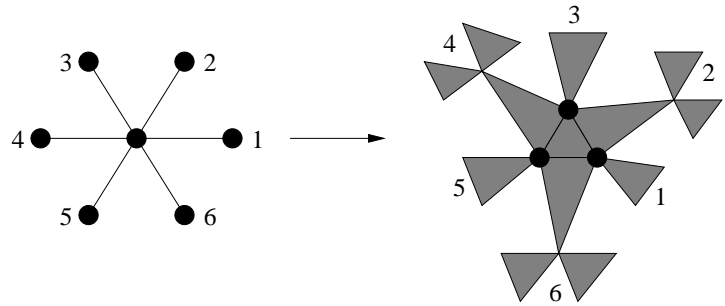}\par
\end{minipage}

\end{center}
\caption{Local configuration around a vertex of the oriented tripling operation.}
\label{LocalConfigurationTripling}
\end{figure}

\begin{definition}
If $G$ is a $6$-regular plane graph, then its dual $G^{*}$ is bipartite.
For each such bipartite class $C$ we define a graph $Or_C(G)$ with
the following properties:
\begin{enumerate}
\item[(i)] $Or_C(G)$ is a $6$-regular plane graph with $3$ times as many vertices.
\item[(ii)] Each vertex of $G$ corresponds to $3$ vertices of $Or_C(G)$ and $4$ triangular faces.
\item[(iii)] Every symmetry of $G$ preserving $C$ also occur as symmetry of $Or_C(G)$.
\end{enumerate}
The local configuration of the operation is shown in
Figure \ref{LocalConfigurationTripling}.
For every face $F$ of $G$, we orient the edges of $F$ counter-clockwise.
Thus for every bipartite class $C$ of $G^{*}$ we get an orientation
of the edges of $G$.
Around a vertex $v$ and its six adjacent vertices, there are three
vertices $w$ to which the edge $\{v,w\}$ is oriented from $v$ to $w$.
They are the vertices $1,3,5$ in Figure \ref{LocalConfigurationTripling}.
\end{definition}

\begin{figure}
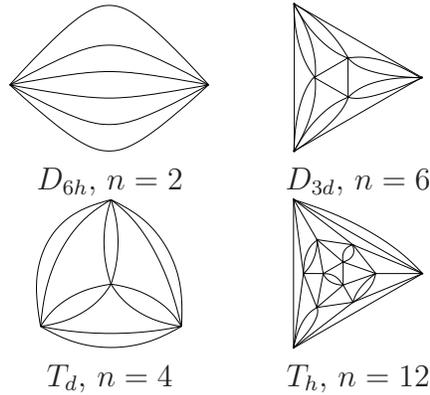

\begin{center}
\begin{minipage}[b]{6.5cm}
\begin{minipage}[b]{3.2cm}
\centering
\epsfig{height=20mm, file=Bundle6.eps}\par
$D_{6h}$, $n=2$
\end{minipage}
\begin{minipage}[b]{3.2cm}
\centering
\epsfig{height=20mm, file=SmallOcta_D3dSec.eps}\par
$D_{3d}$, $n=6$
\end{minipage}
\end{minipage}
\begin{minipage}[b]{6.5cm}
\begin{minipage}[b]{3.2cm}
\centering
\epsfig{height=20mm, file=Example23_6val_2.eps}\par
$T_d$, $n=4$
\end{minipage}
\begin{minipage}[b]{3.2cm}
\centering
\epsfig{height=20mm, file=IcosahedronThSec.eps}\par
$T_h$, $n=12$
\end{minipage}
\end{minipage}

\end{center}
\caption{Two examples of $(\{2,3\},6)$-spheres with unique oriented tripling}
\label{ExamplesTripling}
\end{figure}

\begin{figure}
\begin{center}
\begin{minipage}[b]{2.9cm}
\centering
\epsfig{height=20mm, file=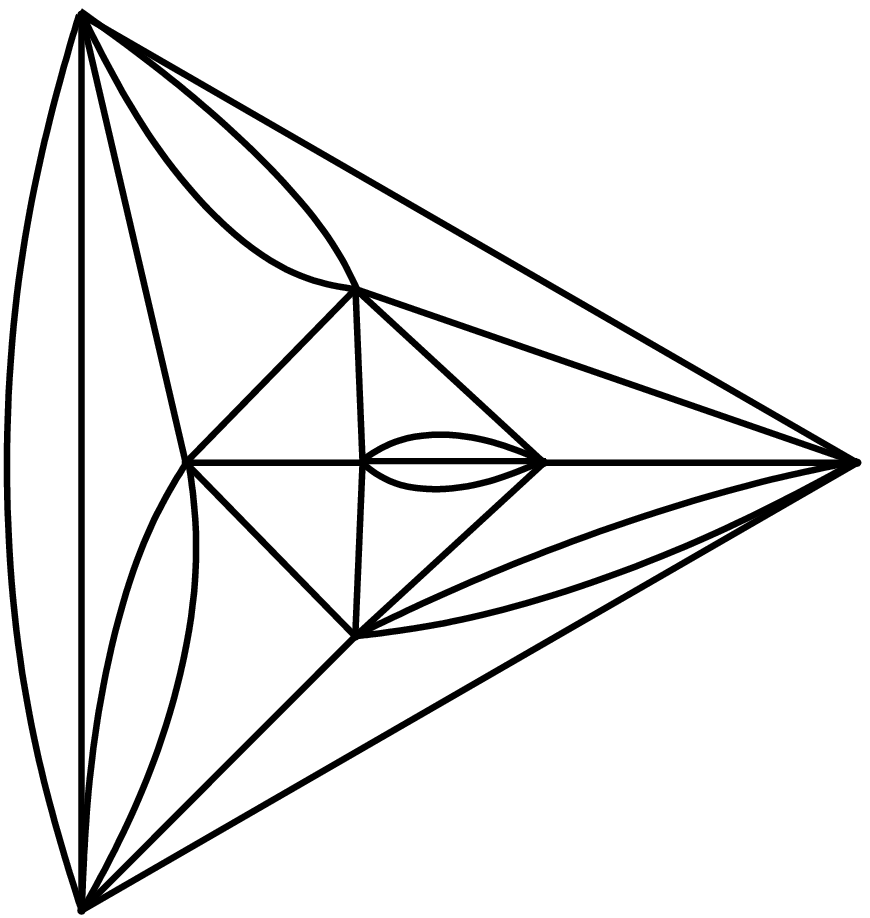}\par
$G$: $C_{1}$, $n=8$
\end{minipage}
\begin{minipage}[b]{3.9cm}
\centering
\epsfig{height=32mm, file=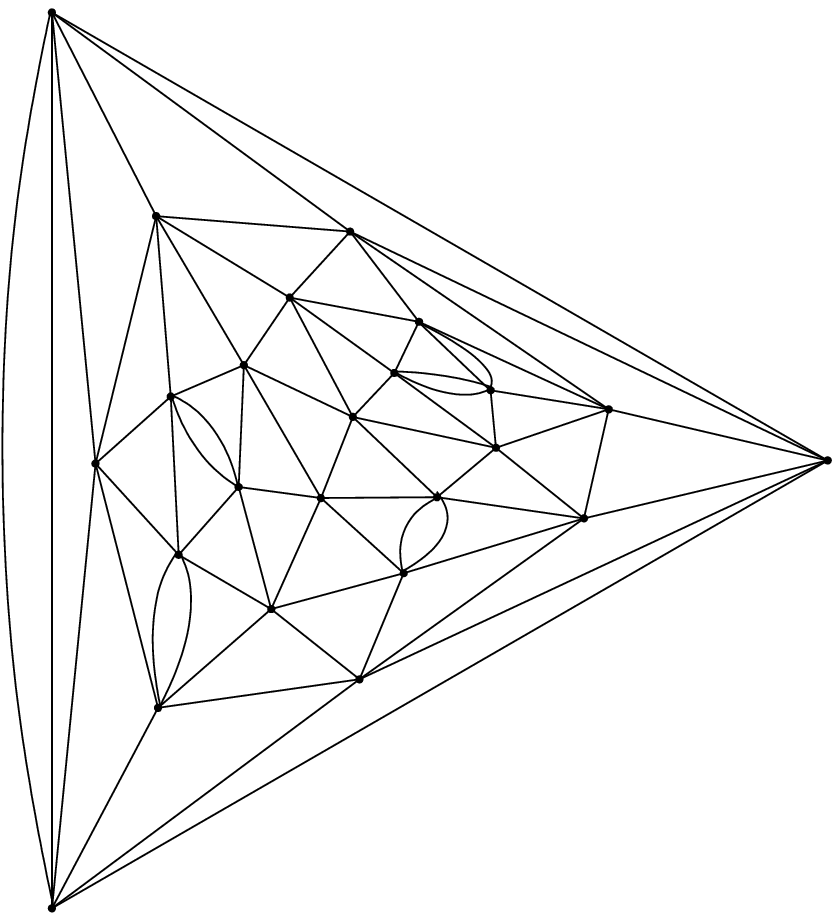}\par
$Or_{C_1}(G)$: $C_{1}$, $n=24$
\end{minipage}
\begin{minipage}[b]{3.9cm}
\centering
\epsfig{height=32mm, file=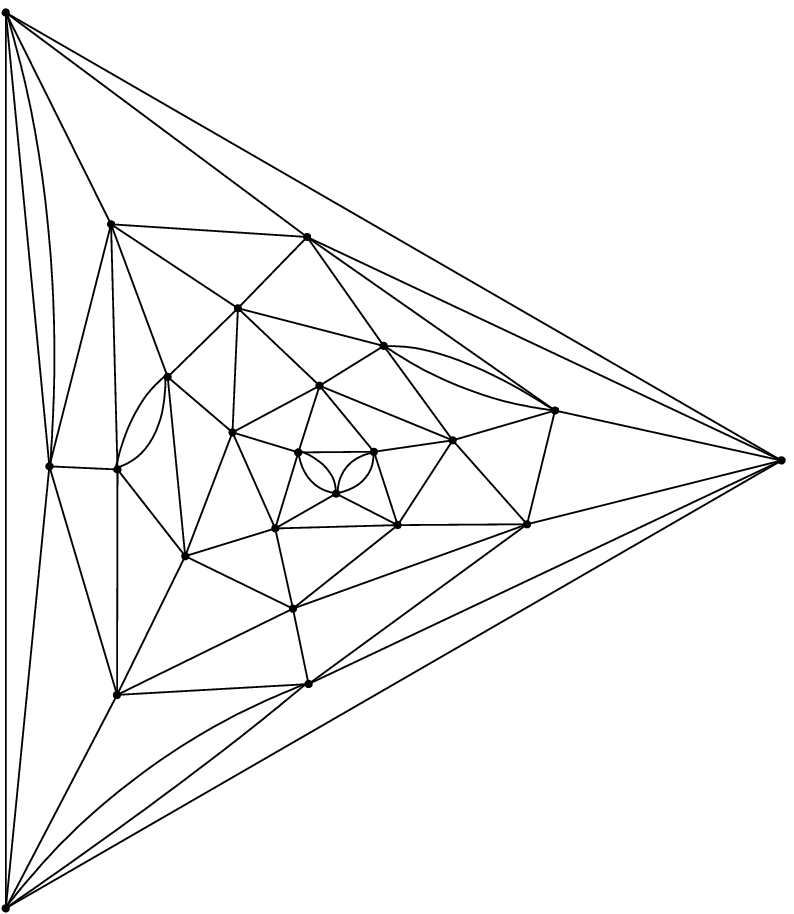}\par
$Or_{C_2}(G)$: $C_{1}$, $n=24$
\end{minipage}

\end{center}
\caption{Smallest $(\{2,3\},6)$-sphere having two non-isomorphic oriented triplings}
\label{ExampleGraphWithNoUniqueOrientedTripling}
\end{figure}

So, if $G$ has two inequivalent bipartite components $C_1$ and $C_2$,
then $Or_{C_1}(G)$ and $Or_{C_2}(G)$ are not necessarily isomorphic and
the smallest such example is shown on
Figure \ref{ExampleGraphWithNoUniqueOrientedTripling}.
In Figure \ref{ExamplesTripling} we give two examples of the action of
the oriented tripling when the obtained graph is unique.

\begin{figure}
\begin{center}
\begin{minipage}[b]{2.9cm}
\centering
\epsfig{height=20mm, file=Class_C3_C3vB.eps}\par
$C_{3v}$, $n=1$
\end{minipage}
\begin{minipage}[b]{2.9cm}
\centering
\epsfig{height=20mm, file=Singular13_3vert.eps}\par
$C_{3h}$, $n=3$
\end{minipage}
\begin{minipage}[b]{2.9cm}
\centering
\epsfig{height=20mm, file=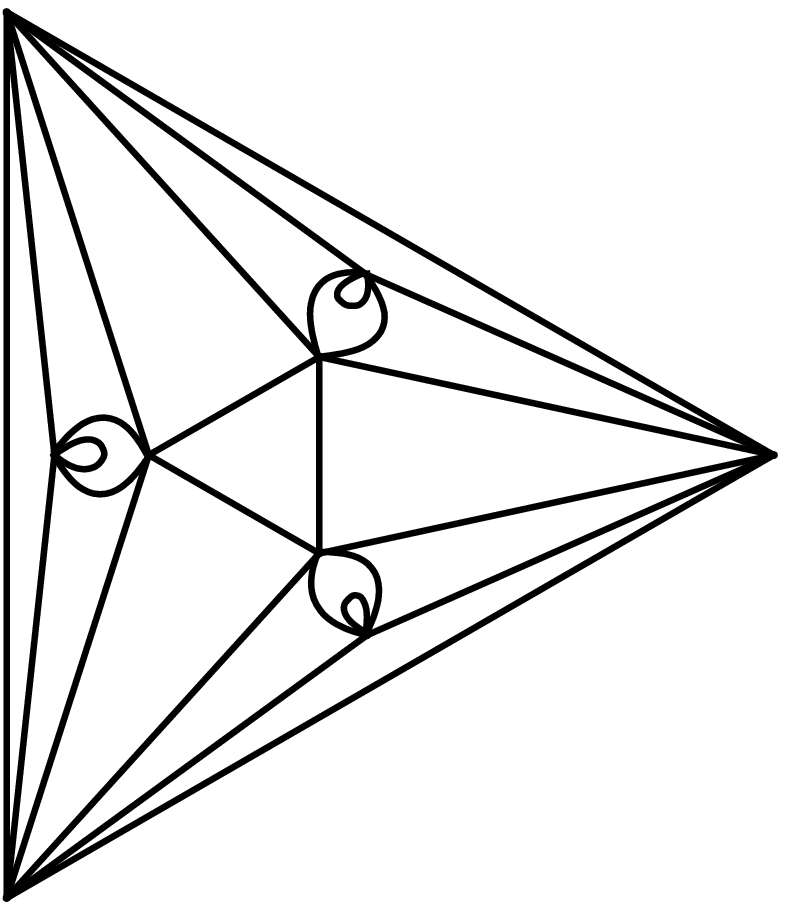}\par
$C_{3v}$, $n=9$
\end{minipage}
\begin{minipage}[b]{2.9cm}
\centering
\epsfig{height=20mm, file=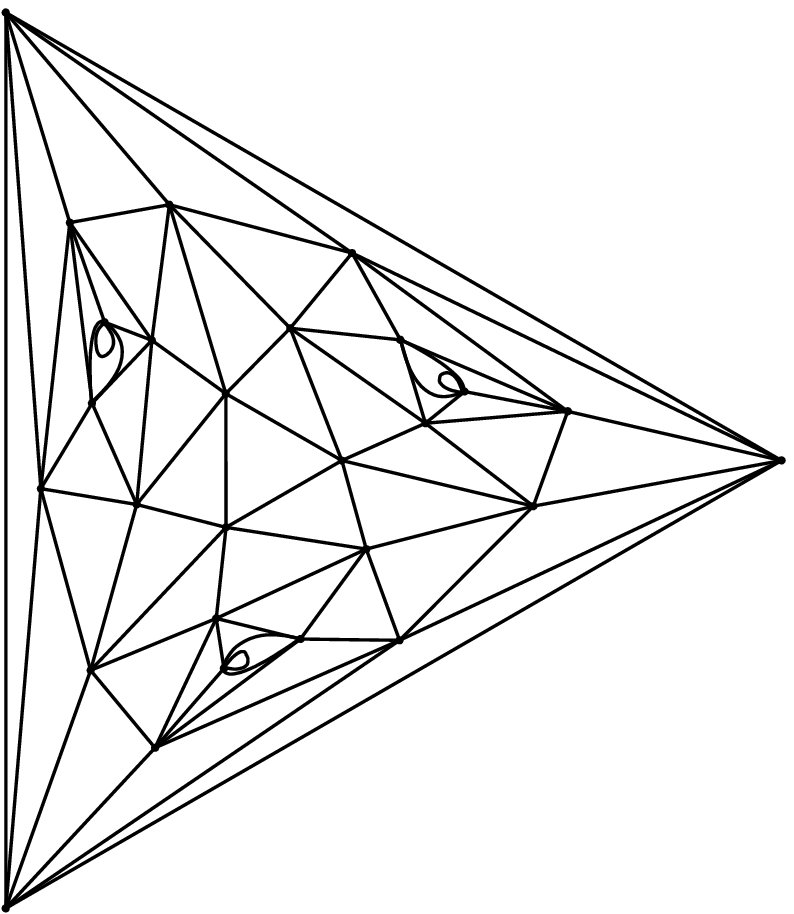}\par
$C_{3h}$, $n=27$
\end{minipage}

\end{center}
\caption{First terms of the infinite sequence of $(\{1,3\}, 6)$-spheres $T_{i}$}
\label{FirstTermInfiniteSequence}
\end{figure}

For the Trifolium $T_1$,
we can define a sequence $T_{i}$ of graphs
with $T_{i+1}$ obtained by applying the oriented tripling to $T_i$.
The first $4$ terms are shown in Figure \ref{FirstTermInfiniteSequence}.

We now introduce the Goldberg-Coxeter construction in the general case.
For a sphere $G$, denote by $Tr(G)$ the truncation of $G$, i.e. 
the sphere obtained by replacing every vertex of degree $k$ of $G$
by a $k$-gonal face.

We will also use the following result: if $G$ is a $3$-regular sphere
with faces of even size, then it is possible to color the faces of $G$
so that any two adjacent faces have different colors. Such a 
coloring is unique up to permutation of the colors.
If $G_0$ is a graph with vertices of even degree, then its dual is bipartite
and the three colors in $Tr(G_0)$ come from the vertices of $G_0$ and
the two classes of faces in $G_0$.

\begin{theorem}\label{GoldbergCoxeterNextGeneration}
For a $6$-regular plane graph $G_0$ and two integers $k,l$ with $k,l\not=0$,
we can define two $6$-regular
spheres $G_1$, $G_2$ with $GC_{k,l}(G_0)=\{G_1, G_2\}$.
This will satisfy to the following properties:

(i) $Tr(G_i) = GC_{k,l}(Tr(G_0))$ for $i=1$, $2$ with
$GC_{k,l}$ being the Goldberg-Coxeter construction for $3$-regular spheres.

(ii) $G_1$ and $G_2$ are $6$-regular plane graphs with $|G_0|(k^2+kl+l^2)$ vertices.

(iii) Every face of $G_0$ corresponds to a face of $G_1$ and $G_2$ with
all new faces of $G_1$ and $G_2$ being $3$-gons.

(iv) $GC_{1,1}(G_0) = \{Or_C(G_0), Or_{C'}(G_0)\}$.

(v) If $k+l j\in L_B$, then $G_1=G_2$.

(vi) If $GC_{k,l}(G_0)=\{G_1, G_2\}$, then $GC_{k', l'}(G_1)=GC_{k_p, l_p}(G_0)$ with $k_p+l_p j = (k+l j) (k'+l' j)$.

\end{theorem}
\proof Let us take a $3$-coloring in white, red, blue of the faces
of $Tr(G_0)$ with white corresponding to the faces coming from vertices
of $G_0$.
The $3$-regular sphere $GC_{k,l}(Tr(G_0))$ has the faces of $G_0$ and
some $6$-gonal faces, thus all its faces are of even size and we can
find a $3$-coloring of them.

One can see directly that all the white faces of $Tr(G_0)$ have the same
color in $GC_{k,l}(Tr(G_0))$; we color them white.
If $k\equiv l\pm 1 \pmod 3$ (this contains the case $k+lj\in L_B$),
then the faces of $Tr(G_0)$ coming from faces
of $G_0$ will not be white in $GC_{k,l}(Tr(G_0))$.
Thus by shrinking the white faces, we get a graph, which is actually the
$GC_{k,l}(G_0)$ defined in Theorem \ref{FirstDefinitionGCkl} if $k+lj\in L_B$.

If $k\equiv l\pmod 3$ then all faces of $Tr(G_0)$ will correspond to
white faces in $GC_{k,l}(Tr(G_0))$. The remaining $6$-gonal faces have color
red and blue. This gives two set of faces that can be shrunk and thus
two possible graphs.
All properties follow easily. \qed

\begin{theorem}\label{FactorizationProperties}
(i) Any $k+l j\not=0$ can be written as $k+l j = (1+j)^s (k'+l' j) j^u$ with $s\geq 0$, $u\in \{0,1\}$ and $k'+l'j\in L_B$.

(ii) The sphere $GC_{k,l}(G_0)$ is obtained by applying the oriented
tripling $s$ times and then the Goldberg-Coxeter construction from
Theorem \ref{FirstDefinitionGCkl}.
\end{theorem}
\proof (i) The ring of Eisenstein integers is a unique factorization domain.
That is every $k+lj\not= 0$ can be factorized by into the relevant primes.
The condition $k\equiv l\pmod 3$ is equivalent to $k+lj$ being divisible
by $1+j$.
Thus by repeated application of this we can write
\begin{equation*}
k+lj = (1+j)^s (k_2 +l_2 j) \mbox{~with~} k_2\equiv l_2\pm 1\pmod 3.
\end{equation*}
If $k_2 \equiv l_2 +1\pmod 3$, then we are done, otherwise we divide by $j$.

(ii) follows from the multiplicativity property (vi) of Theorem \ref{GoldbergCoxeterNextGeneration}. \qed

This idea of using the truncation and resulting $3$-regular
spheres was, perhaps, used for the first time in \cite{GrZa}.
This idea could in principle
be applied to the enumeration of the $(\{2,3\}, 6)$-spheres, since the
$(\{4,6\}, 3)$-spheres can be obtained from the {\tt CPF} program.
But the truncation
multiplies the number of vertices by $6$ and this
makes this method uncompetitive to the one of
Section \ref{GenerationMethod}.

We cannot say much in general for the symmetry groups of $GC_{k,l}(G_0)$.
This is essentially the same situation as for the oriented tripling.
What happens is that for $3$-regular graphs,
Goldberg-Coxeter construction $GC_{k,l}$
preserve all symmetries if $k=0$ or $k=l$ and only rotational symmetries 
otherwise. Thus we get the automorphism group $\Gamma$ of $GC_{k,l}(Tr(G_0))$.
If $\Gamma$ preserves the set of faces of color red and blue, then
$\Gamma$ is a group of symmetries of $GC_{k,l}(G_0)$, otherwise the stabilizer
of the red faces is a group of symmetries of $GC_{k,l}(G_0)$.
But some accidental symmetries can occur and we have thus to work on a 
case-by-case basis.

\begin{theorem}\label{Structure13spheres}
Let $G$ be a $(\{1,3\},6)$-sphere. The following hold:

(i) $G=GC_{k,l}(Trifolium)$ with $0\leq l\leq k$ and has $k^2+kl+l^2$ vertices.

(ii) $G$ has symmetry $C_{3v}$ if $k=0$, $C_{3h}$ if $k=l$ and $C_3$ otherwise.

(iii) $G$ is obtained as $GC_{k,l}(T_i)$ with $k+lj\in L_B$, where $(T_i)_{i\geq 1}$ is the infinite series of $6$-regular graphs obtained by starting from Trifolium.
\end{theorem}
\proof Let us take a $(\{1,3\}, 6)$-sphere.
Then $Tr(G)$ is a $(\{2,6\}, 3)$-sphere. 
Either from \cite{GrZa} or \cite{T}, we know that such spheres
are obtained as $GC_{k,l}(3\times K_2)$ 
with $GC_{k,l}$ denoting here the $3$-regular Goldberg-Coxeter
construction.
Since the faces of $GC_{k,l}(3\times K_2)$ are of even size, it is possible
to define a $3$-coloring of those faces.
The $2$-gonal faces should not be in all $3$ different colorings.
This can occur only if $k\equiv l\pmod 3$.
So, $k+l j$ can be factorized as $(1+j) (k'+l'j)$ and we get
\begin{equation*}
\begin{array}{rcl}
Tr(G) &=& GC_{k,l}(3\times K_2)\\
      &=& GC_{k',l'}(GC_{1,1}(3\times K_2)\\
      &=& GC_{k',l'}(Tr(Trifolium)).
\end{array}
\end{equation*}
Thus we have proved (i).

The symmetry of $GC_{k,l}(3\times K_2)$ is $D_{3h}$ if $k=0$ or $k = l$ and
$D_3$ otherwise.
If $k\equiv l\pmod 3$ then all $2$-gons of $GC_{k,l}(3\times K_2)$ are in the
same color, say white.
The $3$-gonal faces that are not white are of two possible colors red and
blue.
In order for a symmetry of $Tr(G)=GC_{k,l}(3\times K_2)$ to induce a symmetry
of $G$ it is necessary and sufficient that it preserves all $3$ colors
of the coloring. This reduces by a factor of $2$ the symmetry group and
we get $C_3$, $C_{3h}$ and $C_{3v}$ as possible groups. 

Statement (iii) follows from
Theorem \ref{FactorizationProperties} (ii). \qed

Note that the possible number of vertices of $(\{2,6\},3)$-spheres
was already determined in \cite{GrZa}.

Denote by $K_2\times Tetrahedron$ the Tetrahedron with edges doubled.

\begin{theorem}
(i) Any $(\{2,3\}, 6)$-sphere of symmetry $T$, $T_h$ or $T_d$ 
is obtained as $GC_{k,l}(K_2\times Tetrahedron)$.

(ii) The $(\{2,3\},6)$-spheres of symmetry $T_d$, respectively $T_h$,
are of the form $GC_{k,0}(K_2\times Tetrahedron)$, respectively
$GC_{k,k}(K_2\times Tetrahedron)$.

\end{theorem}
\proof Take a $(\{2,3\},6)$-sphere $G$ of symmetry $T$, $T_d$ or $T_h$
and consider their truncation $Tr(G)$.
It is a $(\{4,6\},3)$-sphere which contains a subgroup $T$ of symmetry.
By Theorem 6.2 of \cite{zig2}, this implies that the symmetry group of
$Tr(G)$ is $O$ or $O_h$.
By \cite{goldberg} Theorem 5.2, $Tr(G)$ is described by the Goldberg-Coxeter
construction applied to the cube, i.e. $Tr(G)=GC_{k,l}(Cube)$.
We need now to determinate which graphs $GC_{k,l}(Cube)$ are of the form
$Tr(G)$. For that we need to consider the $3$-coloring of the faces.
It is necessary that the $4$-gonal faces are not in all $3$ colors
of the faces. This implies that $k\equiv l\pmod 3$.
In turn this give us that $k+lj = (1+j) (k'+l'j)$, which then gives us
\begin{equation*}
\begin{array}{rcl}
Tr(G) &=& GC_{k,l}(Cube)\\
      &=& GC_{k',l'}(GC_{1,1}(Cube))\\
      &=& GC_{k',l'}(Tr(K_2\times Tetrahedron)).
\end{array}
\end{equation*}
(i) follows from
Theorem \ref{GoldbergCoxeterNextGeneration}.

If a $(\{2,3\}, 6)$-sphere is of symmetry $T_d$ or $T_h$, then the symmetry
group of the truncation is $O_h$ and such spheres are described as $GC_{k,0}(Cube)$ and $GC_{k,k}(Cube)$. \qed

\begin{theorem}\label{TheoremParameterization}
The number of $(\{1,2,3\}, 6)$-spheres with $i$ $1$-gons and less than $n$ vertices grows like $O(n^{4-i})$.
\end{theorem}
\proof Take $G$ a $(\{1,2,3\}, 6)$-sphere with $n$ vertices and $i$ $1$-gons.
Then $Tr(G)$ is a $(\{2,4,6\}, 3)$-sphere with $i$ $2$-gons, $6-2i$ $4$-gons and $6n$-vertices.
Thus the number of faces of size $2$ or $4$ is $6-i$.
The $3$-regular plane graphs whose faces have size at most $6$ and the set of
faces of size less than $6$ is fixed are described by 
the parametrization theory of Thurston \cite{T}.
By using it \cite{sah94} obtained some upper bound on the number of
geometric fullerenes.
The proof applies just as well for the other classes of graphs and give us
the required upper bound. \qed

Note that in principle, Thurston's theory allows to say more.
First it gives that the $(\{2,4,6\},3)$-spheres with $i$ $2$-gons. 
are described by $4-i$ Eisenstein integers. Not all such spheres
correspond to $(\{1,2,3\},6)$-spheres with $i$ $1$-gons. For that
some congruence have to be satisfied.

\section{Zigzags and central circuits}\label{ZigzagCentralCircuitSection}

For a plane graph $G$ and a zigzag or central circuit,
if we change the orientation, then 
the type of intersection does not change.
Thus, to a zigzag or central circuit of length $l$ with $\alpha_1$
and $\alpha_2$ intersections of type I and II,
we attribute the
symbol $l_{\alpha_1, \alpha_2}$ and we define the
{\em $z$-vector}, respectively, {\em $c$-vector} 
${l_1}^{m_1}_{\alpha_{1,1}, \alpha_{1,2}}, \dots, {l_p}^{m_p}_{\alpha_{p,1}, \alpha_{p,2}}$
to be the vector enumerating such lengths with multiplicities $m_i$.

\begin{figure}
\begin{center}
\begin{minipage}[b]{6.2cm}
\centering
\epsfig{height=30mm, file=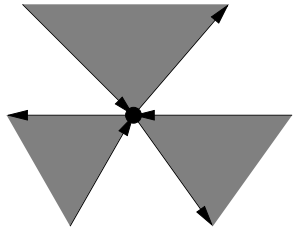}
\end{minipage}

\end{center}
\caption{The orientation of the edges of the face}
\label{TheOrient}
\end{figure}

\begin{theorem}\label{TypeIIalways}
For a $6$-regular plane graph, it is possible to find an orientation
on the zigzags and central circuits such that any edge, respectively,
vertex of intersection is of type II.
\end{theorem}
\proof Let us take a $6$-regular plane graph $G$.
Since every vertex is of even degree and $G$ is planar, the dual
graph $G^{*}$ is bipartite.
Let us take one color $c$ of the faces of $G$ and orient the edges of
the face of color $c$ in such a way that they turn clockwise
around the face (see Figure \ref{TheOrient}).
It is apparent that such orientation carries over to the zigzags and
central circuits of $G$ and that with this orientation all the intersection
are of type II. \qed

For zigzags, this result is not new, see for example \cite{zig2,Sh}.

\begin{figure}
\begin{center}
\begin{minipage}[b]{5.2cm}
\centering
\epsfig{height=22mm, file=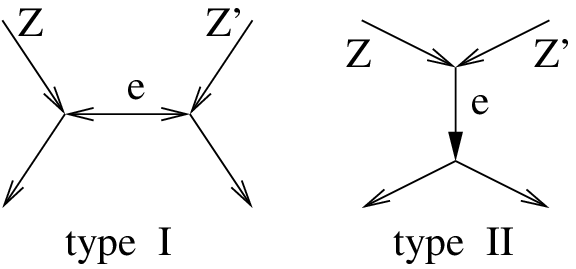}
Zigzag case
\end{minipage}
\begin{minipage}[b]{5.2cm}
\centering
\epsfig{height=22mm, file=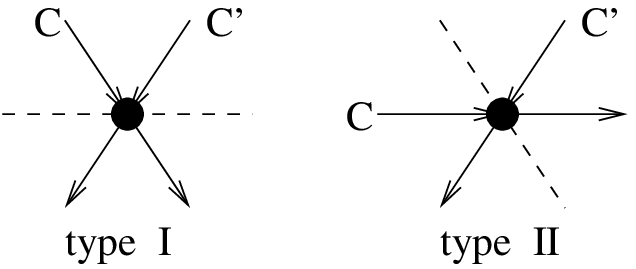}
Central circuit case
\end{minipage}

\end{center}

\caption{Intersection types of zigzags and central circuits}
\label{IntTypes_ZZ_CC}
\end{figure}

\begin{figure}
\begin{center}
\begin{minipage}[b]{6.2cm}
\centering
\epsfig{height=50mm, file=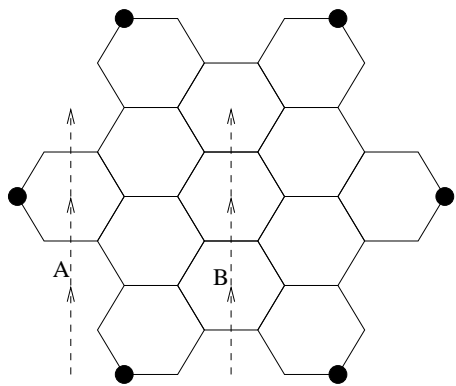}
\end{minipage}

\end{center}
\caption{The local structure of the Goldberg-Coxeter construction $GC_{4,0}$}
\label{LocalStructureGC40}
\end{figure}

\begin{theorem}
Let us take a $6$-regular plane graph $G$ with $z$-vector
$\dots, {l_i}^{a_i}_{\alpha_i, \beta_i}, \dots$ and
$c$-vector $\dots, {k_j}^{b_j}_{\alpha'_j, \beta'_j}, \dots$.
Then the $z$-vector and $c$-vector of $GC_{1+4u, 0}(G)$ are
\begin{equation*}
\dots, \left\{l_i(1+3u)\right\}_{\alpha_i, \beta_i}^{a_i(1+u)},\dots, \dots, \left\{2k_j(1+3u)\right\}_{\alpha'_j, \beta'_j}^{2u b_j}, \dots
\end{equation*}
and
\begin{equation*}
\dots, \left\{l_i\frac{1+3u}{2}\right\}_{\alpha_i,\beta_i}^{u a_i},\dots,
\dots, \left\{k_j(1+3u)\right\}_{\alpha'_j,\beta'_j}^{b_j(1+2u)},\dots
\end{equation*}

\end{theorem}
\proof The proof uses the Goldberg-Coxeter construction previously
built. One goes into the dual and subdivides the hexagons.
The picture in Figure \ref{LocalStructureGC40} shows that any
central circuit of $G$ corresponds to $1+2u$ central circuits 
(named $B$ in the figure) and that the zigzags in $A$ on one 
side correspond to zigzags in $G$.
The result follows similarly for $z$-vector. \qed

\begin{figure}
\begin{center}
\begin{minipage}[b]{3.2cm}
\centering
\epsfig{height=20mm, file=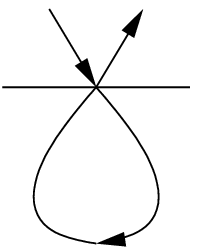}\par
Zigzag case
\end{minipage}
\begin{minipage}[b]{4.2cm}
\centering
\epsfig{height=20mm, file=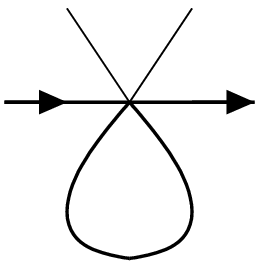}\par
Central circuit case
\end{minipage}

\end{center}
\caption{Self-intersection induced by a $1$-gon}
\label{SelfIntOneGonZigzagCentralCircuit}
\end{figure}

\begin{theorem}
(i) If a $(\{1,2,3\}, 6)$-sphere has a $1$-gon, then it has at
least one self-intersecting central circuit and one
self-intersecting zigzag.

(ii) For a $(\{1,3\}, 6)$-spheres, all central circuits and zigzags
are self-intersecting.
\end{theorem}
\proof (i) The self-intersection is evident from Figure \ref{SelfIntOneGonZigzagCentralCircuit}.

(ii) If a central circuit is simple in a $(\{1,2,3\},6)$-sphere $G$,
then it splits $G$ into two domains $D_1$ and $D_2$. If one denotes
$n_{i,j}$ the number of faces of size $i$ into the domain $D_j$,
then one has obviously $2n_{1,j}+n_{2,j}=3$. So, if $n_{2,j}=0$,
then there is no solution. The proof for zigzags is the same. \qed

\begin{figure}
\begin{center}
\begin{minipage}[b]{4.9cm}
\centering
\epsfig{height=12mm, file=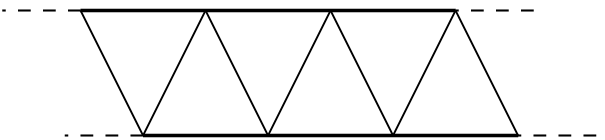}\par
The central circuit case
\end{minipage}
\begin{minipage}[b]{4.9cm}
\centering
\epsfig{height=12mm, file=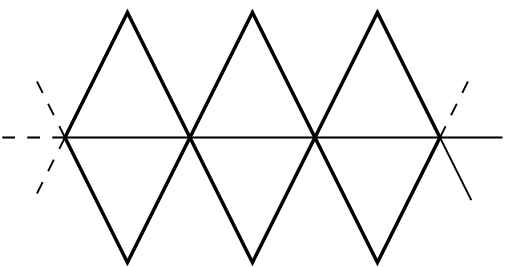}\par
The zigzag case
\end{minipage}

\end{center}
\caption{A $c$-railroad and a $z$-rairoadss bounded by two central circuits, respectively zigzags}
\label{ParallelZZ_CC}
\end{figure}

A $z$-, respectively 
{\em $c$-railroad} is the circuit of $3$-gons bounded
by two parallel zigzags, respectively central circuits.
See Figure \ref{ParallelZZ_CC} for an illustration.

\begin{figure}
\begin{center}
\begin{minipage}[b]{3.5cm}
\centering
\epsfig{height=23mm, file=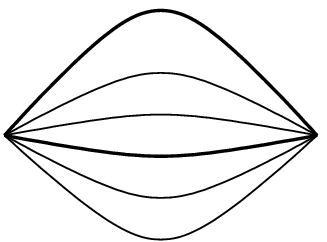}\par
A $c$-tight $(\{2,3\},6)$-sphere
\end{minipage}
\begin{minipage}[b]{3.5cm}
\centering
\epsfig{height=27mm, file=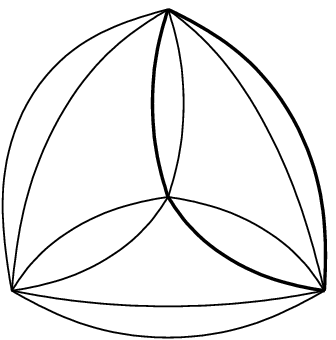}\par
A $c$-weakly tight $(\{2,3\},6)$-sphere
\end{minipage}
\begin{minipage}[b]{4.2cm}
\centering
\epsfig{height=27mm, file=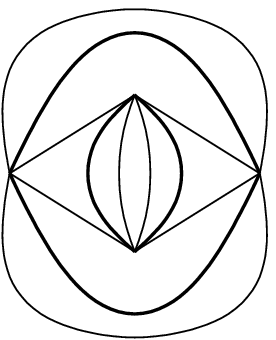}\par
A not $c$-weakly tight $(\{2,3\},6)$-sphere
\end{minipage}

\end{center}
\caption{Illustration of the notions of $c$-tightness}
\label{ExampleCentralCicruitIrrNotion}
\end{figure}

\begin{figure}
\begin{center}
\begin{minipage}[b]{3.5cm}
\centering
\epsfig{height=27mm, file=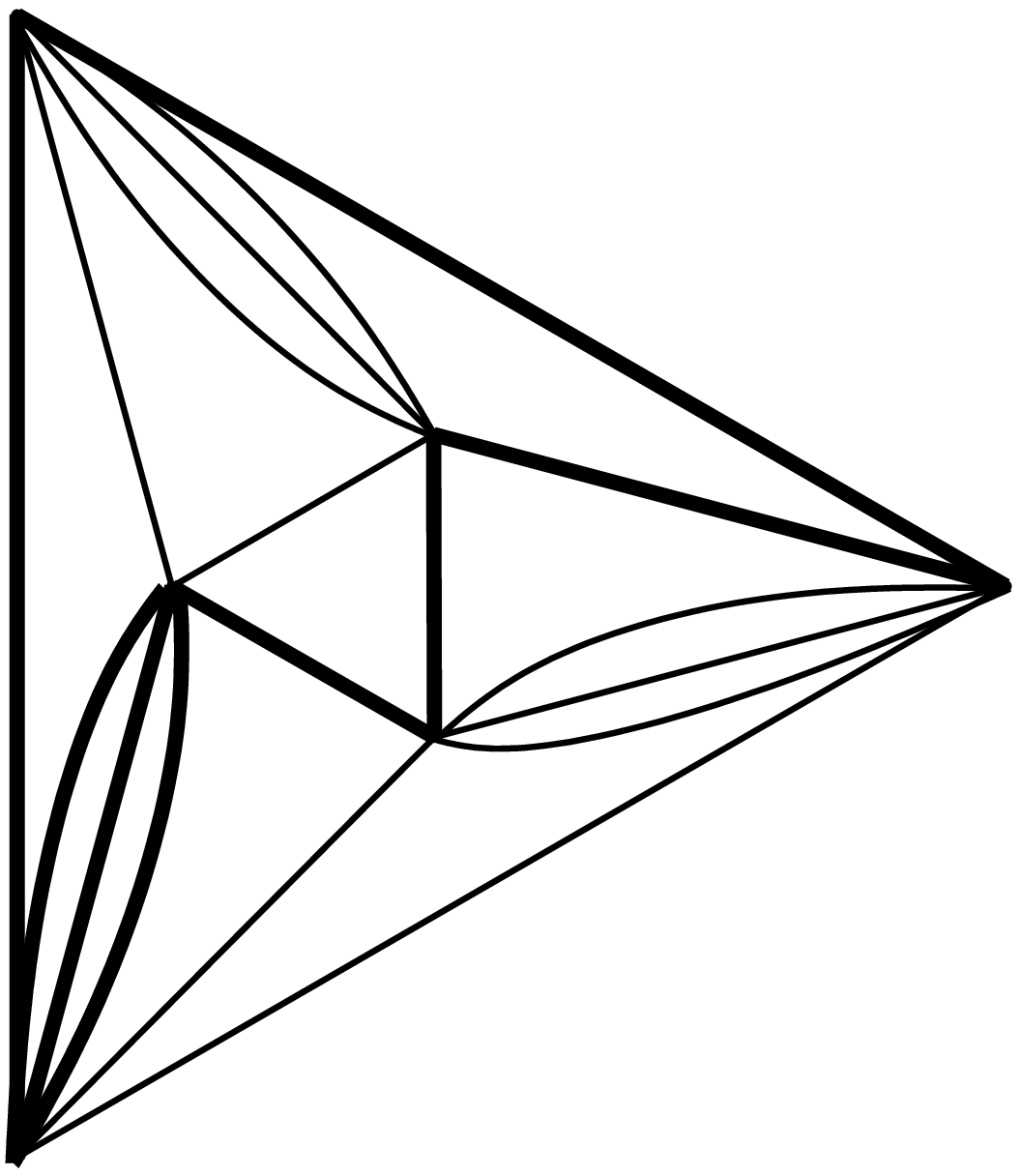}\par
A $z$-tight $(\{2,3\},6)$-sphere
\end{minipage}
\begin{minipage}[b]{3.5cm}
\centering
\epsfig{height=27mm, file=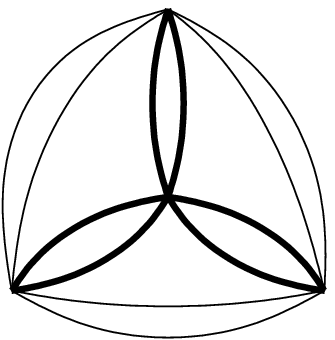}\par
A $z$-weakly tight $(\{2,3\},6)$-sphere
\end{minipage}
\begin{minipage}[b]{4.3cm}
\centering
\epsfig{height=27mm, file=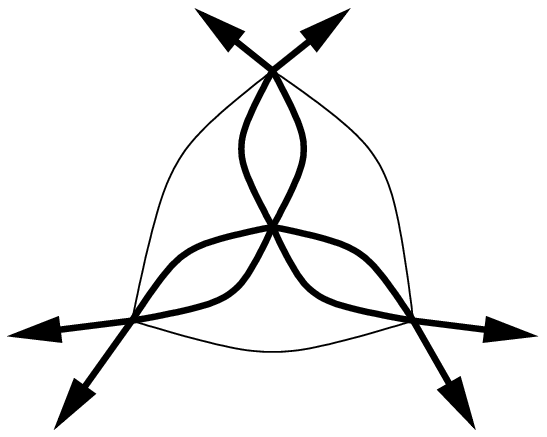}\par
A not $z$-weakly tight $(\{2,3\},6)$-sphere
\end{minipage}

\end{center}
\caption{Illustration of the notions of $z$-tightness}
\label{ExampleZigzagIrrNotion}
\end{figure}

A $(\{1,2,3\}, 6)$-sphere is called {\em $z$-tight}
if for any zigzag there is at least one $1$-gon or $2$-gon
on each of its side of the sphere.
It is called {\em $z$-weakly tight}
if for any zigzag there is no zigzag parallel to it.
We define the corresponding notions for central circuits.
See Figures \ref{ExampleZigzagIrrNotion} and 
\ref{ExampleCentralCicruitIrrNotion} for some illustration
of those notions.

The notion of tightness was introduced in \cite{octa}
for $4$-regular plane graphs (see also \cite{oct2}) and in \cite{zig1,zig2} 
for $3$-regular plane graphs. 
For $4$-regular plane graphs, central circuits were  used.
A central circuit is then called {\em reducible} if on one of its side
there are only $4$-gons. This sequence of $4$-gons can be completely
eliminated to get a reduced graph.
For a $3$-regular plane graph, a zigzag is called {\em reducible} if on
one side there is only $6$-gons. We can reduce the graph by eliminating
those $6$-gons only if the zigzag is simple. Moreover, there are
several possibilities for this reduced graph while in the $4$-regular
case, the reduction is uniquely defined.

For a $(\{1,2,3\},6)$-sphere $G$, let 
$s(G)=p_1(G) + 2p_2(G)$, where $p_i(G)$ is the number of $i$-gonal faces.
If $p\not=3$ a $p$-gon is called {\em incident to}
a zigzag or central circuit
if it share an edge with it.
It is called {\em weakly incident}
if it is not incident to it but still 
prevent the existence of a railroad.

\begin{theorem}\label{TrivialUpperBound}
For a $(\{1,2,3\},6)$-sphere $G$ we have:

(i) If $G$ is $z$-, respectively $c$-tight, then it
has at most $\frac{s(G)}{2}$ zigzags, respectively central circuits.

(ii) If $G$ is $z$-, respectively $c$-weakly tightness,
then it has at most $s(G)$ 
zigzags, respectively central circuits.

\end{theorem}
\proof Suppose that $G$ is $c$-tight, then for any central circuit
$C$ there is at least one face of of size different from $3$
on each sides.
Since the number of sides is $s(G)$ and there is two sides per central
circuit, this gives (i).
Also, the zigzag case is identical.

If $G$ is $c$-weakly tight, then a $p$-gon for $p\not=3$ is incident,
respectively weakly incident, to at most $p$ central
circuits.
Since it is weakly tight on each side of central circuits, there should
be at least one incident or weakly incident central circuit.
Thus the maximal
number of central circuits is $s(G)$ and the proof for zigzags is
identical. \qed

For $(\{1,2,3\}, 6)$-spheres with $i$ $1$-gons ($i=0,1,2,3$),
this gives the upper bounds of $(6, 4, 3, 1)$ for tightness
and $(12, 9, 6, 3)$ for weak tightness.

\begin{theorem}\label{TheoremZC_13spheres}
For a $(\{1,3\},6)$-sphere it holds:

(i) It cannot be $c$-, or $z$-tight.

(ii) Every central circuit correspond in a unique way to a
zigzag of doubled length.

(iii) If it is $c$- or $z$-weakly tight, then the number of 
central circuits, zigzags is $1$ or $3$.
\end{theorem}
\proof By Theorem \ref{GroupSpheres_C1_C2_C3}, all $(\{1,3\}, 6)$-spheres
have symmetry $C_3$, $C_{3v}$ or $C_{3h}$. Hence they have a $3$-fold axis
of rotation and hence the $1$-gons belong to a single orbit under the
group.
The faces of a $6$-regular plane graph are partitioned in two
classes, say $F_1$, $F_2$,
since their dual graph is bipartite.
Clearly, the $1$-gons are all in one partition class, say, $F_1$.
A $c$-, $z$-circuit has two sides, and the faces in those sides all
belong to the same partition class. Thus on one side
of any $zc$-circuit, there is only $3$-gons and so, (i) holds.

For a central circuit $C$, denote by $t_1$, \dots, $t_N$ the triangles
on the side of $F_2$.
Clearly, the set of edges of triangles $t_i$
not contained in $C$, define a zigzag and (ii) holds.

If $C$ is a central circuit in a $c$-weakly tight
$(\{1,3\},6)$-sphere $G$, then on the side of $F_1$
there is a $1$-gon and there is at most $3$ central circuits.
$2$ is excluded by the group action. \qed

\begin{table}
\caption{The maximal number of zigzags and central circuits for both notions of tightness and $4$ types of spheres. Bold numbers are definite answer, while intervals give the possible range}
\label{TableMaximalNrKnownResult}
\begin{center}
\begin{tabular}{|c||c|c||c|c||}
\hline
          & $z$-tig. & $z$-w. tig. & $c$-tig. & $c$-w. tig.\\
\hline
$(\{2,3\},6)$-spheres            & {\bf 6} & {\bf 9} & {\bf 6}  & $[8, 9]$\\
$(\{1,2,3\},6)$-spheres, $p_1=1$ & $[3,4]$ & $[5,7]$ & $[3,4]$  & $[5,7]$\\
$(\{1,2,3\},6)$-spheres, $p_1=2$ & {\bf 3} & {\bf 5} & {\bf 3} & $[4,5]$\\
$(\{1,3\},6)$-spheres &   {\bf 0}     & {\bf 3}  & {\bf 0} & {\bf 3}\\
\hline
\end{tabular}
\end{center}
\end{table}

\begin{figure}
\begin{center}
\begin{minipage}[b]{3.05cm}
\centering
\epsfig{height=23mm, file=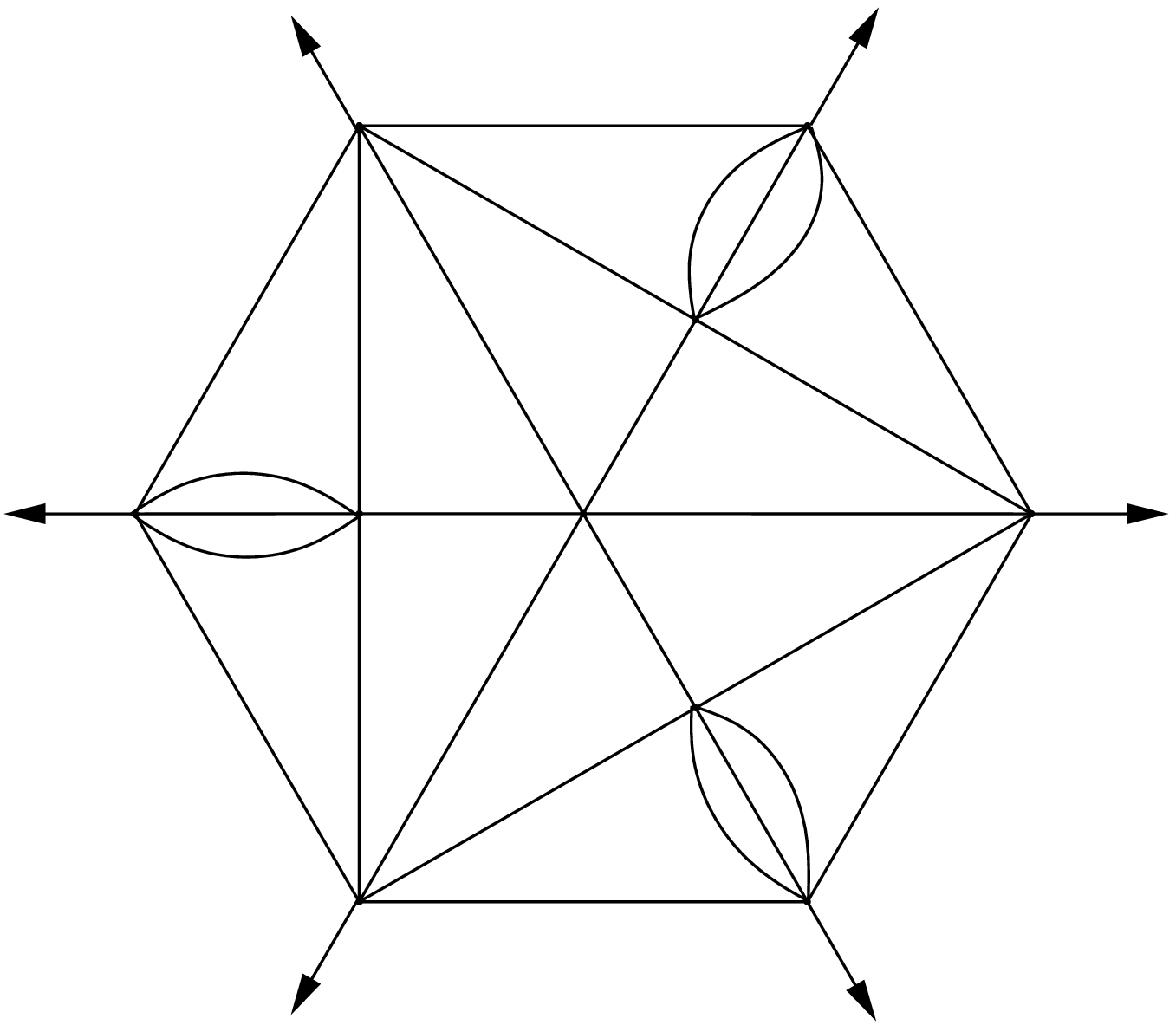}\par
$p_1=0$, $D_{3h}$, $n=11$ $c$-tight\par
$5^3, 6_{0,1}$
\end{minipage}
\begin{minipage}[b]{3.05cm}
\centering
\epsfig{height=23mm, file=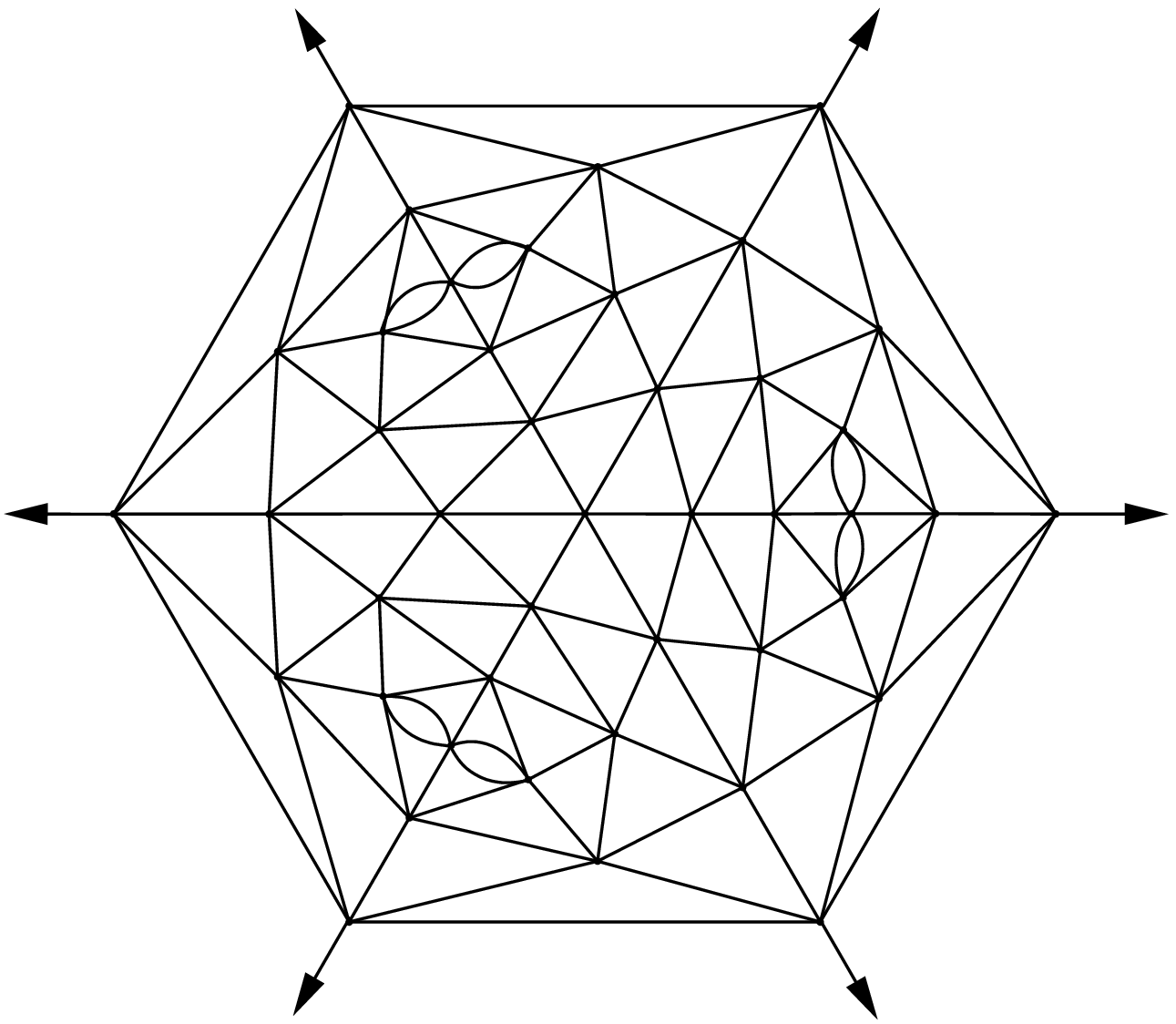}\par
$p_1=0$, $D_{3h}$, $n=44$ $c$-w. tig.\par
$10^3, 11_{0,1}^3, 22_{0,3}^3$
\end{minipage}
\begin{minipage}[b]{3.05cm}
\centering
\epsfig{height=27mm, file=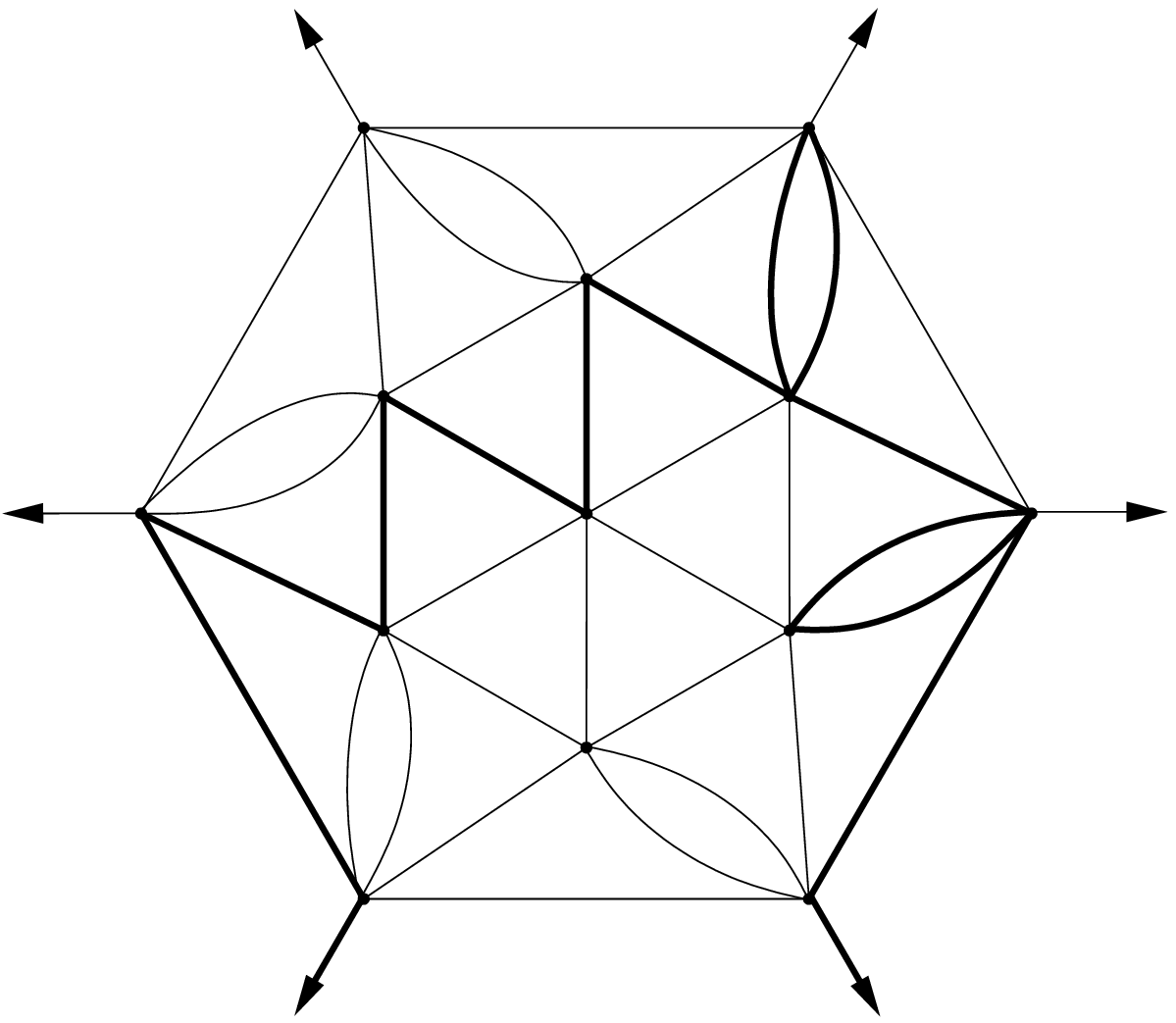}\par
$p_1=0$, $D_6$, $n=14$ $z$-tight\par
$14^6$
\end{minipage}
\begin{minipage}[b]{3.05cm}
\centering
\epsfig{height=23mm, file=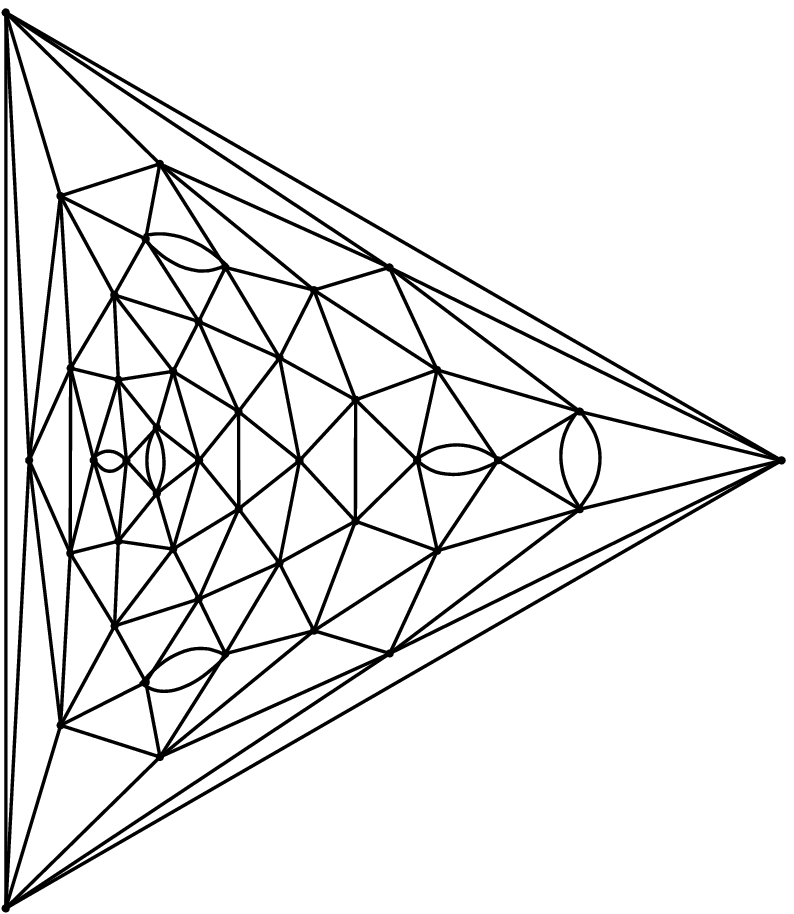}\par
$p_1=0$, $C_{2h}$, $n=44$ $z$-w. tig.\par
$24^4, 30_{0,1}^2, 54_{0,5}^2$
\end{minipage}
\begin{minipage}[b]{3.05cm}
\centering
\epsfig{height=23mm, file=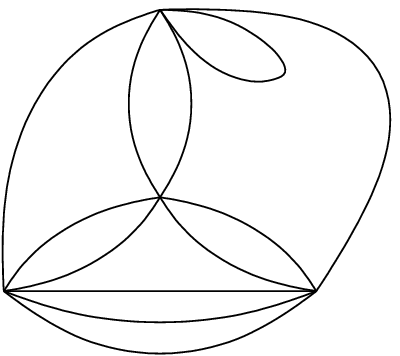}\par
$p_1=1$, $C_{s}$, $n=4$ $c$-tight\par
$3, 4_{0,1}, 5_{0,1}$
\end{minipage}
\begin{minipage}[b]{3.05cm}
\centering
\epsfig{height=23mm, file=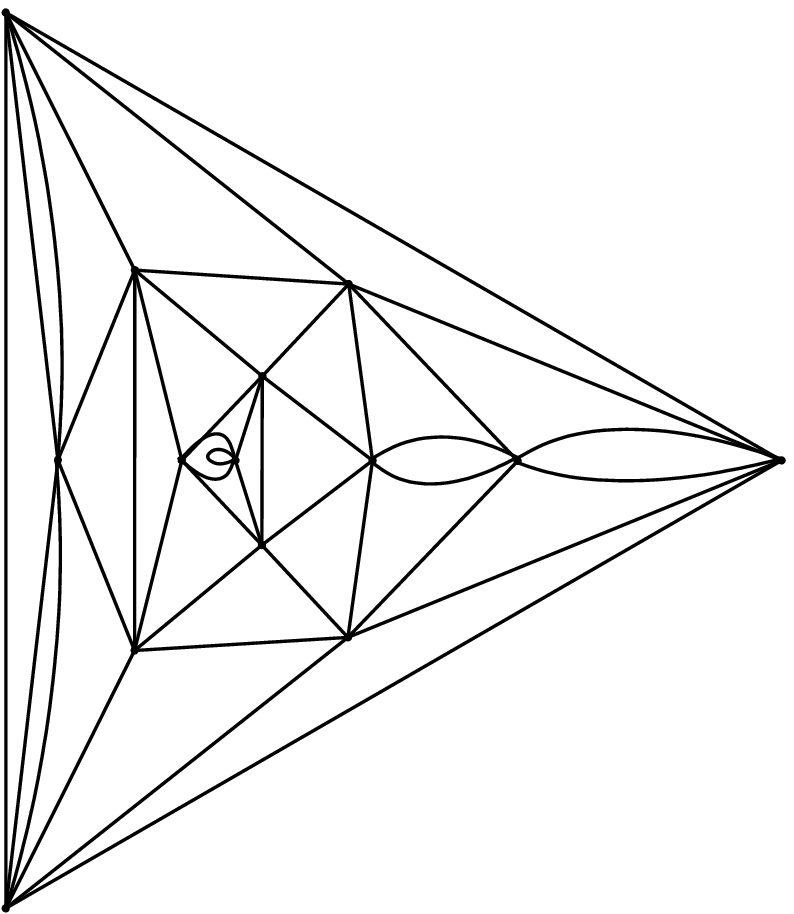}\par
$p_1=1$, $C_{s}$, $n=14$ $c$-w. tig.\par
$5, 11_{0,1}, 12_{0,3}, 7_{0,1}^2$
\end{minipage}
\begin{minipage}[b]{3.05cm}
\centering
\epsfig{height=23mm, file=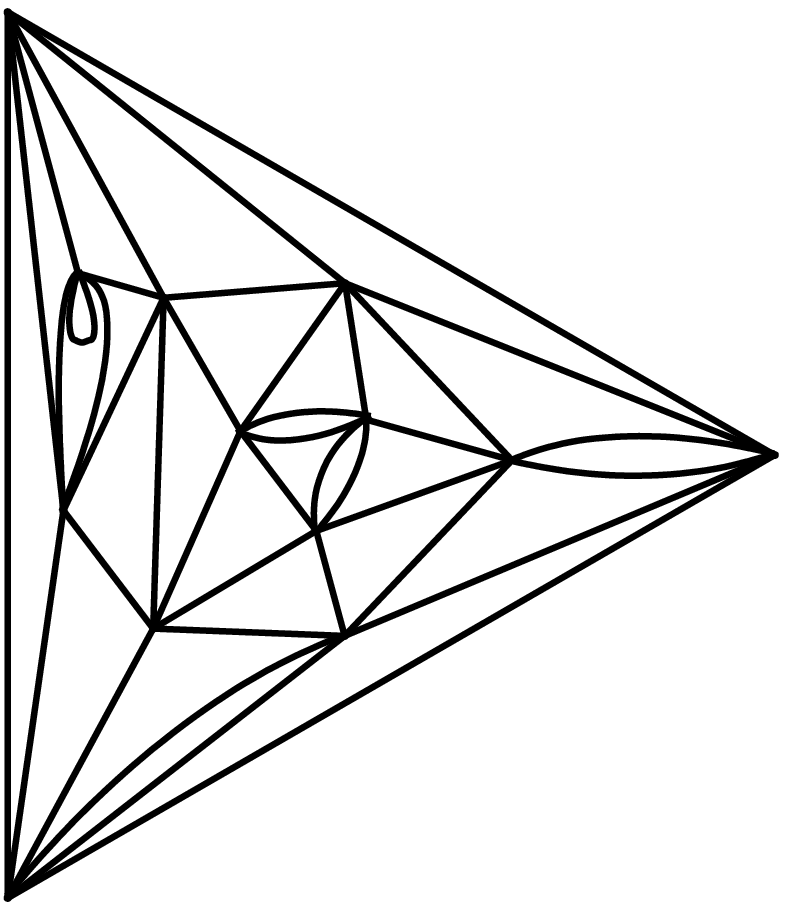}\par
$p_1=1$, $C_{1}$, $n=13$ $z$-tight\par
$16_{0,1}, 20_{0,1}, 43_{0,9}$
\end{minipage}
\begin{minipage}[b]{3.05cm}
\centering
\epsfig{height=23mm, file=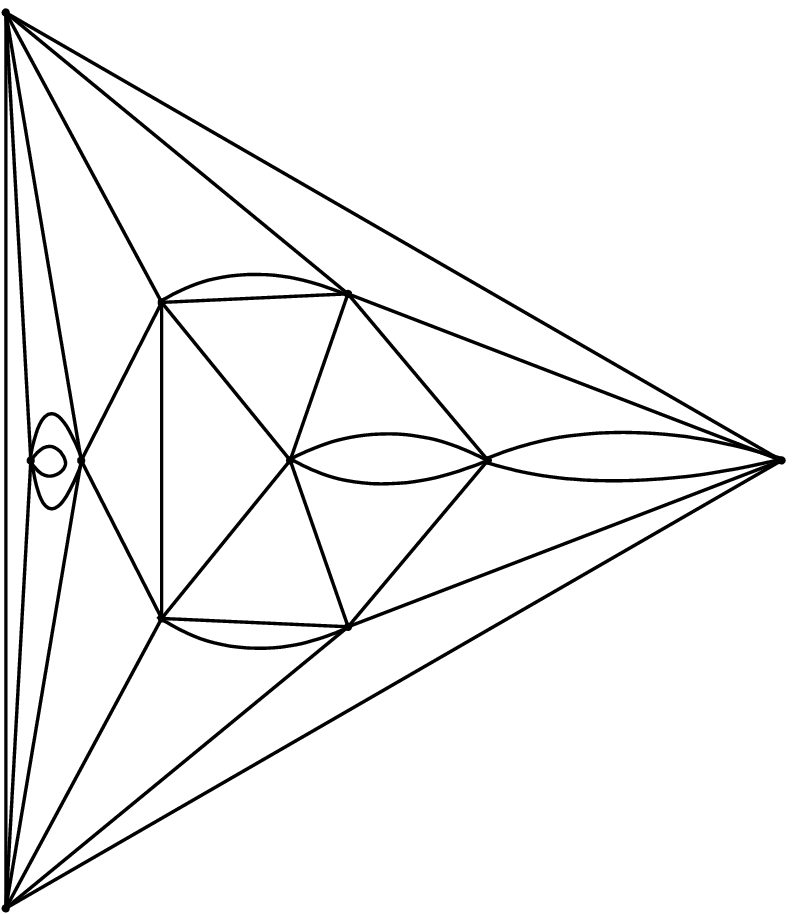}\par
$p_1=1$, $C_{s}$, $n=11$ $z$-w. tig.\par
$10, 12, 14_{0,1}^2, 16_{0,1}$
\end{minipage}
\begin{minipage}[b]{3.05cm}
\centering
\epsfig{height=23mm, file=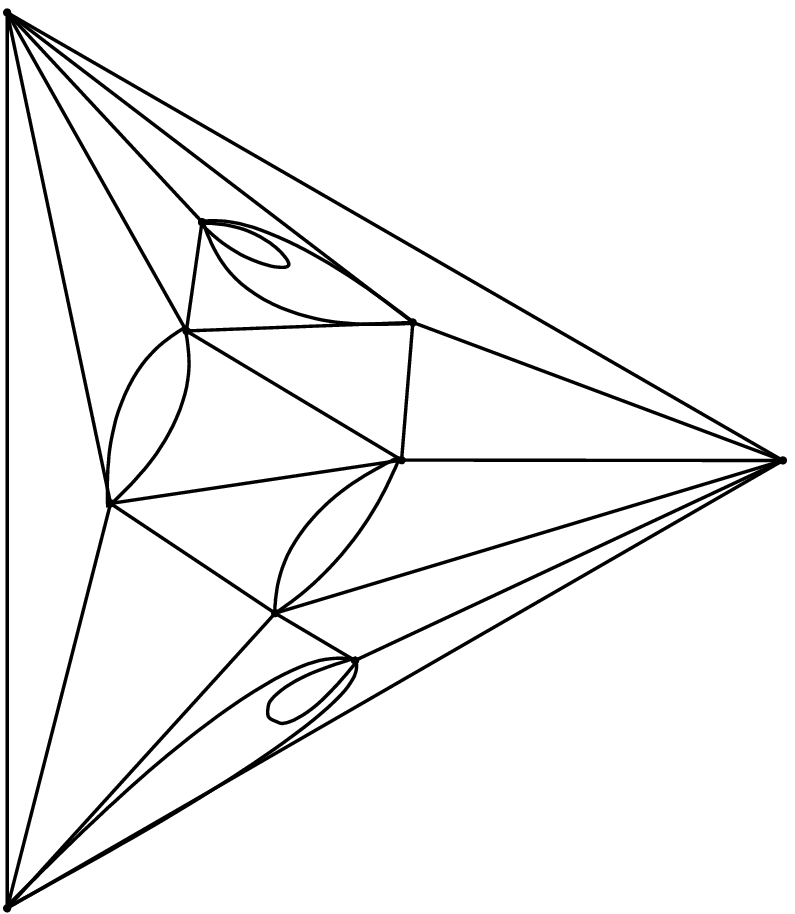}\par
$p_1=2$, $C_{2}$, $n=10$ $c$-tight\par
$8_{0,2}^2, 14_{0,6}$
\end{minipage}
\begin{minipage}[b]{3.05cm}
\centering
\epsfig{height=23mm, file=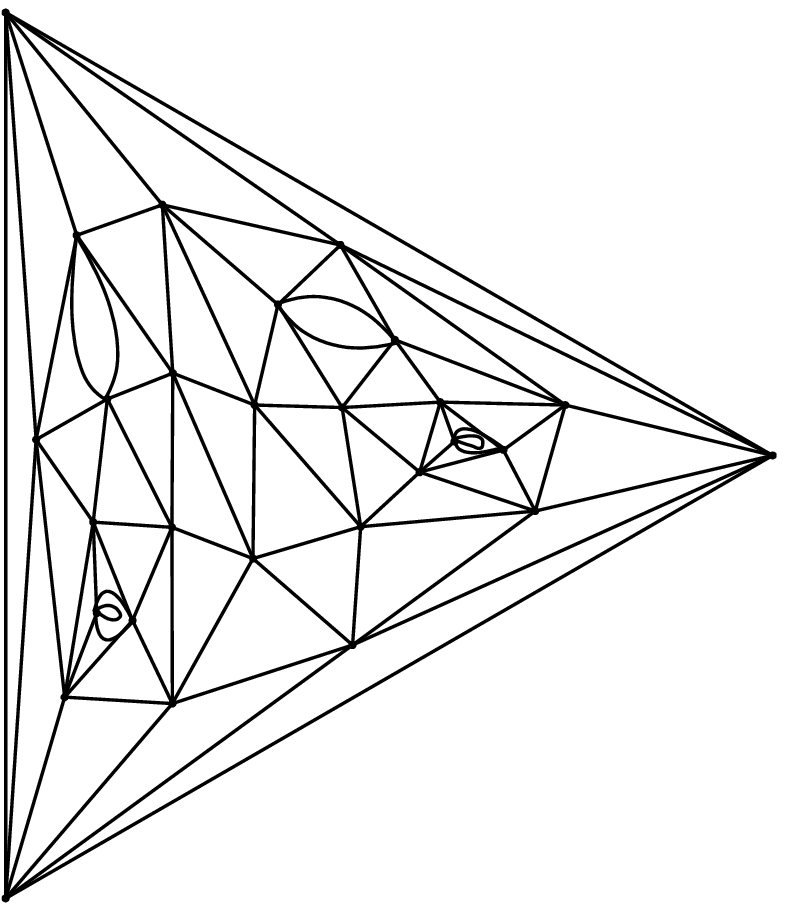}\par
$p_1=2$, $C_{2}$, $n=28$ $c$-w. tig.\par
$14_{0,1}, 14_{0,2}^3, 28_{0,9}$
\end{minipage}
\begin{minipage}[b]{3.05cm}
\centering
\epsfig{height=23mm, file=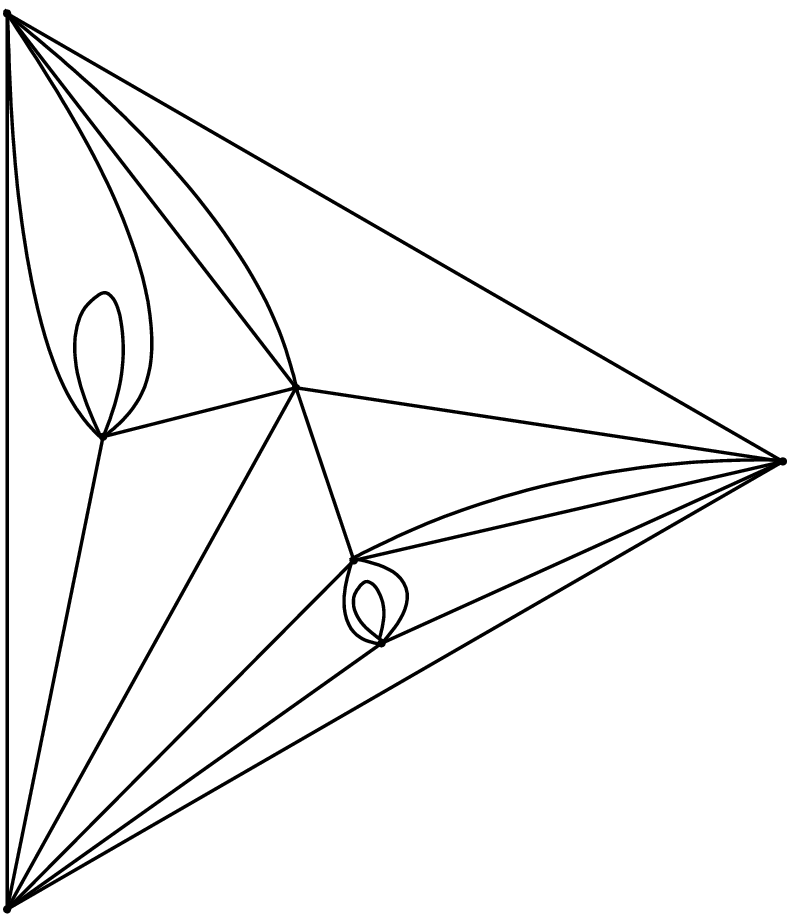}\par
$p_1=2$, $C_{2}$, $n=7$ $z$-tight\par
$14_{0,1}, 14_{0,2}^2$
\end{minipage}
\begin{minipage}[b]{3.05cm}
\centering
\epsfig{height=23mm, file=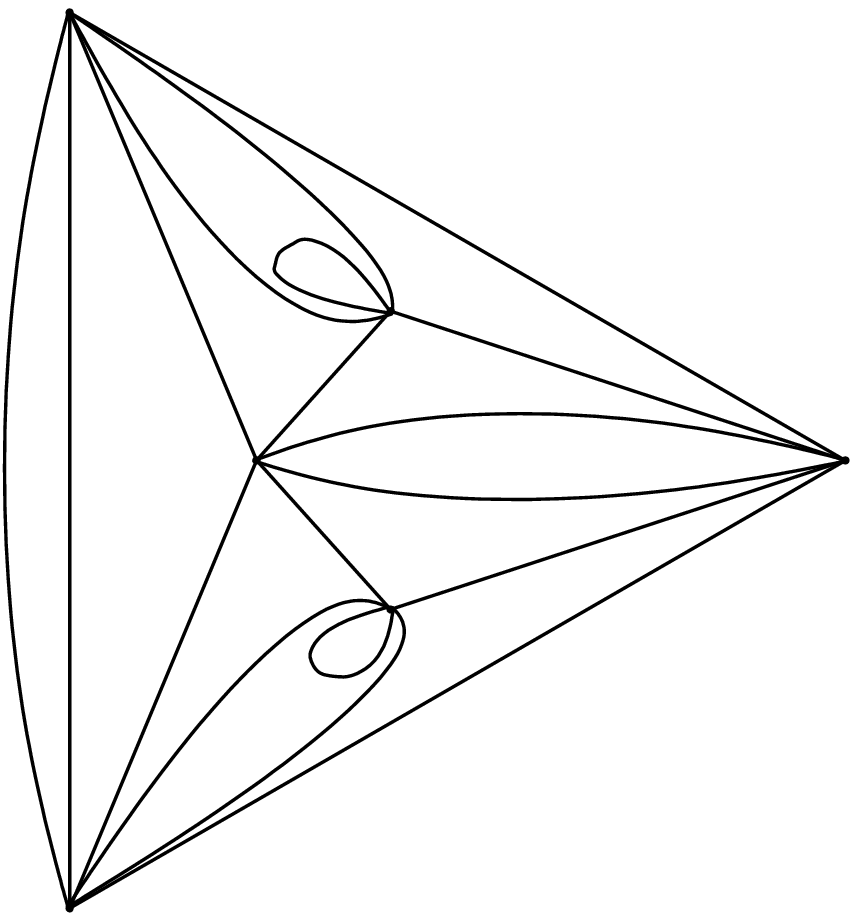}\par
$p_1=2$, $C_{2v}$, $n=6$ $z$-w. tig.\par
$6^2, 12_{0,2}^2$
\end{minipage}
\end{center}
\caption{The smallest weakly tight and tight $(\{1,2,3\},6)$-spheres with the maximal known number of zigzags and central circuits}
\label{WeaklyIrreducibleC0123_zz}
\end{figure}

\begin{theorem}\label{TheoremMaximalZC_irr_wirr}
Table \ref{TableMaximalNrKnownResult} for the maximal number of zigzags and central circuits and both notions of weak tightness and tightness hold.
\end{theorem}
\proof For $(\{1,3\},6)$-spheres, Theorem \ref{TheoremZC_13spheres} resolves
the question.
The existence of specific graphs in Figure \ref{WeaklyIrreducibleC0123_zz}
shows the lower bounds
that are indicated.
Theorem \ref{TrivialUpperBound} shows the required upper bounds for $z$-tightness and $c$-tightness.

For the notion of weak tightness, we have to provide something more.
Let $G$ be a $c$-weakly tight $(\{1,2,3\},6)$-sphere with central circuits
$C_1$, \dots, $C_l$.
The number of $1$-gons and $2$-gons is $p_1=i$, $p_2=6-2i$.
We obtain $2l$ sides since every central circuits has two sides.
A side $S$ is called {\em lonely} if it is incident or weakly incident
to only one $2$-gon.

\begin{figure}
\begin{center}
\begin{minipage}[b]{5.0cm}
\centering
\epsfig{height=27mm, file=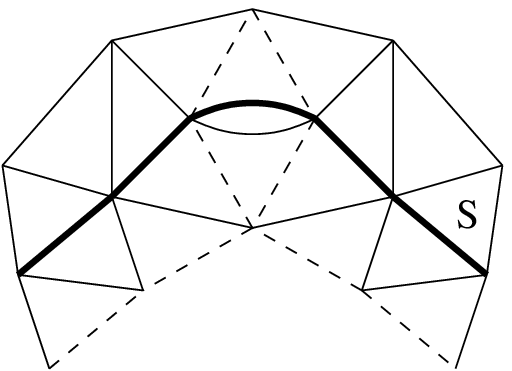}\par
(a) One lonely side $S$
\end{minipage}
\begin{minipage}[b]{5.0cm}
\centering
\epsfig{height=27mm, file=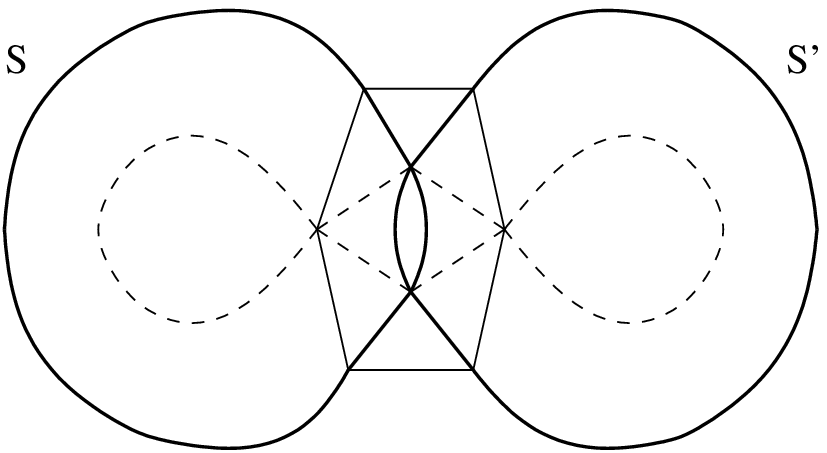}\par
(b) Two lonely sides $S$ and $S'$
\end{minipage}
\end{center}
\caption{Local structure around a side $S$ incident to a $2$-gon}
\label{ExampleLocalStructCC}
\end{figure}

If a side $S$ is incident to exactly one $2$-gon, then Figure \ref{ExampleLocalStructCC}.a shows
that there is a side of parallel central circuit that is weakly incident two
times to this $2$-gon. Moreover, if it is incident exactly two times
then there is another lonely side, see Figure \ref{ExampleLocalStructCC}.b.
A similar structure show up if a side is weakly incident to a $2$-gon.

Call $n_{1a}$ the number of lonely sides in the first case and $n_{1b}$
the number of lonely sides in the second case. Call $n_{1c}$ the number
of sides incident or weakly incident to exactly one $1$-gon.
Also let $n_2$ be the number of sides incident to exactly two $i$-gons
(identical or not). Let $n_3$ be the number of sides incident to at least
$3$ $i$-gons (identical or not).
Obviously, $l=\frac{1}{2}(n_{1a} + n_{1b} + n_{1c} +n_2 +n_3)$.

In case (a), a lonely side $S$ is incident to at least $3$ $i$-gons so
$n_{1a}\leq n_3$. Clearly, $n_{1c}\leq 2i$. Every $2$-gon can be incident
to $0$, $1$ or $2$ lonely sides so $n_{1a} + \frac{n_{1b}}{2} \leq 6-2i$.
By an enumeration of incidence we get
\begin{equation*}
n_{1a} +n_{1b} + n_{1c} + 2n_2 + 3 n_3 \leq s(G)=2i + 4 (6-2i)=24 - 6i.
\end{equation*}
Denote by ${\mathcal P}_i$ the $5$-dimensional polytope defined by
these inequalities and $n_{1a},\dots,n_3\geq 0$. We optimize the quantity
$l$ over ${\mathcal P}_i$ by using {\tt cdd} \cite{cdd}, which uses exact
arithmetic, and found the optimal value to be $9-2i$ for $i\leq 0\leq 3$.
The proof for zigzags is identical. \qed

In the rest of this section, we give a local Euler formula for central
circuits in order to enumerate the $(\{2,3\},6)$-spheres which are
$c$-weakly tight and with simple central circuits. The method for
zigzags is very similar and thus not indicated.

Let $G$ be a $(\{2,3\},6)$-sphere. Consider a patch $A$ in $G$, 
which is bounded by $t$ arcs (i.e., sections of central circuits)
belonging to central circuits (different or coinciding).

\begin{figure}
\begin{center}
\begin{minipage}[b]{8.0cm}
\centering
\epsfig{height=27mm, file=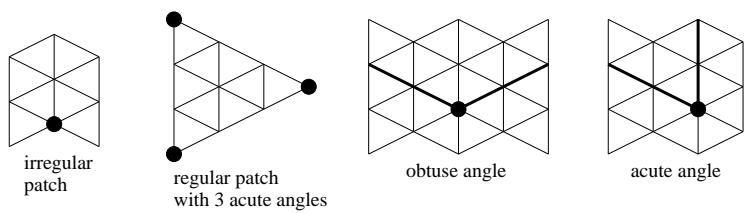}\par
\end{minipage}
\end{center}
\caption{Example of patches and their angles}
\label{fig:Regular_unregular}
\end{figure}

We admit also $0$-gonal patch $A$, i.e., just the interior of a simple 
central circuit.
Suppose that the patch $A$ is {\em regular}, i.e., the continuation 
of any of its bounding arcs (on the 
central circuit, to which it belongs) lies outside of the patch
(see Figure  \ref{fig:Regular_unregular}).
Let $p'_2(A)$ be the number of $2$-gonal faces in $A$.

There are two types of intersections of arcs on the boundary of a regular 
patch: either intersection in an edge of the boundary, or intersection in 
a vertex of the boundary. Let us call these types of intersections
{\em obtuse} and {\em acute}, respectively (see Figure
\ref{fig:Regular_unregular}); denote by $t_{ob}$ and $t_{ac}$ the 
respective number of obtuse and acute intersections. Clearly, 
$t_{ob}+t_{ac}=t$, where $t$ is the number of arcs forming the patch.
The following formula can easily be verified:
\begin{equation}\label{Local-Euler-FormulaCC}
6-t_{ob}-2t_{ac}=2p'_2(A).
\end{equation}

\begin{theorem}\label{TheoremIntersectionTwoSimpleZZ}
The intersection of every two simple central 
circuits, respectively zigzags,
of a $(\{2,3\},6)$-sphere, if non-empty, 
has one of the following forms (and so, its size is $2$, $4$ or $6$):
\begin{center}
\epsfig{height=25mm, file=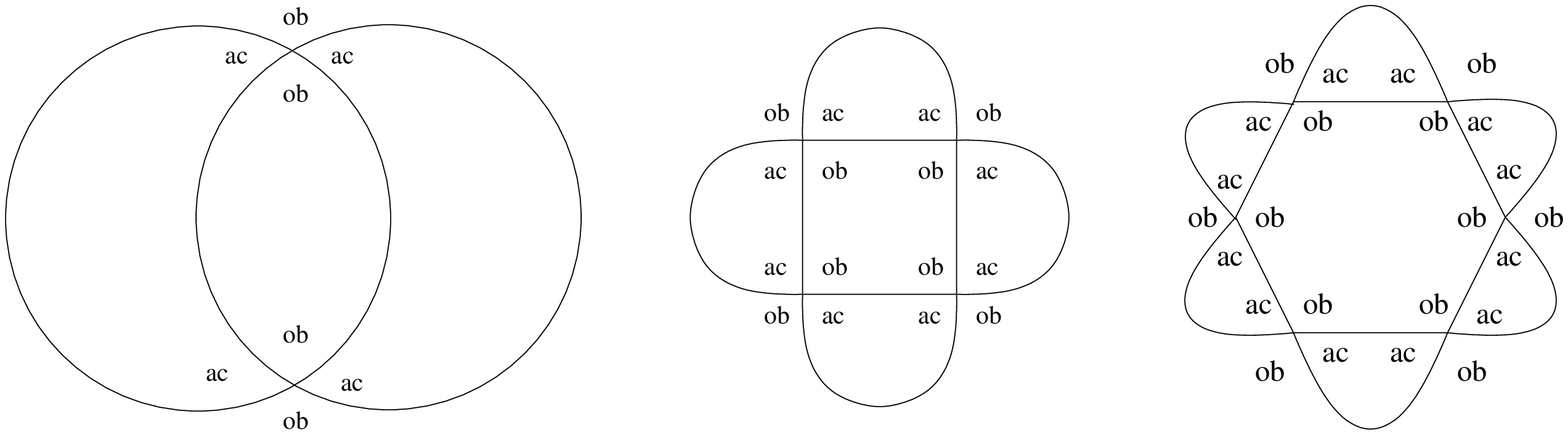}
\end{center}

\end{theorem}
\proof Let us consider the central circuit case, the zigzag case being
identical.
Define $H$ to be the graph, whose vertices are edges of
intersection between simple central circuits $C$ and $C'$, with
two vertices
being adjacent if they are linked by a path belonging to one of $C$, $C'$.
The graph $H$ is a plane $4$-regular graph and $C$, $C'$ 
define two
central circuits in $H$.
Since $C$ and $C'$ are simple,
the faces of $H$ are $t$-gons with even $t$.

Applying Formula \eqref{Local-Euler-FormulaCC} to a $t$-gonal face $F$ of $H$, we obtain that the number $p'_2(F)$ of $2$-gons in $F$ satisfies $6-t_{ob}-2t_{ac}=2p'_2(F)$.
So, the numbers $t_{ob}$ and $t_{ac}$ are even, since
$t=t_{ob}+t_{ac}$. Also, $6-t_{ob}-2t_{ac}\geq 0$.
So, $t\leq 6$.

We obtain the following five possibilities for the faces of $H$: 
$2$-gons with two acute angles, $2$-gons with two obtuse angles,
$4$-gons with four obtuse angles, $4$-gons with two acute and two
obtuse angles, $6$-gons with six obtuse angles.

\begin{figure}
\begin{center}
\epsfig{height=14mm, file=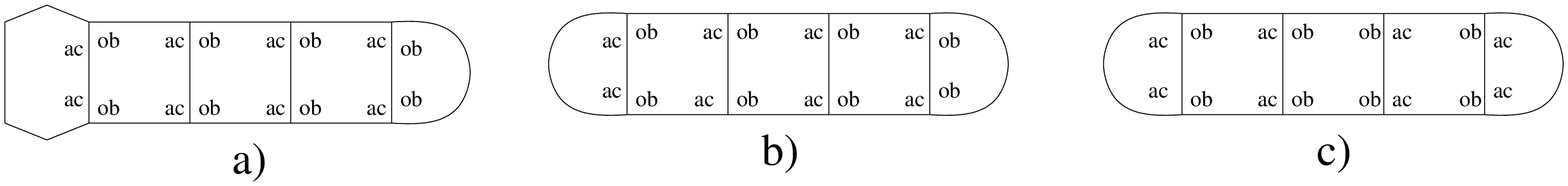}
\end{center}
\caption{Three cases for sequence of $4$-gons}
\label{fig:CaseSequence}
\end{figure}

Take an edge $e$ of a $6$-gon in $H$ and consider the sequence
(possibly, empty) of adjacent $4$-gons of $H$ emanating from this
edge. This sequence will stop at a $2$-gon or a $6$-gon; the
case-by-case analysis of angles yields that this sequence has to stop
at a $2$-gon (see Figure \ref{fig:CaseSequence}.a)).

Take an edge of a $2$-gon in $H$ and consider the same construction. 
If the angles are both obtuse, then the construction is identical and 
the sequence will terminate at a $2$-gon or a $6$-gon. If the angles are 
both acute, then cases b), c) of Figure \ref{fig:CaseSequence} are possible.

In the first case, all $4$-gons contain two obtuse angles and two acute
angles; so, the sequence of $4$-gons finishes with an edge of two obtuse
angles. In the second case, there is a $4$-gon, whose angles are all
obtuse; this $4$-gon is unique in the sequence and its position is
arbitrary. Every pair of opposite edges of a $4$-gon belongs to a 
sequence of $4$-gons considered above. So, all angles of a $4$-gon 
are the same, i.e., obtuse. This fact restricts 
the possibilities of intersections to the three cases of the theorem. \qed

The proof of this theorem is similar to the one of
Theorem 6.4 given in \cite{zig2}.

\begin{figure}
\begin{center}
\begin{minipage}[b]{3.2cm}
\centering
\epsfig{height=23mm, file=Bundle6irr.eps}\par
$D_{6h}$, $n=2$\par
$c=2^3$, tight
\end{minipage}
\begin{minipage}[b]{3.2cm}
\centering
\epsfig{height=27mm, file=Example23_6val_2wirr.eps}\par
$T_d$, $n=4$\par
$c=3^4$, weakly tight
\end{minipage}
\begin{minipage}[b]{3.2cm}
\centering
\epsfig{height=27mm, file=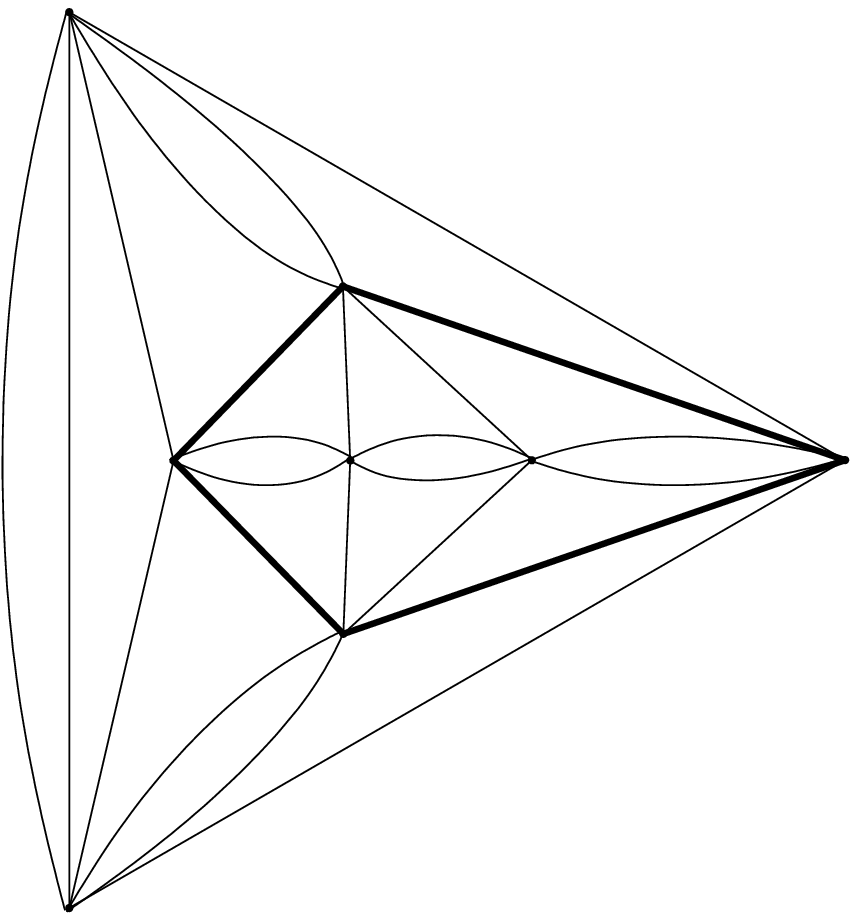}\par
$D_{2d}$, $n=8$\par
$c=4, 5^4$, weakly tight
\end{minipage}
\begin{minipage}[b]{3.2cm}
\centering
\epsfig{height=27mm, file=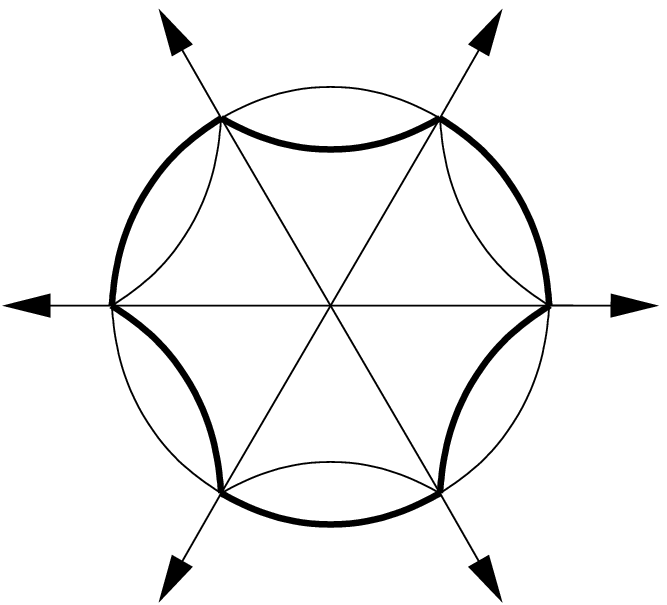}\par
$D_{6h}$, $n=8$\par
$c=4^3, 6^2$, weakly tight
\end{minipage}
\begin{minipage}[b]{3.2cm}
\centering
\epsfig{height=27mm, file=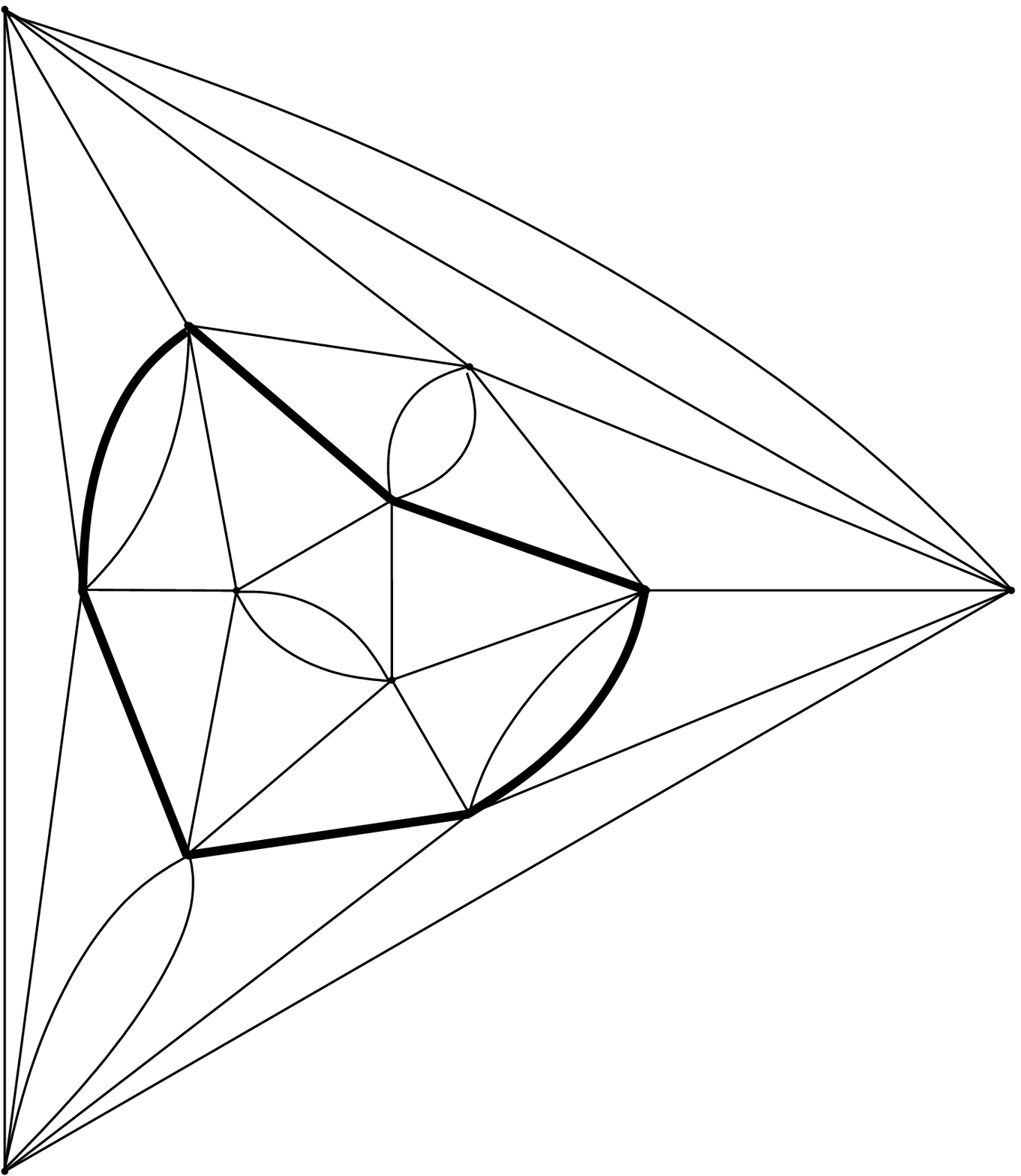}\par
$T_h$, $n=12$\par
$c=6^6$, weakly tight
\end{minipage}

\end{center}
\caption{The $c$-weakly tight $(\{2,3\}, 6)$-spheres with simple central circuits}
\label{WeaklyIrreducibleCentralCircuit}
\end{figure}

\begin{figure}
\begin{center}
\begin{minipage}[b]{3.2cm}
\centering
\epsfig{height=23mm, file=Bundle6.eps}\par
$D_{6h}$, $n=2$\par
$z=6^2$, tight
\end{minipage}
\begin{minipage}[b]{3.2cm}
\centering
\epsfig{height=27mm, file=Example23_6val_2_z_wirr.eps}\par
$T_d$, $n=4$\par
$z=6^4$, weakly tight
\end{minipage}
\begin{minipage}[b]{3.2cm}
\centering
\epsfig{height=27mm, file=SmallOcta_D3_irrSec.eps}\par
$D_3$, $n=6$\par
$z=8^3, 12$, tight
\end{minipage}
\begin{minipage}[b]{3.2cm}
\centering
\epsfig{height=27mm, file=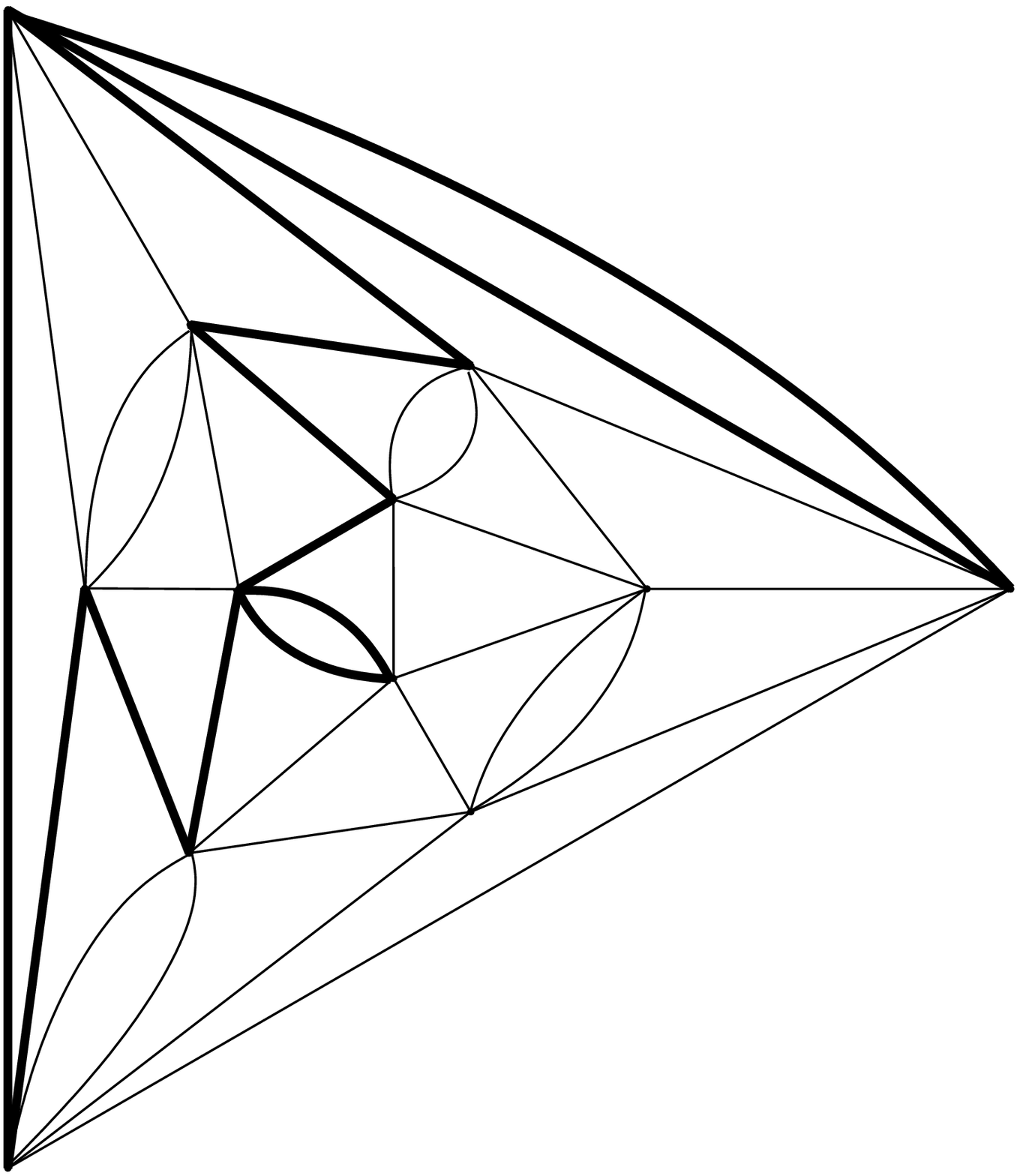}\par
$T_h$, $n=12$\par
$z=12^6$, weakly tight
\end{minipage}
\begin{minipage}[b]{3.2cm}
\centering
\epsfig{height=27mm, file=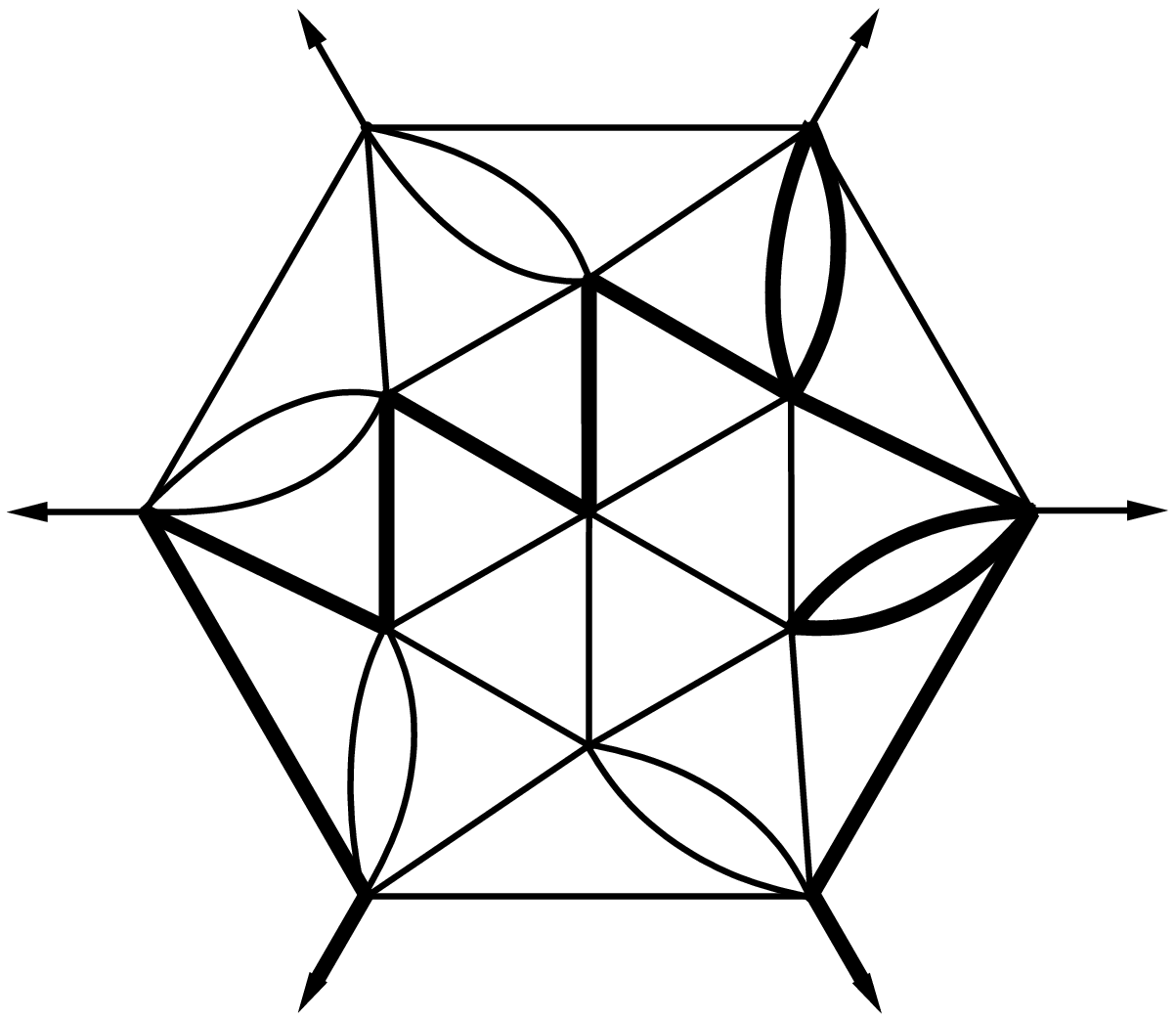}\par
$D_6$, $n=14$\par
$z=14^6$, tight
\end{minipage}

\end{center}
\caption{The $z$-weakly tight $(\{2,3\}, 6)$-spheres with simple zigzags}
\label{WeaklyIrreducibleZigzag}
\end{figure}

\begin{theorem}
The only weakly tight $(\{2,3\},6)$-spheres,
having only simple zigzags, respectively simple central circuits,
are the ones of Figure \ref{WeaklyIrreducibleZigzag} and
\ref{WeaklyIrreducibleCentralCircuit}.
\end{theorem}
\proof Let us consider first the central circuit case.
By Theorem \ref{TheoremIntersectionTwoSimpleZZ}, every two
simple central circuits intersect in at most six vertices.
If a $(\{2,3\}, 6)$-sphere has $t$ central circuits, this gives
an upper bound of $6 \frac{t(t-1)}{2}$ on the number of vertices
of intersection. Since any vertex can be the intersection of only
$3$ central circuits we get the upper bound of $t(t-1)$ on the number
of vertices.
If one uses the upper bound of Table \ref{TableMaximalNrKnownResult}
on $t$ for weakly tight $(\{2,3\},6)$-spheres, then one gets $t\leq 9$
and the upper bound $72$ on the number of vertices, which is too large
for the enumeration done in Table \ref{Enumeration123_6spheres_C0123}.
If one looks at the proof of Theorem \ref{TheoremMaximalZC_irr_wirr},
then one sees that a lonely side implies a self-intersection
of a parallel central circuit.
So, there is no lonely sides in
$(\{2,3\},6)$-spheres with only simple central circuits.
This gives the upper bound $t\leq 6$ on the number of central circuits
and then $30$ on the number of vertices.

For zigzags we have the upper bound $3t(t-1)$ on the number of edges
and this gives the same upper bound of $t(t-1)$ on the number of vertices.
The enumeration result shown in the Figures follow from the
determination results of Section \ref{GenerationMethod}. \qed

An interesting problem is to determine all $(\{2,3\},6)$-spheres
with simple zigzags and/or central circuits.

A $(\{1,2,3\},6)$-sphere is called {\em $z$-} or 
{\em $c$-knotted} if it has only one zigzag of central circuit.

\begin{conjecture}
(i) A $z$-, $c$-knotted $(\{1,2,3\},6)$-sphere, except Trifolium $C_{3v}$, 
has symmetry $C_1$, $C_2$, $C_3$, $D_2$ or $D_3$.

(ii) The $(\{2,3\},6)$-spheres with only simple central circuits have 
symmetry $T_d$, $T_h$, $D_{6h}$, $D_{3d}$, $D_{2d}$, $D_{2h}$, $D_3$,
$C_{2h}$ and $C_{3v}$.

(iii) A $(\{1,2,3\},6)$-sphere of symmetry 
$D_{6h}$, $T_h$, $T_d$ have only simple central circuits and zigzags.

(iv) The $(\{2,3\}, 6)$-spheres of symmetry $T_d$ have $v=4x^2$ vertices,
$c$-vector $(3x)^{4x}$ and $z$-vector $(6x)^{4x}$.

(v) The $(\{2,3\}, 6)$-spheres of symmetry $T_h$ have $12 x^2$ vertices,
$c$-vector $(6x)^{6x}$ and $z$-vector $(12x)^{6x}$.
\end{conjecture}

\begin{conjecture}
Let $f_i(v)$ denote the maximal number of central circuits 
in a $(\{1,2,3\}, 6)$-sphere with $i$ $1$-gons and $v$ vertices.
We conjecture:

(i) $f_2(v)=v+1$. It is realized exactly by the series (one for each $v\geq 1$)
having symmetry $C_{2h}$ and $c=1^v, (2v)_{0,v}$.

(ii) $f_1(v)=\frac{v-1}{2}+1$, $\frac{v-1}{2}+2$ for $v\equiv 3, 1\pmod 4$.
It is realized exactly by the series (one for each odd $v\ge 3$, all of symmetry $C_s$) 
with $c=2^{(v-1)/2}, (2v+1)_{0,v+2}$ if $v\equiv 3\pmod 4$ and
$2^{(v-1)/2}, v_{0,\frac{v-1}{4}}, (v+1)_{0,\frac{v+1}{4}}$ if $v\equiv 1\pmod 4$.

For even $v$,  
$f_1(v)=\lfloor \frac{v-1}{3}\rfloor +2$. It is realized for $v\ge 4$
by series of symmetry $C_s$ (two spheres for $v\equiv 2 \pmod 6$ and unique for other 
even $v$)
with $c=3^{\lfloor \frac{v-1}{3}\rfloor},  
(\frac{v}{2}+2+3\lfloor\frac{v}{18}\rfloor)_{0,2\lfloor\frac{v}{18}\rfloor+1}, 
(v+1+3z(v))_{0,4z(v)+1}$, 
where $z(v)=2\lfloor\frac{v}{18}\rfloor+1$ if $v\equiv 6,8,10 \pmod {18}$ and 
$z(v)=2\lfloor\frac{v+6}{18}\rfloor$ if $v$ is
other even number.

(iii) $f_0(v)=\frac{v}{2}+1$, $\frac{v}{2}+2$ for $v\equiv 0, 2\pmod 4$.
It is realized exactly by the series (one for each $v\ge 6$) of symmetry
$D_{2d}$ with $c=2^{\frac{v}{2}}, 2v_{0,v}$ if $v\equiv 0\pmod 4$
and of symmetry $D_{2h}$ with $c=2^{\frac{v}{2}}, (v_{0,\frac{v-2}{4}})^2$ if $v\equiv 
2\pmod 4$.

(iv) For add $v$, $f_0$ is $\lfloor \frac{v}{3} \rfloor+3$ if
$v\equiv 2,4,6,\pmod 9$ and $\lfloor\frac{v}{3}\rfloor +1$, otherwise.
Define $t_v$ by $\frac{v-t_v}{3}=\lfloor \frac{v}{3}\rfloor$.
$f_0(v)$ is realized by the series of symmetry $C_{3v}$ if
$v\equiv 1\pmod 3$ and $D_{3h}$, otherwise.
$c$-vector is $3^{\lfloor \frac{v}{3}\rfloor},(2\lfloor \frac{v}{3}\rfloor+t_v)^3_{0,\lfloor\frac{v-2t_v}{9}\rfloor}$ if $v\equiv 2,4,6\pmod 9$
and $3^{\lfloor \frac{v}{3}\rfloor}, (2v+t_v)_{0,v+2t_v}$, otherwise.

\end{conjecture}

\section{Acknowledgment}
Second author has been supported by the Croatian Ministry of Science, Education and Sport under contract 098-0982705-2707.


\begin{thebibliography}{99}

\bibitem[BDDH97]{CaGe}
G. Brinkmann, O. Delgado-Friedrichs, A. Dress and T. Harmuth, {\em CaGe -- a virtual environment for studying some special classes of large molecules}, MATCH {\bf 36} (1997) 233--237.

\bibitem[BrHaHe03]{Heidemeier}
G. Brinkmann, T. Harmuth and O. Heidemeier,
{\em The construction of cubic and quartic planar maps with prescribed face degrees},
Discrete Applied Mathematics, {\bf 128-2,3} (2003) 541--554.

\bibitem[Cox71]{Cox71}
H.S.M. Coxeter, {\em Virus macromolecules and geodesic domes}, in {\em A spectrum of mathematics}; Edited by J.C. Butcher, Oxford University Press/Auckland University Press: Oxford, U.K./Auckland, New-Zealand (1971) 98--107.

\bibitem[DDFo09]{zig1}
M. Deza, M. Dutour Sikiri\'c and P. Fowler, {\em The symmetries of cubic polyhedral graphs with face size no larger than $6$}, MATCH {\bf 61} (2009) 589-602

\bibitem[DeDu05]{zig2}
M. Deza and M. Dutour, {\em Zigzag structure of Simple Two-faced Polyhedra}, Combinatorics, Probability \& Computing {\bf 14} (2005) 31--57.

\bibitem[DeDu08]{book3}
M. Deza and M. Dutour Sikiri\'c, {\em Geometry of chemical graphs: polycycles and two-faced maps}, Cambridge University Press, Encyclopedia of mathematics and its applications, {\bf 119}, 2008.

\bibitem[DeDuSh03]{oct2}
M. Deza, M. Dutour and M. Shtogrin, {\em $4$-valent plane graphs with $2$-, $3$- and $4$-gonal faces}, ``Advances in Algebra and Related Topics'' (in memory of B.H.Neumann; Proceedings of ICM Satellite Conference on Algebra and Combinatorics, Hong Kong 2002), World Scientific Publ. Co. (2003) 73--97.

\bibitem[DeSt03]{octa}
M. Deza and M. Shtogrin,
{\em Octahedrites}, Symmetry, Special Issue 
{\bf 11-1,2,3,4} ``Polyhedra in Science and Art'' (2003) 27--64.

\bibitem[DHL02]{covcent}
M. Deza, T. Huang and K.-W. Lih, {\em Central circuit coverings of octahedrites and medial polyhedra}, Journal of Math. Research and Exposition {\bf 22-1} (2002) 49--66.

\bibitem[Dut04a]{pointgroup}
M. Dutour, {\em Point Groups}, \url{http://www.liga.ens.fr/~dutour/PointGroups/}, 2004

\bibitem[DuDe04]{goldberg}
M. Dutour and M. Deza, {\em Goldberg-Coxeter construction for $3$- or $4$-valent plane graphs}, Electronic Journal of Combinatorics {\bf 11} (2004) R20.

\bibitem[DuDe10]{selfdual}
M. Dutour Sikiri\'c and M. Deza, {\em 4-regular and self-dual analogs of fullerenes}, Mathematics and Topology of Fullerenes, Lecture notes in physics, 2010, to appear.

\bibitem[Fu]{cdd}
K. Fukuda, {\em The cdd program}, \url{http://www.ifor.math.ethz.ch/~fukuda/cdd_home/cdd.html}

\bibitem[Gold37]{Gold}
M. Goldberg, {\em A class of multisymmetric polyhedra}, Tohoku Mathematical Journal, {\bf 43} (1937) 104--108.

\bibitem[GrZa74]{GrZa}
B. Gr\"{u}nbaum and J. Zaks, {\em The existence of certain planar maps}, Discrete Mathematics, {\bf 10} (1974) 93-115.

\bibitem[Sah94]{sah94}
C.H. Sah, {\em A generalized leapfrog for fullerene structures}, Fullerenes Science and Technology {\bf 2-4} (1994) 445--458.

\bibitem[Sh75]{Sh}
H. Shank, {\em The theory of left-right paths}, in Combinatorial 
Mathematics III,
Proceedings of 3rd Australian Conference, St. Lucia 1974, Lecture Notes in
Mathematics 452, Springer-Verlag, New York (1975) 42--54.

\bibitem[Thur98]{T}
W.P. Thurston, {\em Shapes of polyhedra and triangulations of the sphere},
in Geometry and Topology Monographs 1, The Epstein Birthday Schrift,
(J. Rivin, C. Rourke and C. Series, eds.) Geom. Topol. Publ., Coventry,
1998, 511--549.

\bibitem[YHZQ10]{chinese4valent}
H. Yao, G. Hu, H. Zhang, W.Y. Qiu,
{\em The construction of $4$-regular polyhedra containing triangles, quadrilaterals and pentagons}, MATCH {\bf 64} (2010) 345--358.

\end{thebibliography}
\end{document}